\pdfoutput=1  %% COMMENT OUT IF "lualatex" IS USED

\documentclass[12pt,a4paper]{article}

\usepackage{amsfonts,amssymb,amsthm,amsmath,latexsym}
\usepackage{eqsection,indent}
\usepackage{subeqnarray, eepicemu, url,cite,bm}
\usepackage{multirow}  % See http://en.wikibooks.org/wiki/LaTeX/Tables
\usepackage{tensor}  % To help do 2F1, 2\phi1, etc
\usepackage{graphicx,graphbox}
\usepackage{tikz,xcolor}
\usepackage{float}  % http://tex.stackexchange.com/questions/2275/keeping-tables-figures-close-to-where-they-are-mentioned
\usepackage[scr=dutchcal,scrscaled=1.05]{mathalfa}

\usepackage{color}
\def\blue{\textcolor{blue}}
\def\red{\textcolor{red}}
\def\green{\textcolor{green}}

\date{April 13, 2023 \\[2mm] 
Revised April 2, 2024
}

\begin{document}

\title{\hspace*{-8mm}Continued fractions for cycle-alternating permutations}

\author{ \\
      \hspace*{-1cm}
      {\large Bishal Deb${}^{1,2}$ and Alan D.~Sokal${}^{1,3}$}
   \\[5mm]
     \hspace*{-1.3cm}
      \normalsize
           ${}^1$Department of Mathematics, University College London,
                    London WC1E 6BT, UK   \\[1mm]
     \hspace*{-1.15cm}
      \normalsize
           ${}^2$Sorbonne Universit\'e and Universit\'e Paris Cit\'e, CNRS,
         Laboratoire de Probabilit\'es, \\[-0.5mm]
     \hspace*{-7cm}
      \normalsize
         Statistique et Mod\'elisation, Paris, France \\[1mm]
     \hspace*{-2.9cm}
      \normalsize
	   ${}^3$Department of Physics, New York University,
                    New York, NY 10003, USA     \\[1mm]
       \\
     \hspace*{-0.5cm}
     {\tt bishal@gonitsora.com}, {\tt sokal@nyu.edu}  \\[1cm]
}

\maketitle
\thispagestyle{empty}   % Suppress page number on front page.

\begin{abstract}
A permutation is said to be cycle-alternating
if it has no cycle double rises, cycle double falls or fixed points;
thus each index $i$ is either a cycle valley ($\sigma^{-1}(i)>i<\sigma(i)$)
or a cycle peak ($\sigma^{-1}(i)<i>\sigma(i)$).
We find Stieltjes-type continued fractions for some multivariate
polynomials that enumerate cycle-alternating permutations
with respect to a large (sometimes infinite) number
of simultaneous statistics that measure cycle status, record status,
crossings and nestings along with the parity of the indices.
Our continued fractions are specializations of more general
continued fractions of Sokal and Zeng.
We then introduce alternating Laguerre digraphs,
which are generalization of cycle-alternating permutations,
and find exponential generating functions
for some polynomials enumerating them.
We interpret the Stieltjes--Rogers and Jacobi--Rogers matrices
associated to some of our continued fractions
in terms of alternating Laguerre digraphs.
\end{abstract}


\bigskip
\bigskip
\noindent
{\bf Mathematics Subject Classification (MSC 2020) codes:}
05A19 (Primary);
05A05, 05A15, 05A30, 11B68, 30B70 (Secondary).

%%\clearpage
\vspace*{1cm}

\newtheorem{theorem}{Theorem}[section]
\newtheorem{proposition}[theorem]{Proposition}
\newtheorem{lemma}[theorem]{Lemma}
\newtheorem{corollary}[theorem]{Corollary}
\newtheorem{definition}[theorem]{Definition}
\newtheorem{conjecture}[theorem]{Conjecture}
\newtheorem{question}[theorem]{Question}
\newtheorem{problem}[theorem]{Problem}
\newtheorem{openproblem}[theorem]{Open Problem}
\newtheorem{example}[theorem]{Example}

\renewcommand{\theenumi}{\alph{enumi}}
\renewcommand{\labelenumi}{(\theenumi)}
\def\eop{\hbox{\kern1pt\vrule height6pt width4pt
depth1pt\kern1pt}\medskip}
\def\prf{\par\noindent{\bf Proof.\enspace}\rm}
\def\rmk{\par\medskip\noindent{\bf Remark\enspace}\rm}

\newcommand{\textbfit}[1]{\textbf{\textit{#1}}}

\newcommand{\bigdash}{%
\smallskip\begin{center} \rule{5cm}{0.1mm} \end{center}\smallskip}

% DO NOT USE \par WITHIN A \footnote;  USE \safepar INSTEAD.
\newcommand{\safepar}{ {\protect\hfill\protect\break\hspace*{5mm}} }

\newcommand{\be}{\begin{equation}}
\newcommand{\ee}{\end{equation}}
\newcommand{\<}{\langle}
\renewcommand{\>}{\rangle}
\newcommand{\widebar}{\overline}
\def\reff#1{(\protect\ref{#1})}
\def\spose#1{\hbox to 0pt{#1\hss}}
\def\ltapprox{\mathrel{\spose{\lower 3pt\hbox{$\mathchar"218$}}
    \raise 2.0pt\hbox{$\mathchar"13C$}}}
\def\gtapprox{\mathrel{\spose{\lower 3pt\hbox{$\mathchar"218$}}
    \raise 2.0pt\hbox{$\mathchar"13E$}}}
\def\textprime{${}^\prime$}
\def\proof{\par\medskip\noindent{\sc Proof.\ }}
\def\firstproof{\par\medskip\noindent{\sc First Proof.\ }}
\def\secondproof{\par\medskip\noindent{\sc Second Proof.\ }}
\def\alternateproof{\par\medskip\noindent{\sc Alternate Proof.\ }}
\def\algebraicproof{\par\medskip\noindent{\sc Algebraic Proof.\ }}
\def\combinatorialproof{\par\medskip\noindent{\sc Combinatorial Proof.\ }}
\def\proofof#1{\bigskip\noindent{\sc Proof of #1.\ }}
\def\firstproofof#1{\bigskip\noindent{\sc First Proof of #1.\ }}
\def\secondproofof#1{\bigskip\noindent{\sc Second Proof of #1.\ }}
\def\thirdproofof#1{\bigskip\noindent{\sc Third Proof of #1.\ }}
\def\algebraicproofof#1{\bigskip\noindent{\sc Algebraic Proof of #1.\ }}
\def\combinatorialproofof#1{\bigskip\noindent{\sc Combinatorial Proof of #1.\ }}
\def\sketchofproof{\par\medskip\noindent{\sc Sketch of proof.\ }}
\renewcommand{\qed}{ $\square$ \bigskip}
\newcommand{\myendremark}{ $\blacksquare$ \bigskip}
\def\half{ {1 \over 2} }
\def\third{ {1 \over 3} }
\def\twothird{ {2 \over 3} }
\def\smfrac#1#2{{\textstyle{#1\over #2}}}
\def\smhalf{ {\smfrac{1}{2}} }
\newcommand{\real}{\mathop{\rm Re}\nolimits}
\renewcommand{\Re}{\mathop{\rm Re}\nolimits}
\newcommand{\imag}{\mathop{\rm Im}\nolimits}
\renewcommand{\Im}{\mathop{\rm Im}\nolimits}
\newcommand{\sgn}{\mathop{\rm sgn}\nolimits}
\newcommand{\tr}{\mathop{\rm tr}\nolimits}
\newcommand{\supp}{\mathop{\rm supp}\nolimits}
\newcommand{\disc}{\mathop{\rm disc}\nolimits}
\newcommand{\diag}{\mathop{\rm diag}\nolimits}
\newcommand{\tridiag}{\mathop{\rm tridiag}\nolimits}
\newcommand{\AZ}{\mathop{\rm AZ}\nolimits}
\newcommand{\NC}{\mathop{\rm NC}\nolimits}
\newcommand{\PF}{{\rm PF}}
\newcommand{\rk}{\mathop{\rm rk}\nolimits}
\newcommand{\perm}{\mathop{\rm perm}\nolimits}
\def\hboxscript#1{ {\hbox{\scriptsize\em #1}} }
\renewcommand{\emptyset}{\varnothing}
\newcommand{\eqdef}{\stackrel{\rm def}{=}}

\newcommand{\restrict}{\upharpoonright}

\newcommand{\compinv}{{\langle -1 \rangle}}   % Compositional inverse

\newcommand{\scra}{{\mathcal{A}}}
\newcommand{\scrb}{{\mathcal{B}}}
\newcommand{\scrc}{{\mathcal{C}}}
\newcommand{\bfscra}{{\bm{\mathcal{A}}}}
\newcommand{\bfscrb}{{\bm{\mathcal{B}}}}
\newcommand{\bfscrc}{{\bm{\mathcal{C}}}}
\newcommand{\bfscrap}{{\bm{\mathcal{A}'}}}
\newcommand{\bfscrbp}{{\bm{\mathcal{B}'}}}
\newcommand{\bfscrcp}{{\bm{\mathcal{C}'}}}
\newcommand{\bfscrapp}{{\bm{\mathcal{A}''}}}
\newcommand{\bfscrbpp}{{\bm{\mathcal{B}''}}}
\newcommand{\bfscrcpp}{{\bm{\mathcal{C}''}}}
\newcommand{\scrd}{{\mathcal{D}}}
\newcommand{\scre}{{\mathcal{E}}}
\newcommand{\scrf}{{\mathcal{F}}}
\newcommand{\scrg}{{\mathcal{G}}}
\newcommand{\scrgg}{{\mathscr{g}}}  % Lower-case calligraphic g (using Dutchcal)
\newcommand{\scrh}{{\mathcal{H}}}
\newcommand{\scri}{{\mathcal{I}}}
\newcommand{\scrj}{{\mathcal{J}}}
\newcommand{\scrk}{{\mathcal{K}}}
\newcommand{\scrl}{{\mathcal{L}}}
\newcommand{\scrm}{{\mathcal{M}}}
\newcommand{\scrn}{{\mathcal{N}}}
\newcommand{\scro}{{\mathcal{O}}}
\newcommand\scroo{
  \mathchoice
    {{\scriptstyle\mathcal{O}}}% \displaystyle
    {{\scriptstyle\mathcal{O}}}% \textstyle
    {{\scriptscriptstyle\mathcal{O}}}% \scriptstyle
    {\scalebox{0.6}{$\scriptscriptstyle\mathcal{O}$}}%\scriptscriptstyle
  }
%% Taken from https://tex.stackexchange.com/questions/191059/how-to-get-a-small-letter-version-of-mathcalo
\newcommand{\scrp}{{\mathcal{P}}}
\newcommand{\scrq}{{\mathcal{Q}}}
\newcommand{\scrr}{{\mathcal{R}}}
\newcommand{\scrs}{{\mathcal{S}}}
\newcommand{\scrss}{{\mathscr{s}}}  % Lower-case calligraphic s (using Dutchcal)
\newcommand{\scrt}{{\mathcal{T}}}
\newcommand{\scrv}{{\mathcal{V}}}
\newcommand{\scrw}{{\mathcal{W}}}
\newcommand{\scrz}{{\mathcal{Z}}}
\newcommand{\SP}{{\mathcal{SP}}}
\newcommand{\ST}{{\mathcal{ST}}}

\newcommand{\bfa}{{\mathbf{a}}}
\newcommand{\bfb}{{\mathbf{b}}}
\newcommand{\bfc}{{\mathbf{c}}}
\newcommand{\bfd}{{\mathbf{d}}}
\newcommand{\bfe}{{\mathbf{e}}}
\newcommand{\bfh}{{\mathbf{h}}}
\newcommand{\bfj}{{\mathbf{j}}}
\newcommand{\bfi}{{\mathbf{i}}}
\newcommand{\bfk}{{\mathbf{k}}}
\newcommand{\bfl}{{\mathbf{l}}}
\newcommand{\bfL}{{\mathbf{L}}}
\newcommand{\bfm}{{\mathbf{m}}}
\newcommand{\bfn}{{\mathbf{n}}}
\newcommand{\bfp}{{\mathbf{p}}}
\newcommand{\bfr}{{\mathbf{r}}}
\newcommand{\bfu}{{\mathbf{u}}}
\newcommand{\bfv}{{\mathbf{v}}}
\newcommand{\bfw}{{\mathbf{w}}}
\newcommand{\bfx}{{\mathbf{x}}}
\newcommand{\bfy}{{\mathbf{y}}}
\newcommand{\bfz}{{\mathbf{z}}}
\renewcommand{\k}{{\mathbf{k}}}
\newcommand{\n}{{\mathbf{n}}}
\newcommand{\vv}{{\mathbf{v}}}
\newcommand{\bv}{{\mathbf{v}}}
\newcommand{\w}{{\mathbf{w}}}
\newcommand{\x}{{\mathbf{x}}}
\newcommand{\y}{{\mathbf{y}}}
\newcommand{\cc}{{\mathbf{c}}}
\newcommand{\zero}{{\mathbf{0}}}
\newcommand{\one}{{\mathbf{1}}}
\newcommand{\bmm}{{\mathbf{m}}}

\newcommand{\ahat}{{\widehat{a}}}
\newcommand{\vhat}{{\widehat{v}}}
\newcommand{\yhat}{{\widehat{y}}}
\newcommand{\phat}{{\widehat{p}}}
\newcommand{\qhat}{{\widehat{q}}}
\newcommand{\Zhat}{{\widehat{Z}}}
\newcommand{\myhat}{{\widehat{\;}}}
\newcommand{\vtilde}{{\widetilde{v}}}
\newcommand{\ytilde}{{\widetilde{y}}}

\newcommand{\C}{{\mathbb C}}
\newcommand{\D}{{\mathbb D}}
\newcommand{\Z}{{\mathbb Z}}
\newcommand{\N}{{\mathbb N}}
\newcommand{\Q}{{\mathbb Q}}
\newcommand{\PP}{{\mathbb P}}
\newcommand{\R}{{\mathbb R}}
\newcommand{\RR}{{\mathbb R}}
\newcommand{\E}{{\mathbb E}}

\newcommand{\Sym}{{\mathfrak{S}}}
\newcommand{\SymB}{{\mathfrak{B}}}
\newcommand{\Cyc}{{\mathfrak{C}}}
\newcommand{\Altcyc}{{\scra\mathfrak{C}}}
\newcommand{\Alt}{{\mathrm{Alt}}}
\newcommand{\dperm}{{\mathfrak{D}}}
\newcommand{\dcycle}{{\mathfrak{DC}}}

\newcommand{\germanA}{{\mathfrak{A}}}
\newcommand{\germanB}{{\mathfrak{B}}}
\newcommand{\germanQ}{{\mathfrak{Q}}}
\newcommand{\germanh}{{\mathfrak{h}}}

\newcommand{\myle}{\preceq}
\newcommand{\myge}{\succeq}
\newcommand{\mygt}{\succ}

\newcommand{\B}{{\sf B}}
\newcommand{\OB}{B^{\rm ord}}
\newcommand{\OS}{{\sf OS}}
\newcommand{\OO}{{\sf O}}
\newcommand{\OSP}{{\sf OSP}}
\newcommand{\Eu}{{\sf Eu}}
\newcommand{\ERR}{{\sf ERR}}
\newcommand{\sfB}{{\sf B}}
\newcommand{\sfD}{{\sf D}}
\newcommand{\sfE}{{\sf E}}
\newcommand{\sfG}{{\sf G}}
\newcommand{\sfJ}{{\sf J}}
\newcommand{\sfL}{{\sf L}}
\newcommand{\sfLhat}{{\widehat{{\sf L}}}}
\newcommand{\sfLcheck}{{\widecheck{{\sf L}}}}
\newcommand{\sfLtilde}{{\widetilde{{\sf L}}}}
\newcommand{\sfP}{{\sf P}}
\newcommand{\sfQ}{{\sf Q}}
\newcommand{\sfS}{{\sf S}}
\newcommand{\sfT}{{\sf T}}
\newcommand{\sfW}{{\sf W}}
\newcommand{\sfMV}{{\sf MV}}
\newcommand{\AMV}{{\sf AMV}}
\newcommand{\BM}{{\sf BM}}
\newcommand{\emIB}{B^{\rm irr}}
\newcommand{\emIP}{P^{\rm irr}}
\newcommand{\emOB}{B^{\rm ord}}
\newcommand{\emCB}{B^{\rm cyc}}
\newcommand{\emSC}{P^{\rm cyc}}

\newcommand{\lev}{{\rm lev}}
\newcommand{\stat}{{\rm stat}}
\newcommand{\cyc}{{\rm cyc}}
\newcommand{\mysteryone}{{\rm mys1}}
\newcommand{\mysterytwo}{{\rm mys2}}
\newcommand{\Asc}{{\rm Asc}}
\newcommand{\asc}{{\rm asc}}
\newcommand{\Des}{{\rm Des}}
\newcommand{\des}{{\rm des}}
\newcommand{\Exc}{{\rm Exc}}

\newcommand{\EArec}{{\rm EArec}}
\newcommand{\earec}{{\rm earec}}
\newcommand{\recarec}{{\rm recarec}}
\newcommand{\erec}{{\rm erec}}
\newcommand{\nonrec}{{\rm nonrec}}
\newcommand{\nrar}{{\rm nrar}}
\newcommand{\Ereccval}{{\rm Ereccval}}
\newcommand{\ereccval}{{\rm ereccval}}
\newcommand{\ereccvalodd}{{\rm ereccvalodd}}
\newcommand{\ereccvaleven}{{\rm ereccvaleven}}
\newcommand{\ereccdrise}{{\rm ereccdrise}}
\newcommand{\Eareccpeak}{{\rm Eareccpeak}}
\newcommand{\eareccpeak}{{\rm eareccpeak}}
\newcommand{\Eareccpeakodd}{{\rm Eareccpeakodd}}
\newcommand{\Eareccpeakeven}{{\rm Eareccpeakeven}}
\newcommand{\eareccpeakodd}{{\rm eareccpeakodd}}
\newcommand{\eareccpeakeven}{{\rm eareccpeakeven}}
\newcommand{\eareccdfall}{{\rm eareccdfall}}
\newcommand{\eareccval}{{\rm eareccval}}
\newcommand{\ereccpeak}{{\rm ereccpeak}}
\newcommand{\Nrpeak}{{\rm Nrpeak}}
\newcommand{\rar}{{\rm rar}}
\newcommand{\evenrar}{{\rm evenrar}}
\newcommand{\oddrar}{{\rm oddrar}}
\newcommand{\nrcpeak}{{\rm nrcpeak}}
\newcommand{\nrcpeakodd}{{\rm nrcpeakodd}}
\newcommand{\nrcpeakeven}{{\rm nrcpeakeven}}
\newcommand{\nrcval}{{\rm nrcval}}
\newcommand{\nrcvalodd}{{\rm nrcvalodd}}
\newcommand{\nrcvaleven}{{\rm nrcvaleven}}
\newcommand{\nrcdrise}{{\rm nrcdrise}}
\newcommand{\nrcdfall}{{\rm nrcdfall}}
\newcommand{\nrfix}{{\rm nrfix}}
\newcommand{\Evenfix}{{\rm Evenfix}}
\newcommand{\evenfix}{{\rm evenfix}}
\newcommand{\Oddfix}{{\rm Oddfix}}
\newcommand{\oddfix}{{\rm oddfix}}
\newcommand{\evennrfix}{{\rm evennrfix}}
\newcommand{\oddnrfix}{{\rm oddnrfix}}
\newcommand{\Cpeak}{{\rm Cpeak}}
\newcommand{\cpeak}{{\rm cpeak}}
\newcommand{\Cval}{{\rm Cval}}
\newcommand{\cval}{{\rm cval}}
\newcommand{\cvaleven}{{\rm cvaleven}}
\newcommand{\cvalodd}{{\rm cvalodd}}
\newcommand{\Cdasc}{{\rm Cdasc}}
\newcommand{\cdasc}{{\rm cdasc}}
\newcommand{\Cddes}{{\rm Cddes}}
\newcommand{\cddes}{{\rm cddes}}
\newcommand{\Cdrise}{{\rm Cdrise}}
\newcommand{\cdrise}{{\rm cdrise}}
\newcommand{\Cdfall}{{\rm Cdfall}}
\newcommand{\cdfall}{{\rm cdfall}}

\newcommand{\maxpeak}{{\rm maxpeak}}
\newcommand{\nmaxpeak}{{\rm nmaxpeak}}
\newcommand{\minval}{{\rm minval}}
\newcommand{\nminval}{{\rm nminval}}

\newcommand{\exc}{{\rm exc}}
\newcommand{\excee}{{\rm excee}}
\newcommand{\exceo}{{\rm exceo}}
\newcommand{\excoe}{{\rm excoe}}
\newcommand{\excoo}{{\rm excoo}}
\newcommand{\erecexcoe}{{\rm erecexcoe}}
\newcommand{\erecexcoo}{{\rm erecexcoo}}
\newcommand{\nrexcoe}{{\rm nrexcoe}}
\newcommand{\nrexcoo}{{\rm nrexcoo}}
\newcommand{\aexc}{{\rm aexc}}
\newcommand{\aexcee}{{\rm aexcee}}
\newcommand{\aexceo}{{\rm aexceo}}
\newcommand{\aexcoe}{{\rm aexcoe}}
\newcommand{\aexcoo}{{\rm aexcoo}}
\newcommand{\earecaexcee}{{\rm earecaexcee}}
\newcommand{\earecaexceo}{{\rm earecaexceo}}
\newcommand{\nraexcee}{{\rm nraexcee}}
\newcommand{\nraexceo}{{\rm nraexceo}}
\newcommand{\Fix}{{\rm Fix}}
\newcommand{\fix}{{\rm fix}}
\newcommand{\fixe}{{\rm fixe}}
\newcommand{\fixo}{{\rm fixo}}
\newcommand{\rare}{{\rm rare}}
\newcommand{\raro}{{\rm raro}}
\newcommand{\nrfixe}{{\rm nrfixe}}
\newcommand{\nrfixo}{{\rm nrfixo}}
\newcommand{\xee}{x_{\rm ee}}
\newcommand{\xeo}{x_{\rm eo}}
\newcommand{\uee}{u_{\rm ee}}
\newcommand{\ueo}{u_{\rm eo}}
\newcommand{\yoo}{y_{\rm oo}}
\newcommand{\yoe}{y_{\rm oe}}
\newcommand{\voo}{v_{\rm oo}}
\newcommand{\voe}{v_{\rm oe}}
\newcommand{\zo}{z_{\rm o}}
\newcommand{\ze}{z_{\rm e}}
\newcommand{\wo}{w_{\rm o}}
\newcommand{\we}{w_{\rm e}}
\newcommand{\so}{s_{\rm o}}
\newcommand{\se}{s_{\rm e}}

\newcommand{\xo}{x_{\rm o}}
\newcommand{\xe}{x_{\rm e}}
\newcommand{\yo}{y_{\rm o}}
\newcommand{\ye}{y_{\rm e}}
\newcommand{\uo}{u_{\rm o}}
\newcommand{\ue}{u_{\rm e}}
\newcommand{\vo}{v_{\rm o}}
\newcommand{\ve}{v_{\rm e}}

\newcommand{\Wex}{{\rm Wex}}
\newcommand{\wex}{{\rm wex}}
\newcommand{\lrmax}{{\rm lrmax}}
\newcommand{\rlmax}{{\rm rlmax}}
\newcommand{\Rec}{{\rm Rec}}
\newcommand{\rec}{{\rm rec}}
\newcommand{\Arec}{{\rm Arec}}
\newcommand{\Arecpeak}{{\rm Arecpeak}}
\newcommand{\arec}{{\rm arec}}
\newcommand{\Even}{{\rm Even}}
\newcommand{\Odd}{{\rm Odd}}
\newcommand{\ERec}{{\rm ERec}}
\newcommand{\Val}{{\rm Val}}
\newcommand{\Peak}{{\rm Peak}}
\newcommand{\dasc}{{\rm dasc}}
\newcommand{\ddes}{{\rm ddes}}
\newcommand{\inv}{{\rm inv}}
\newcommand{\maj}{{\rm maj}}
\newcommand{\rs}{{\rm rs}}
\newcommand{\cross}{{\rm cr}}
\newcommand{\crosshat}{{\widehat{\rm cr}}}
\newcommand{\nest}{{\rm ne}}
\newcommand{\ucross}{{\rm ucross}}
\newcommand{\ucrosscval}{{\rm ucrosscval}}
\newcommand{\ucrosscpeak}{{\rm ucrosscpeak}}
\newcommand{\ucrosscdrise}{{\rm ucrosscdrise}}
\newcommand{\lcross}{{\rm lcross}}
\newcommand{\lcrosscpeak}{{\rm lcrosscpeak}}
\newcommand{\lcrosscval}{{\rm lcrosscval}}
\newcommand{\lcrosscdfall}{{\rm lcrosscdfall}}
\newcommand{\unest}{{\rm unest}}
\newcommand{\unestcval}{{\rm unestcval}}
\newcommand{\unestcpeak}{{\rm unestcpeak}}
\newcommand{\unestcdrise}{{\rm unestcdrise}}
\newcommand{\lnest}{{\rm lnest}}
\newcommand{\lnestcpeak}{{\rm lnestcpeak}}
\newcommand{\lnestcval}{{\rm lnestcval}}
\newcommand{\lnestcdfall}{{\rm lnestcdfall}}
\newcommand{\ulev}{{\rm ulev}}
\newcommand{\llev}{{\rm llev}}
\newcommand{\ujoin}{{\rm ujoin}}
\newcommand{\ljoin}{{\rm ljoin}}
\newcommand{\psnest}{{\rm psnest}}
\newcommand{\upsnest}{{\rm upsnest}}
\newcommand{\lpsnest}{{\rm lpsnest}}
\newcommand{\epsnest}{{\rm epsnest}}
\newcommand{\opsnest}{{\rm opsnest}}
\newcommand{\rodd}{{\rm rodd}}
\newcommand{\reven}{{\rm reven}}
\newcommand{\lodd}{{\rm lodd}}
\newcommand{\leven}{{\rm leven}}
\newcommand{\sg}{{\rm sg}}
\newcommand{\bl}{{\rm bl}}
\newcommand{\tran}{{\rm tr}}
\newcommand{\area}{{\rm area}}
\newcommand{\ret}{{\rm ret}}
\newcommand{\peaks}{{\rm peaks}}
\newcommand{\hl}{{\rm hl}}
\newcommand{\sll}{{\rm sl}}
\newcommand{\negg}{{\rm neg}}
\newcommand{\imp}{{\rm imp}}
\newcommand{\osg}{{\rm osg}}
\newcommand{\ons}{{\rm ons}}
\newcommand{\isg}{{\rm isg}}
\newcommand{\ins}{{\rm ins}}
\newcommand{\LL}{{\rm LL}}
\newcommand{\height}{{\rm ht}}
\newcommand{\as}{{\rm as}}

\newcommand{\ba}{{\bm{a}}}
\newcommand{\bahat}{{\widehat{\bm{a}}}}
\newcommand{\bb}{{\bm{b}}}
\newcommand{\bc}{{\bm{c}}}
\newcommand{\bchat}{{\widehat{\bm{c}}}}
\newcommand{\bd}{{\bm{d}}}
\newcommand{\bee}{{\bm{e}}}
\newcommand{\beh}{{\bm{eh}}}
\newcommand{\bff}{{\bm{f}}}
\newcommand{\bg}{{\bm{g}}}
\newcommand{\bh}{{\bm{h}}}
\newcommand{\bll}{{\bm{\ell}}}
\newcommand{\bp}{{\bm{p}}}
\newcommand{\br}{{\bm{r}}}
\newcommand{\bs}{{\bm{s}}}
\newcommand{\bu}{{\bm{u}}}
\newcommand{\bw}{{\bm{w}}}
\newcommand{\bx}{{\bm{x}}}
\newcommand{\by}{{\bm{y}}}
\newcommand{\bz}{{\bm{z}}}
\newcommand{\bA}{{\bm{A}}}
\newcommand{\bB}{{\bm{B}}}
\newcommand{\bC}{{\bm{C}}}
\newcommand{\bE}{{\bm{E}}}
\newcommand{\bF}{{\bm{F}}}
\newcommand{\bG}{{\bm{G}}}
\newcommand{\bH}{{\bm{H}}}
\newcommand{\bI}{{\bm{I}}}
\newcommand{\bJ}{{\bm{J}}}
\newcommand{\bM}{{\bm{M}}}
\newcommand{\bN}{{\bm{N}}}
\newcommand{\bP}{{\bm{P}}}
\newcommand{\bQ}{{\bm{Q}}}
\newcommand{\bR}{{\bm{R}}}
\newcommand{\bS}{{\bm{S}}}
\newcommand{\bT}{{\bm{T}}}
\newcommand{\bW}{{\bm{W}}}
\newcommand{\bX}{{\bm{X}}}
\newcommand{\bY}{{\bm{Y}}}
\newcommand{\bIB}{{\bm{B}^{\rm irr}}}
\newcommand{\bOB}{{\bm{B}^{\rm ord}}}
\newcommand{\bOS}{{\bm{OS}}}
\newcommand{\bERR}{{\bm{ERR}}}
\newcommand{\bSP}{{\bm{SP}}}
\newcommand{\bMV}{{\bm{MV}}}
\newcommand{\bBM}{{\bm{BM}}}
\newcommand{\balpha}{{\bm{\alpha}}}
\newcommand{\balphapre}{{\bm{\alpha}^{\rm pre}}}
\newcommand{\bbeta}{{\bm{\beta}}}
\newcommand{\bgamma}{{\bm{\gamma}}}
\newcommand{\bdelta}{{\bm{\delta}}}
\newcommand{\bkappa}{{\bm{\kappa}}}
\newcommand{\bmu}{{\bm{\mu}}}
\newcommand{\bomega}{{\bm{\omega}}}
\newcommand{\bsigma}{{\bm{\sigma}}}
\newcommand{\btau}{{\bm{\tau}}}
\newcommand{\bphi}{{\bm{\phi}}}
\newcommand{\bpsi}{{\bm{\psi}}}
\newcommand{\bzeta}{{\bm{\zeta}}}
\newcommand{\bone}{{\bm{1}}}
\newcommand{\bzero}{{\bm{0}}}

\newcommand{\sfa}{{{\sf a}}}
\newcommand{\sfb}{{{\sf b}}}
\newcommand{\sfc}{{{\sf c}}}
\newcommand{\sfd}{{{\sf d}}}
\newcommand{\sfe}{{{\sf e}}}
\newcommand{\sff}{{{\sf f}}}
\newcommand{\sfg}{{{\sf g}}}
\newcommand{\sfh}{{{\sf h}}}
\newcommand{\sfi}{{{\sf i}}}
\newcommand{\bsfa}{{\mbox{\textsf{\textbf{a}}}}}
\newcommand{\bsfb}{{\mbox{\textsf{\textbf{b}}}}}
\newcommand{\bsfc}{{\mbox{\textsf{\textbf{c}}}}}
\newcommand{\bsfd}{{\mbox{\textsf{\textbf{d}}}}}
\newcommand{\bsfe}{{\mbox{\textsf{\textbf{e}}}}}
\newcommand{\bsff}{{\mbox{\textsf{\textbf{f}}}}}
\newcommand{\bsfg}{{\mbox{\textsf{\textbf{g}}}}}
\newcommand{\bsfh}{{\mbox{\textsf{\textbf{h}}}}}
\newcommand{\bsfi}{{\mbox{\textsf{\textbf{i}}}}}

\newcommand{\Cbar}{{\overline{C}}}
\newcommand{\Dbar}{{\overline{D}}}
\newcommand{\dbar}{{\overline{d}}}
\def\Ctilde{{\widetilde{C}}}
\def\Ftilde{{\widetilde{F}}}
\def\Gtilde{{\widetilde{G}}}
\def\Htilde{{\widetilde{H}}}
\def\Ptilde{{\widetilde{P}}}
\def\Chat{{\widehat{C}}}
\def\ctilde{{\widetilde{c}}}
\def\zbar{{\overline{Z}}}
\def\pitilde{{\widetilde{\pi}}}
\def\omegahat{{\widehat{\omega}}}

\newcommand{\LD}{{\mathbf{LD}}}
\newcommand{\e}{{\rm e}}
\newcommand{\ecyc}{{\rm ecyc}}
\newcommand{\epa}{{\rm epa}}
\newcommand{\iv}{{\rm iv}}
\newcommand{\pa}{{\rm pa}}
\newcommand{\pk}{{\rm p}}
\newcommand{\val}{{\rm v}}
\newcommand{\da}{{\rm da}}
\newcommand{\dd}{{\rm dd}}
\newcommand{\fp}{{\rm fp}}
\newcommand{\pkcyc}{{\rm pcyc}}
\newcommand{\valcyc}{{\rm vcyc}}
\newcommand{\dacyc}{{\rm dacyc}}
\newcommand{\ddcyc}{{\rm ddcyc}}
\newcommand{\pkpa}{{\rm ppa}}
\newcommand{\valpa}{{\rm vpa}}
\newcommand{\dapa}{{\rm dapa}}
\newcommand{\ddpa}{{\rm ddpa}}

\newcommand{\yp}{{y_\pk}}
\newcommand{\yptilde}{{\widetilde{y}_\pk}}
\newcommand{\yphat}{{\widehat{y}_\pk}}
\newcommand{\yv}{{y_\val}}
\newcommand{\yiv}{{y_\iv}}
\newcommand{\yda}{{y_\da}}
\newcommand{\ydatilde}{{\widetilde{y}_\da}}
\newcommand{\ydd}{{y_\dd}}
\newcommand{\yddtilde}{{\widetilde{y}_\dd}}
\newcommand{\yfp}{{y_\fp}}
\newcommand{\zp}{{z_\pk}}
\newcommand{\zv}{{z_\val}}
\newcommand{\zda}{{z_\da}}
\newcommand{\zdd}{{z_\dd}}

\newcommand{\sech}{{\rm sech}}

\newcommand{\sinv}{\sigma^{-1}}

%
% Jacobian and Dixonian elliptic functions
%
\newcommand{\sn}{{\rm sn}}
\newcommand{\cn}{{\rm cn}}
\newcommand{\dn}{{\rm dn}}
\newcommand{\sm}{{\rm sm}}
\newcommand{\cm}{{\rm cm}}

%
% Commands for hypergeometric series
%
\newcommand{\zfz}{ {{}_0 \! F_0} }
\newcommand{\zfo}{ {{}_0  F_1} }
\newcommand{\ofz}{ {{}_1 \! F_0} }
\newcommand{\ofo}{ {{}_1 \! F_1} }
\newcommand{\oft}{ {{}_1 \! F_2} }
%\newcommand{\tfo}{ {{}_2 \! F_1} }

%
% Hypergeometric functions, using "tensor" package
%
\newcommand{\FHyper}[2]{ {\tensor[_{#1 \!}]{F}{_{#2}}\!} }
\newcommand{\FHYPER}[5]{ {\FHyper{#1}{#2} \!\biggl(
   \!\!\begin{array}{c} #3 \\[1mm] #4 \end{array}\! \bigg|\, #5 \! \biggr)} }
\newcommand{\tfo}{ {\FHyper{2}{1}} }
\newcommand{\tfz}{ {\FHyper{2}{0}} }
\newcommand{\threefz}{ {\FHyper{3}{0}} }
\newcommand{\FHYPERbottomzero}[3]{ {\FHyper{#1}{0} \hspace*{-0mm}\biggl(
   \!\!\begin{array}{c} #2 \\[1mm] \hbox{---} \end{array}\! \bigg|\, #3 \! \biggr)} }
\newcommand{\FHYPERtopzero}[3]{ {\FHyper{0}{#1} \hspace*{-0mm}\biggl(
   \!\!\begin{array}{c} \hbox{---} \\[1mm] #2 \end{array}\! \bigg|\, #3 \! \biggr)} }

\newcommand{\phiHyper}[2]{ {\tensor[_{#1}]{\phi}{_{#2}}} }
\newcommand{\psiHyper}[2]{ {\tensor[_{#1}]{\psi}{_{#2}}} }
\newcommand{\PhiHyper}[2]{ {\tensor[_{#1}]{\Phi}{_{#2}}} }
\newcommand{\PsiHyper}[2]{ {\tensor[_{#1}]{\Psi}{_{#2}}} }
\newcommand{\phiHYPER}[6]{ {\phiHyper{#1}{#2} \!\left(
   \!\!\begin{array}{c} #3 \\ #4 \end{array}\! ;\, #5, \, #6 \! \right)\!} }
\newcommand{\psiHYPER}[6]{ {\psiHyper{#1}{#2} \!\left(
   \!\!\begin{array}{c} #3 \\ #4 \end{array}\! ;\, #5, \, #6 \! \right)} }
\newcommand{\PhiHYPER}[5]{ {\PhiHyper{#1}{#2} \!\left(
   \!\!\begin{array}{c} #3 \\ #4 \end{array}\! ;\, #5 \! \right)\!} }
\newcommand{\PsiHYPER}[5]{ {\PsiHyper{#1}{#2} \!\left(
   \!\!\begin{array}{c} #3 \\ #4 \end{array}\! ;\, #5 \! \right)\!} }
\newcommand{\zerophizero}{ {\phiHyper{0}{0}} }
\newcommand{\ophizero}{ {\phiHyper{1}{0}} }
\newcommand{\zphio}{ {\phiHyper{0}{1}} }
\newcommand{\ophio}{ {\phiHyper{1}{1}} }
\newcommand{\tphio}{ {\phiHyper{2}{1}} }
\newcommand{\tphiz}{ {\phiHyper{2}{0}} }
\newcommand{\tPhio}{ {\PhiHyper{2}{1}} }
\newcommand{\opsio}{ {\psiHyper{1}{1}} }

%
% Variants of \binom  (defined using the AMS "genfrac" command)
%
\newcommand{\stirlingsubset}[2]{\genfrac{\{}{\}}{0pt}{}{#1}{#2}}
\newcommand{\stirlingcycle}[2]{\genfrac{[}{]}{0pt}{}{#1}{#2}}
\newcommand{\assocstirlingsubset}[3]{{\genfrac{\{}{\}}{0pt}{}{#1}{#2}}_{\! \ge #3}}
\newcommand{\genstirlingsubset}[4]{{\genfrac{\{}{\}}{0pt}{}{#1}{#2}}_{\! #3,#4}}
\newcommand{\irredstirlingsubset}[2]{{\genfrac{\{}{\}}{0pt}{}{#1}{#2}}^{\!\rm irr}}
\newcommand{\euler}[2]{\genfrac{\langle}{\rangle}{0pt}{}{#1}{#2}}
\newcommand{\eulergen}[3]{{\genfrac{\langle}{\rangle}{0pt}{}{#1}{#2}}_{\! #3}}
\newcommand{\eulersecond}[2]{\left\langle\!\! \euler{#1}{#2} \!\!\right\rangle}
\newcommand{\eulersecondgen}[3]{{\left\langle\!\! \euler{#1}{#2} \!\!\right\rangle}_{\! #3}}
\newcommand{\binomvert}[2]{\genfrac{\vert}{\vert}{0pt}{}{#1}{#2}}
\newcommand{\binomsquare}[2]{\genfrac{[}{]}{0pt}{}{#1}{#2}}
\newcommand{\doublebinom}[2]{\left(\!\! \binom{#1}{#2} \!\!\right)}

% Array for subscripts

\newenvironment{sarray}{
             \textfont0=\scriptfont0
             \scriptfont0=\scriptscriptfont0
             \textfont1=\scriptfont1
             \scriptfont1=\scriptscriptfont1
             \textfont2=\scriptfont2
             \scriptfont2=\scriptscriptfont2
             \textfont3=\scriptfont3
             \scriptfont3=\scriptscriptfont3
           \renewcommand{\arraystretch}{0.7}
           \begin{array}{l}}{\end{array}}

\newenvironment{scarray}{
             \textfont0=\scriptfont0
             \scriptfont0=\scriptscriptfont0
             \textfont1=\scriptfont1
             \scriptfont1=\scriptscriptfont1
             \textfont2=\scriptfont2
             \scriptfont2=\scriptscriptfont2
             \textfont3=\scriptfont3
             \scriptfont3=\scriptscriptfont3
           \renewcommand{\arraystretch}{0.7}
           \begin{array}{c}}{\end{array}}

% Circled math symbols:
% From http://latex-community.org/forum/viewtopic.php?f=44&t=22367

%\usepackage{tikz}
\newcommand*\circled[1]{\tikz[baseline=(char.base)]{
  \node[shape=circle,draw,inner sep=1pt] (char) {#1};}}
\newcommand{\ostar}{{\circledast}}
\newcommand{\ostarN}{{\,\circledast_{\vphantom{\dot{N}}N}\,}}
\newcommand{\ostarPsi}{{\,\circledast_{\vphantom{\dot{\Psi}}\Psi}\,}}
\newcommand{\starN}{{\,\ast_{\vphantom{\dot{N}}N}\,}}
\newcommand{\starpsi}{{\,\ast_{\vphantom{\dot{\bpsi}}\!\bpsi}\,}}
\newcommand{\starone}{{\,\ast_{\vphantom{\dot{1}}1}\,}}
\newcommand{\startwo}{{\,\ast_{\vphantom{\dot{2}}2}\,}}
\newcommand{\starinfty}{{\,\ast_{\vphantom{\dot{\infty}}\infty}\,}}
\newcommand{\starT}{{\,\ast_{\vphantom{\dot{T}}T}\,}}

%% For scaling equations (uses "graphicx" package):  see
%% http://tex.stackexchange.com/questions/60453/reducing-font-size-in-equation
\newcommand*{\Scale}[2][4]{\scalebox{#1}{$#2$}}

\newcommand*{\Scaletext}[2][4]{\scalebox{#1}{#2}} %% THIS DOESN'T SEEM TO WORK

\clearpage

\tableofcontents

\clearpage

\section{Introduction}

A permutation $\sigma$ is called \textbfit{cycle-alternating}
if it has no cycle double rises, cycle double falls, or fixed points;
thus, each cycle of $\sigma$ is of even length (call it~$2m$)
and consists of $m$ cycle valleys and $m$ cycle peaks in alternation.\footnote{
{\bf Warning:} The phrase ``cycle-alternating permutations'' 
has previously been used by \mbox{Josuat-Verg\`es}
in \cite{Josuat-Verges_14} to describe a different class of objects.
We think that his objects should instead be called 
cycle-alternating {\em signed} permutations.
}
Deutsch and Elizalde \cite[Proposition~2.2]{Deutsch_11}
showed that the number of cycle-alternating permutations of $[2n]$
is the secant number $E_{2n}$
(see also Dumont \cite[pp.~37, 40]{Dumont_86}
 and Biane \cite[section~6]{Biane_93}).
In this paper we would like to present some new continued fractions
for multivariate polynomials that enumerate cycle-alternating permutations
with respect to a large set of statistics (8 or 16 variables).

But first it may be appropriate to place this project in a wider context.
In a recent article \cite{Sokal-Zeng_masterpoly}, Zeng and one of us
presented Stieltjes-type and Jacobi-type continued fractions
for some ``master polynomials''
that enumerate permutations, set partitions or perfect matchings
with respect to a large (sometimes infinite) number of independent statistics.
These polynomials systematize what we think of as the ``linear family'':
namely, sequences in which
the Stieltjes continued-fraction coefficients $(\alpha_n)_{n \ge 1}$
grow linearly in $n$.
For instance, Euler \cite[section~21]{Euler_1760} showed in 1746 that
\be
   \sum_{n=0}^\infty n! \: t^n
   \;=\;
   \cfrac{1}{1 - \cfrac{1t}{1 - \cfrac{1t}{1 - \cfrac{2t}{1- \cfrac{2t}{1-\cdots}}}}}
   \;.
 \label{eq.nfact.contfrac}
\ee
In other words,
the ordinary generating function (ogf) of the sequence $a_n = n!$
has a Stieltjes-type continued fraction with coefficients
$\alpha_{2k-1} = \alpha_{2k} = k$.
This inspired us to introduce the polynomials $P_n(x,y,u,v)$ defined
by the continued fraction
\be
   \sum_{n=0}^\infty P_n(x,y,u,v) \: t^n
   \;=\;
   \cfrac{1}{1 - \cfrac{xt}{1 - \cfrac{yt}{1 - \cfrac{(x+u)t}{1- \cfrac{(y+v)t}{1 - \cfrac{(x+2u)t}{1 - \cfrac{(y+2v)t}{1-\cdots}}}}}}}
 \label{eq.eulerian.fourvar.contfrac}
\ee
with coefficients
\begin{subeqnarray}
   \alpha_{2k-1}  & = &  x + (k-1) u \\
   \alpha_{2k}    & = &  y + (k-1) v
 \label{def.weights.eulerian.fourvar}
\end{subeqnarray}
Clearly $P_n(x,y,u,v)$ is a homogeneous polynomial of degree $n$,
with nonnegative integer coefficients;
and since $P_n(1,1,1,1) = n!$,
which enumerates permutations of an $n$-element set,
it was plausible to expect that $P_n(x,y,u,v)$ enumerates permutations of $[n]$
according to some natural fourfold classification of indices in $[n]$.
Our first theorem \cite[Theorem~2.1]{Sokal-Zeng_masterpoly}
provided two alternative versions of that classification.

But there are other sequences of integers,
and hence other combinatorial models, that also belong to the ``linear family''.
For instance, in that same paper
Euler also showed \cite[section~29]{Euler_1760}
(see also \cite{Euler_1788})
that the generating function of the odd semifactorials can be represented
as a Stieltjes-type continued fraction
\be
   \sum_{n=0}^\infty (2n-1)!! \: t^n
   \;=\;
   \cfrac{1}{1 - \cfrac{1t}{1 - \cfrac{2t}{1 - \cfrac{3t}{1- \cdots}}}}
 \label{eq.2n-1semifact.contfrac}
\ee
with coefficients $\alpha_n = n$.
Inspired by \reff{eq.2n-1semifact.contfrac},
we introduced the polynomials $M_n(x,y,u,v)$ defined
by the continued fraction
\be
   \sum_{n=0}^\infty M_n(x,y,u,v) \: t^n
   \;=\;
   \cfrac{1}{1 - \cfrac{xt}{1 - \cfrac{(y+v)t}{1 - \cfrac{(x+2u)t}{1- \cfrac{(y+3v)t}{1 - \cfrac{(x+4u)t}{1 - \cfrac{(y+5v)t}{1-\cdots}}}}}}}
 \label{eq.matching.fourvar.contfrac}
\ee
with coefficients
\begin{subeqnarray}
   \alpha_{2k-1}  & = &  x + (2k-2) u \\
   \alpha_{2k}    & = &  y + (2k-1) v
 \label{def.weights.matching.fourvar}
\end{subeqnarray}
Clearly $M_n(x,y,u,v)$ is a homogeneous polynomial of degree $n$,
with nonnegative integer coefficients;
and since $M_n(1,1,1,1) = (2n-1)!!$,
which enumerates perfect matchings of a $2n$-element set,
it was plausible to expect that $M_n(x,y,u,v)$ enumerates
perfect matchings of $[2n]$
(or equivalently, fixed-point-free involutions of $[2n]$)
according to some natural fourfold classification
associated to half of the indices in $[2n]$.
In \cite[Theorem~4.1]{Sokal-Zeng_masterpoly}
we identified such a classification,
involving even and odd cycle peaks and their antirecord status.

In a very recent paper \cite{Deb-Sokal_genocchi},
the two present authors took one step up in complexity,
to consider the ``quadratic family'',
in~which the $(\alpha_n)_{n \ge 1}$ grow quadratically in $n$.
For instance, we could consider
\begin{subeqnarray}
   \alpha_{2k-1}  & = &  [x_1 + (k-1) u_1] \, [x_2 + (k-1) u_2] \\
   \alpha_{2k}    & = &  [y_1 + (k-1) v_1] \, [y_2 + (k-1) v_2]
 \label{eq.quadratic_family}
\end{subeqnarray}
which can be thought of as the product of two independent sets
of coefficients \reff{def.weights.eulerian.fourvar}.
With all parameters set to 1,
these coefficients $\alpha_{2k-1} = \alpha_{2k} = k^2$
correspond to the continued fraction
\cite[eq.~(9.7)]{Viennot_81} \cite[p.~V-15]{Viennot_83}
for the median Genocchi numbers \cite[A005439]{OEIS};
so it was natural to seek a combinatorial model
that is enumerated by the median Genocchi numbers.
In \cite{Deb-Sokal_genocchi} we considered D-permutations
\cite{Lazar_22,Lazar_20}
and various subclasses thereof (D-semiderangements and D-derangements),
and found Stieltjes-type and Thron-type continued fractions
for some multivariate polynomials enumerating these objects
with respect to a large set of statistics (12 or 22 variables).

But just as with the ``linear family'',
there are other sequences of integers,
and hence other combinatorial models,
that also belong to the ``quadratic family''.
For instance, the secant numbers $E_{2n}$ \cite[A000364]{OEIS}
can be represented by the continued fraction
\cite[p.~H9]{Stieltjes_1889}
\cite[p.~77]{Rogers_07}
\cite[Theorem~3B(iii)]{Flajolet_80}
\be
   \sum_{n=0}^\infty E_{2n} t^n
   \;=\;
   \cfrac{1}{1 - \cfrac{1^2 t}{1 - \cfrac{2^2 t}{1 -  \cfrac{3^2 t}{1- \cdots}}}}
 \label{eq.cfrac.secant}
\ee
with coefficients $\alpha_n = n^2$.
In \cite[Section~2.15]{Sokal-Zeng_masterpoly}
this was generalized by introducing the polynomials $Q_n(x,y,u,v)$
whose ogf has an S-fraction with coefficients
$\alpha_n = [x+ (n-1)u] \, [y + (n-1)v]$,
and they were interpreted as enumerating cycle-alternating permutations
of $[2n]$ according to the antirecord status of cycle peaks and
the record status of cycle valleys.

In the present paper we shall go farther, and introduce polynomials
$Q_n(x',y',u',v',$ $x'',y'',u'',v'')$
whose ogf has an S-fraction with coefficients
\begin{subeqnarray}
   \alpha_{2k-1}  & = &  [x' + (2k-2) u'] \: [y' + (2k-2) v'] \\[1mm]
   \alpha_{2k}    & = &  [x'' + (2k-1) u''] \: [y'' + (2k-1) v'']
 \label{eq.quadratic_family_bis}
\end{subeqnarray}
which can be thought of as the product of two independent sets
of coefficients \reff{def.weights.matching.fourvar}.
When $x'=x''$, $y'=y''$, $u'=u''$, $v'=v''$
this reduces to the previous polynomials.
We will interpret these 8-variable polynomials
as enumerating cycle-alternating permutations of $[2n]$
according to the antirecord status of cycle peaks and
the record status of cycle valleys,
{\em separately for even and odd indices}\/
(Theorem~\ref{thm.first.Sfrac}).\footnote{
   This answers the question posed in
   \cite[last sentence of Section~2.15]{Sokal-Zeng_masterpoly}.
}
We will also obtain a generalization of this result
to 16-variable polynomials that include four pairs of variables $(p,q)$
that count crossings and nestings
(Theorem~\ref{thm.first.Sfrac.pqgen}).
The proofs will in fact be easy applications
of the ``master'' S-fraction for cycle-alternating permutations
\cite[Theorem~2.20]{Sokal-Zeng_masterpoly}.
Finally, we will obtain variants that include the counting of cycles,
at the expense of renouncing taking account of
the record status of cycle valleys
(Theorems~\ref{thm.second.Sfrac} and \ref{thm.second.Sfrac.pqgen}).

This leaves us with the following analogy between
the linear and quadratic families:

\medskip

\begin{center}
\begin{tabular}{rcc}
   ``Linear family'':  & Permutations & Perfect matchings \\
          & $\alpha_{2k-1} = \alpha_{2k} = k$  &  $\alpha_n = n$ \\
          & Generalized to \reff{def.weights.eulerian.fourvar}
            & Generalized to \reff{def.weights.matching.fourvar} \\[2mm]
          &  $\Scale[2]{\downarrow}$   &   $\Scale[2]{\downarrow}$  \\[2mm]
   ``Quadratic family'':  & D-permutations & Cycle-alternating permutations \\
          & $\alpha_{2k-1} = \alpha_{2k} = k^2$  &  $\alpha_n = n^2$ \\
          & Generalized to \reff{eq.quadratic_family}
            & Generalized to \reff{eq.quadratic_family_bis}
\end{tabular}
\end{center}

\medskip

\noindent
The only slight flaw in the analogy is that the permutations are on $[n]$,
while the other three models are on $[2n]$.

This suggests to pursue the analogy one step further:
the third combinatorial model in the ``linear family''
is set partitions, with $a_n = B_n$ (Bell number)
and $\alpha_{2k-1} = 1$, $\alpha_{2k} = k$.
The corresponding member of the ``quadratic family''
would have $\alpha_{2k-1} = 1$, $\alpha_{2k} = k^2$,
which leads to the sequence
\be
   (a_n)_{n \ge 0}
   \;=\;
   1, 1, 2, 5, 17, 78, 461, 3417, 31134, 340217, 4365673, \ldots
   \;.
\ee
Alas, this sequence is not at present in the OEIS \cite{OEIS},
and we have no idea what combinatorial model it might represent.

\bigskip

We can, in fact, place our results on cycle-alternating permutations
in a wider context: namely, that of Laguerre digraphs.
A \textbfit{Laguerre digraph}
\cite{Foata_84,Sokal_multiple_laguerre,latpath_laguerre}
is a digraph in which each vertex has out-degree 0 or~1 and in-degree 0 or~1.
Therefore, each weakly connected component of a Laguerre digraph
is either a directed path of some length $\ell \ge 0$
(where a path of length 0 is an isolated vertex)
or else a directed cycle of some length $\ell \ge 1$
(where a cycle of length 1 is a loop).
A Laguerre digraph with no paths is simply a collection of directed cycles,
i.e.~the digraph associated to a permutation;
so Laguerre digraphs generalize permutations.
We shall then define a subclass of Laguerre digraphs,
called \textbfit{alternating Laguerre digraphs},
that similarly generalizes cycle-alternating permutations.
We shall show that the generalized Stieltjes--Rogers polynomials
associated to the cycle-counting variant of the
S-fraction \reff{eq.quadratic_family_bis},
as well as its generalizations, can be interpreted as enumerating
alternating Laguerre digraphs with suitable weights.

\bigskip

The plan of this paper is as follows:
In Sections~\ref{sec.prelim1} and \ref{sec.prelim2}
we review some needed definitions and facts
concerning continued fractions and permutation statistics, respectively.
In Section~\ref{sec.S-fraction1}
we state and prove our ``first'' S-fractions
for cycle-alternating permutations
(namely, the ones not including the counting of cycles);
and in Section~\ref{sec.S-fraction2}
we do the same for our ``second'' S-fractions
(namely, those including the counting of cycles).
In Section~\ref{sec.alternating} we put cycle-alternating permutations
into the wider context of alternating Laguerre digraphs,
and then use this to interpret the generalized Stieltjes--Rogers polynomials
associated to some of our continued fractions.

\section{Preliminaries: Continued fractions}   \label{sec.prelim1}

In this section we define the objects associated to continued fractions
that will be used in the remainder of the paper,
and review some of their properties.

\subsection{Motzkin and Dyck paths}

A \textbfit{Motzkin path}
is a path in the upper half-plane $\Z \times \N$,
starting and ending on the horizontal axis,
using the steps $(1,1)$ [``rise''], $(1,0)$ [``level step'']
and $(1,-1)$ [``fall''].
More generally, a \textbfit{partial Motzkin path}
is a path in the upper half-plane $\Z \times \N$,
starting on the horizontal axis but ending anywhere,
using the steps $(1,1)$, $(1,0)$ and $(1,-1)$.

A \textbfit{Dyck path}
is a path in the upper half-plane $\Z \times \N$,
starting and ending on the horizontal axis,
using the steps $(1,1)$ [``rise''] and $(1,-1)$ [``fall''];
that is, it is a Motzkin path that has no level steps.
Obviously a Dyck path must have even length.
More generally, a \textbfit{partial Dyck path}
is a path in the upper half-plane $\Z \times \N$,
starting on the horizontal axis but ending anywhere,
using the steps $(1,1)$ and $(1,-1)$.

We shall usually take all these paths to start at the origin $(0,0)$.
Then any partial Dyck path must stay on the even sublattice
$\{(n,k) \colon\, n+k \textrm{ is even} \}$.

\subsection{J-fractions and Jacobi--Rogers polynomials}

Let $\bbeta = (\beta_i)_{i \ge 1}$ and $\bgamma = (\gamma_i)_{i \ge 0}$
be indeterminates.
Following Flajolet \cite{Flajolet_80},
we define the \textbfit{Jacobi--Rogers polynomials} $J_n(\bbeta,\bgamma)$
to be the Taylor coefficients of the generic J-fraction:
\be
   \sum_{n=0}^\infty J_n(\bbeta,\bgamma) \, t^n
   \;=\;
   \cfrac{1}{1 - \gamma_0 t - \cfrac{\beta_1 t^2}{1 - \gamma_1 t - \cfrac{\beta_2 t^2}{1 - \cdots}}}
   \;\;.
 \label{eq.Jtype.cfrac}
\ee
Clearly $J_n(\bbeta,\bgamma)$ is a polynomial in $\bbeta$ and $\bgamma$,
with nonnegative integer coefficients.
Flajolet \cite{Flajolet_80} showed that
the Jacobi--Rogers polynomial $J_n(\bbeta,\bgamma)$
is the generating polynomial for Motzkin paths of length $n$,
in which each rise gets weight~1,
each fall from height~$i$ gets weight $\beta_i$,
and each level step at height~$i$ gets weight $\gamma_i$.

Now define the \textbfit{generalized Jacobi--Rogers polynomial}
$J_{n,k}(\bbeta,\bgamma)$ to be the generating polynomial
for partial Motzkin paths from $(0,0)$ to $(n,k)$,
in which each rise gets weight~1,
each fall from height~$i$ gets weight $\beta_i$,
and each level step at height~$i$ gets weight $\gamma_i$.
Obviously $J_{n,k}$ is nonvanishing only for $0 \le k \le n$,
so we have an infinite lower-triangular array
$\sfJ = \big( J_{n,k}(\bbeta,\bgamma) \big)_{\! n,k \ge 0}$
in which the zeroth column displays
the ordinary Jacobi--Rogers polynomials $J_{n,0} = J_n$.
On the diagonal we have $J_{n,n} = 1$,
and on the first subdiagonal we have $J_{n,n-1} = \sum_{i=0}^{n-1} \gamma_i$.
By considering the last step of the path, we see that
the polynomials $J_{n,k}(\bbeta,\bgamma)$ satisfy the recurrence
\be
   J_{n+1,k}
   \;=\;
   J_{n,k-1}  \:+\:  \gamma_k J_{n,k} \:+\:  \beta_{k+1} J_{n,k+1}
 \label{eq.Jnk.recursion}
\ee
with the initial condition $J_{0,k} = \delta_{k0}$
(where of course we set $J_{n,-1} = 0$).

A redundant generalization of these polynomials,
in which we also weight rises in the Motzkin path,
will be useful later.
Let $\bsfa = (\sfa_i)_{i \ge 0}, \bsfb = (\sfb_i)_{i \ge 1}$ 
and $\bsfc = (\sfc_i)_{i \ge 0}$ be indeterminates.
Define $J_{n,k}(\bsfa,\bsfb,\bsfc)$ to be the generating polynomial
for partial Motzkin paths from $(0,0)$ to $(n,k)$,
in which each rise from height~$i$ gets weight~$\sfa_i$,
each fall from height~$i$ gets weight $\sfb_i$,
and each level step at height~$i$ gets weight $\sfc_i$.
By considering the last step of the path, we see that these polynomials
satisfy the recurrence
\be
   J_{n+1,k}
   \;=\;
   \sfa_{k-1} J_{n,k-1}  \:+\:  \sfc_k J_{n,k} \:+\:  \sfb_{k+1} J_{n,k+1}
   \;.
 \label{eq.Jnk.prime.recursion}
\ee
The recurrences \reff{eq.Jnk.recursion} and \reff{eq.Jnk.prime.recursion} 
then easily imply the following proposition:

\begin{proposition}
   \label{prop.JR.bothkinds}
The polynomials $J_{n,k}(\bsfa, \bsfb, \bsfc)$
are related to the polynomials $J_{n,k}(\bbeta,\bgamma)$
with weights
\begin{subeqnarray}
   \beta_n &=& \sfa_{n-1} \sfb_n\\
   \gamma_n &=& \sfc_{n}
   \label{eq.prop.JRpoly.bothkinds}
\end{subeqnarray}
by
\be
     J_{n,k}(\bsfa, \bsfb, \bsfc)
     \;=\;
     \left(\prod\limits_{i=0}^{k-1}\sfa_i \right) J_{n,k}(\bbeta, \bgamma)
     \;.
  \label{eq.prop.JRmat.bothkinds}
\ee
\end{proposition}

\noindent
\begin{sloppypar}
We remark that since $J_{n,k}(\bbeta, \bgamma)$
is by definition a polynomial,
it follows that the polynomial $J_{n,k}(\bsfa, \bsfb, \bsfc)$
necessarily has a factor $\prod\limits_{i=0}^{k-1}\sfa_i$,
and that
$J_{n,k}(\bbeta, \bgamma) 
	= J_{n,k}(\bsfa, \bsfb, \bsfc)
 \Big/\prod\limits_{i=0}^{k-1}\sfa_i$.
\end{sloppypar}

\medskip

\proof
Let us use the shorthands $J_{n,k} = J_{n,k}(\bbeta, \bgamma)$
and $J'_{n,k} = J_{n,k}(\bsfa, \bsfb, \bsfc)$;
and define $p_k = \prod\limits_{i=0}^{k-1} \sfa_i$.
The recurrence \reff{eq.Jnk.recursion},
with weights \reff{eq.prop.JRpoly.bothkinds}, is
\be
   J_{n+1,k} \;=\;
   J_{n,k-1} \,+\, \sfc_k J_{n,k} \,+\, \sfa_k\sfb_{k+1} J_{n,k+1}
   \;.
 \label{eq.genJRpoly.rec.bothkinds}
\ee
Multiplying both sides of \reff{eq.genJRpoly.rec.bothkinds}
with $p_{k}$ yields
\begin{subeqnarray}
J_{n+1,k} \, p_k & = & p_k \left( J_{n,k-1} \,+\,
   \sfc_k J_{n,k} \,+\, \sfa_k\sfb_{k+1} J_{n,k+1}
   \right)\\
   & = & \sfa_{k-1} \left( J_{n,k-1} \, p_{k-1} \right) \,+\,
   \sfc_k \left( J_{n,k} \, p_k\right) \,+\,
   \sfb_{k+1} \left( J_{n,k+1} \, p_{k+1} \right)
 \slabel{eq.rec.Jnkck.bothkinds}
\end{subeqnarray}
Comparing the coefficients in \reff{eq.rec.Jnkck.bothkinds} and 
\reff{eq.Jnk.prime.recursion},
we notice that the recurrences for the polynomials 
$J_{n,k} \, p_k$ and $J_{n,k}'$ are exactly the same;
and of course the initial conditions are the same as well.
This shows that $J_{n,k}' = J_{n,k} \, p_k$,
as claimed in \reff{eq.prop.JRmat.bothkinds}.
\qed

\bigskip

One important algebraic tool for studying J-fractions is
\textbfit{Rogers' addition formula}
\cite{Rogers_07} \cite[Theorem~53.1]{Wall_48},
which we state in the form given in \cite[Theorem~9.9]{Sokal_alg_contfrac}:

\begin{theorem}[Rogers' addition formula]
   \label{thm.rogers}
The column exponential generating functions of
the matrix of generalized Jacobi--Rogers polynomials,
\be
   \scrj_k(t;\bbeta,\bgamma)
   \;\eqdef\;
   \sum_{n=k}^\infty J_{n,k}(\bbeta,\bgamma) \, {t^n \over n!}
   \;,
\ee
satisfy
\be
   \scrj_0(t+u;\bbeta,\bgamma)
   \;=\;
   \sum\limits_{k=0}^\infty \beta_1 \cdots \beta_k \,
           \scrj_k(t;\bbeta,\bgamma) \, \scrj_k(u;\bbeta,\bgamma)
   \;.
\ee

And conversely, if $A(t)$ and $F_0(t), F_1(t), \ldots$ are formal power series
(with elements in a commutative ring $R$ containing the rationals)
satisfying
\be
   A(t) \;=\; 1 + O(t) \,,\qquad
   F_k(t) \;=\; {t^k \over k!} \,+\, \mu_k {t^{k+1} \over (k+1)!} \,+\,
                  O(t^{k+2})
 \label{eq.thm.rogers.1}
\ee
and
\be
   A(t+u)
   \;=\;
   \sum\limits_{k=0}^\infty \beta_1 \cdots \beta_k \, F_k(t) \, F_k(u)
 \label{eq.thm.rogers.2}
\ee
for some {\em regular} elements $\bbeta = (\beta_k)_{k \ge 1}$,
then $A(t) = F_0(t)$ and $F_k(t) = \scrj_k(t;\bbeta,\bgamma)$
with the given $\bbeta$
and with $\gamma_k = \mu_k - \mu_{k-1}$ (where $\mu_{-1} \eqdef 0$).
\end{theorem}

\noindent
Here an element of a commutative ring $R$ is called
\textbfit{regular} if it is neither zero nor a divisor of zero.

\subsection{S-fractions and Stieltjes--Rogers polynomials}

Let $\balpha = (\alpha_i)_{i \ge 1}$ be indeterminates.
Again following Flajolet \cite{Flajolet_80},
we define the \textbfit{Stieltjes--Rogers polynomials} $S_n(\balpha)$
to be the Taylor coefficients of the generic S-fraction:
\be
   \sum_{n=0}^\infty S_n(\balpha) \, t^n
   \;=\;
   \cfrac{1}{1 - \cfrac{\alpha_1 t}{1 - \cfrac{\alpha_2 t}{1- \cdots}}}
   \;\;.
 \label{eq.Stype.cfrac}
\ee
Clearly $S_n(\balpha)$ is a homogeneous polynomial of degree $n$ in $\balpha$,
with nonnegative integer coefficients.
Flajolet \cite{Flajolet_80} showed that
the Stieltjes--Rogers polynomial $S_n(\balpha)$
is the generating polynomial for Dyck paths of length $2n$,
in which each rise gets weight~1
and each fall from height~$i$ gets weight $\alpha_i$.

Now define the
\textbfit{generalized Stieltjes--Rogers polynomial of the first kind}
$S_{n,k}(\balpha)$
to be the generating polynomial for Dyck paths
starting at $(0,0)$ and ending at $(2n,2k)$,
in which each rise gets weight 1
and each fall from height~$i$ gets weight $\alpha_i$.
Obviously $S_{n,k}$ is nonvanishing only for $0 \le k \le n$,
so we have an infinite lower-triangular array
$\sfS = (S_{n,k}(\balpha))_{n,k \ge 0}$
in which the zeroth column displays
the ordinary Stieltjes--Rogers polynomials $S_{n,0} = S_n$.
We have $S_{n,n} = 1$ and $S_{n,n-1} = \sum_{i=1}^{2n-1} \alpha_i$.

Likewise, let us define
the \textbfit{generalized Stieltjes--Rogers polynomial of the second kind}
$S'_{n,k}(\balpha)$
to be the generating polynomial for Dyck paths
starting at $(0,0)$ and ending at $(2n+1,2k+1)$,
in which again each rise gets weight 1
and each fall from height~$i$ gets weight $\alpha_i$.
Since $S'_{n,k}$ is nonvanishing only for $0 \le k \le n$,
we obtain a second infinite lower-triangular array
$\sfS' = (S'_{n,k}(\balpha))_{n,k \ge 0}$.
We have $S'_{n,n} = 1$ and $S'_{n,n-1} = \sum_{i=1}^{2n} \alpha_i$.

The polynomials $S_{n,k}(\balpha)$ and $S'_{n,k}(\balpha)$
manifestly satisfy the joint recurrence
\begin{subeqnarray}
   S'_{n,k}   & = &  S_{n,k} \:+\: \alpha_{2k+2} \, S_{n,k+1}
       \\[2mm]
   S_{n+1,k}  & = &  S'_{n,k-1} \:+\: \alpha_{2k+1} \, S'_{n,k}
 \slabel{eq.SnlSprimenl.recurrence.b}
 \label{eq.SnlSprimenl.recurrence}
\end{subeqnarray}
for $n,k \ge 0$,
with the initial conditions $S_{0,k} = \delta_{k 0}$ and $S'_{n,-1} = 0$.
It follows that the $S_{n,k}$ satisfy the recurrence
\be
   S_{n+1,k}
   \;=\;
   S_{n,k-1}
      \:+\:  (\alpha_{2k} + \alpha_{2k+1}) \, S_{n,k}
      \:+\:  \alpha_{2k+1} \alpha_{2k+2} \, S_{n,k+1}
 \label{eq.Snl.recursion}
\ee
%% with the initial condition $S_{0,k} = \delta_{k 0}$
%% (where we set $S_{n,-1} = 0$ and $\alpha_0 = 0$),
(where $S_{n,-1} = 0$ and $\alpha_0 = 0$),
while the $S'_{n,k}$ satisfy the recurrence
\be
   S'_{n+1,k}
   \;=\;
   S'_{n,k-1}
      \:+\:  (\alpha_{2k+1} + \alpha_{2k+2}) \, S'_{n,k}
      \:+\:  \alpha_{2k+2} \alpha_{2k+3} \, S'_{n,k+1}
   \;.
 \label{eq.Snlprime.recursion}
\ee
%% with the initial condition $S'_{0,k} = \delta_{k 0}$
%% (where we set $S'_{n,-1} = 0$).
Note that \reff{eq.Snl.recursion} and \reff{eq.Snlprime.recursion}
have the same form as \reff{eq.Jnk.recursion},
when $\bbeta$ and $\bgamma$ are defined suitably in terms of the $\balpha$:
these correspondences are examples of \textbfit{contraction formulae}
\cite[p.~21]{Wall_48} \cite[p.~V-31]{Viennot_83}
that express an S-fraction as an equivalent J-fraction.
In particular we have $S_{n,k}(\balpha) = J_{n,k}(\bbeta,\bgamma)$ where
\begin{subeqnarray}
   \gamma_0  & = &  \alpha_1
       \slabel{eq.contraction_even.coeffs.a}   \\
   \gamma_n  & = &  \alpha_{2n} + \alpha_{2n+1}  \qquad\hbox{for $n \ge 1$}
       \hspace*{-2cm}
       \slabel{eq.contraction_even.coeffs.b}   \\
   \beta_n  & = &  \alpha_{2n-1} \alpha_{2n}
       \slabel{eq.contraction_even.coeffs.c}
 \label{eq.contraction_even.coeffs}
\end{subeqnarray}
and $S'_{n,k}(\balpha) = J_{n,k}(\bbeta',\bgamma')$ where
\begin{subeqnarray}
   \gamma'_n  & = &  \alpha_{2n+1} + \alpha_{2n+2}  \\
   \beta'_n  & = &  \alpha_{2n} \alpha_{2n+1}
 \label{eq.contraction_odd.coeffs}
\end{subeqnarray}
The recurrences
\reff{eq.Jnk.recursion}/\reff{eq.Snl.recursion}/\reff{eq.Snlprime.recursion}
define implicitly the (tridiagonal) \textbfit{production matrices}
for $\sfJ$, $\sfS$ and $\sfS'$:
see \cite{Deutsch_05,Deutsch_09,latpath_SRTR}.
Some workers call the arrays $\sfJ$, $\sfS$ and/or $\sfS'$
the \textbfit{Stieltjes table}.
%% ({\em tableau de Stieltjes}\/).

More trivially, a partial Dyck path is simply a partial Motzkin path
with no level steps.
Therefore, the generalized Stieltjes--Rogers polynomials
with weights $\balpha$
can be related to the generalized Jacobi--Rogers polynomials
with weights $(\bbeta,\bgamma) = (\balpha,\bzero)$ by
\begin{subeqnarray}
   S_{n,k}(\balpha)  & = &  J_{2n,2k}(\balpha,\bzero)   \\[2mm]
   S'_{n,k}(\balpha)  & = &  J_{2n+1,2k+1}(\balpha,\bzero)
 \label{eq.Snk.Jnk}
\end{subeqnarray}
and
\be
   J_{n,k}(\balpha,\bzero) \;=\; 0
      \quad\hbox{if $n+k$ is odd}  
   \;.
\ee

\section{Preliminaries: Permutation statistics}   \label{sec.prelim2}

In this section we define the permutation statistics
that will be used in the remainder of the paper.

\subsection{The record-and-cycle classification}
     \label{subsec.statistics.1}

Given a permutation $\sigma \in \Sym_N$, an index $i \in [N]$ is called an
\begin{itemize}
   \item {\em excedance}\/ (exc) if $i < \sigma(i)$;
   \item {\em anti-excedance}\/ (aexc) if $i > \sigma(i)$;
   \item {\em fixed point}\/ (fix) if $i = \sigma(i)$.
\end{itemize}
Clearly every index $i$ belongs to exactly one of these three types;
we call this the \textbfit{excedance classification}.
We also say that $i$ is a {\em weak excedance}\/ if $i \le \sigma(i)$,
and a {\em weak anti-excedance}\/ if $i \ge \sigma(i)$.

A more refined classification is as follows:
an index $i \in [N]$ is called a
\begin{itemize}
   \item {\em cycle peak}\/ (cpeak) if $\sigma^{-1}(i) < i > \sigma(i)$;
   \item {\em cycle valley}\/ (cval) if $\sigma^{-1}(i) > i < \sigma(i)$;
   \item {\em cycle double rise}\/ (cdrise) if $\sigma^{-1}(i) < i < \sigma(i)$;
   \item {\em cycle double fall}\/ (cdfall) if $\sigma^{-1}(i) > i > \sigma(i)$;
   \item {\em fixed point}\/ (fix) if $\sigma^{-1}(i) = i = \sigma(i)$.
\end{itemize}
Clearly every index $i$ belongs to exactly one of these five types;
we refer to this classification as the \textbfit{cycle classification}.
Obviously, excedance = cycle valley or cycle double rise,
and anti-excedance = cycle peak or cycle double fall.
We write
\be
   \Cpeak(\sigma)
   \;=\;
   \{ i \colon\: \sigma^{-1}(i) < i > \sigma(i) \}
\ee
for the set of cycle peaks and
\be
   \cpeak(\sigma)  \;=\;  |\Cpeak(\sigma)|
\ee
for its cardinality, and likewise for the others.

On the other hand, an index $i \in [N]$ is called a
\begin{itemize}
   \item {\em record}\/ (rec) (or {\em left-to-right maximum}\/)
         if $\sigma(j) < \sigma(i)$ for all $j < i$
      [note in particular that the indices 1 and $\sigma^{-1}(N)$
       are always records];
   \item {\em antirecord}\/ (arec) (or {\em right-to-left minimum}\/)
         if $\sigma(j) > \sigma(i)$ for all $j > i$
      [note in particular that the indices $N$ and $\sigma^{-1}(1)$
       are always antirecords];
   \item {\em exclusive record}\/ (erec) if it is a record and not also
         an antirecord;
   \item {\em exclusive antirecord}\/ (earec) if it is an antirecord
         and not also a record;
   \item {\em record-antirecord}\/ (rar) (or {\em pivot}\/)
      if it is both a record and an antirecord;
   \item {\em neither-record-antirecord}\/ (nrar) if it is neither a record
      nor an antirecord.
\end{itemize}
Every index $i$ thus belongs to exactly one of the latter four types;
we refer to this classification as the \textbfit{record classification}.
We write
\be
   \Rec(\sigma)
   \;=\;
   \{ i \colon\: \sigma(j) < \sigma(i) \hbox{ for all } j < i \}
\ee
for the set of record indices and
\be
   \rec(\sigma)  \;=\;  |\Rec(\sigma)|
\ee
for its cardinality, and likewise for antirecords.

The record and cycle classifications of indices are related as follows:
\begin{quote}
\begin{itemize}
   \item[(a)]  Every record is a weak excedance,
      and every exclusive record is an excedance.
   \item[(b)]  Every antirecord is a weak anti-excedance,
      and every exclusive antirecord is an anti-excedance.
   \item[(c)]  Every record-antirecord is a fixed point.
\end{itemize}
\end{quote}
Therefore, by applying the record and cycle classifications simultaneously,
we obtain 10 (not~20) disjoint categories \cite{Sokal-Zeng_masterpoly}:
%we obtain 10~disjoint categories \cite{Sokal-Zeng_masterpoly}:
\medskip
\begin{center}
\begin{tabular}{c|c|c|c|c|c|}
                & cpeak & cval & cdrise & cdfall & fix \\
        \hline
        erec &   & ereccval & ereccdrise & &\\
        earec & eareccpeak & &  & eareccdfall &\\
        rar & & & & & rar \\
        nrar & nrcpeak & nrcval & nrcdrise & nrcdfall & nrfix\\
        \hline
\end{tabular}
\end{center}
\medskip
%\begin{itemize}
%   \item ereccval:  exclusive records that are also cycle valleys;
%   \item ereccdrise:  exclusive records that are also cycle double rises;
%   \item eareccpeak:  exclusive antirecords that are also cycle peaks;
%   \item eareccdfall:  exclusive antirecords that are also cycle double falls;
%   \item rar:  record-antirecords (these are always fixed points);
%   \item nrcpeak:  neither-record-antirecords that are also cycle peaks;
%   \item nrcval:  neither-record-antirecords that are also cycle valleys;
%   \item nrcdrise:  neither-record-antirecords that are also cycle double rises;
%   \item nrcdfall:  neither-record-antirecords that are also cycle double falls;
%   \item nrfix:  neither-record-antirecords that are also fixed points.
%\end{itemize}
Clearly every index $i$ belongs to exactly one of these 10~types;
we call this the \textbfit{record-and-cycle classification}.
The corresponding sets are $\Ereccval = \Rec \cap \Cval$, etc.

\subsection{Crossings and nestings}
     \label{subsec.statistics.2}

We now define (following \cite{Sokal-Zeng_masterpoly})
some permutation statistics that count
\textbfit{crossings} and \textbfit{nestings}.

First we associate to each permutation $\sigma \in \Sym_N$
a pictorial representation (Figure~\ref{fig.pictorial})
by placing vertices $1,2,\ldots,N$ along a horizontal axis
and then drawing an arc from $i$ to $\sigma(i)$
above (resp.\ below) the horizontal axis
in case $\sigma(i) > i$ [resp.\ $\sigma(i) < i$];
if $\sigma(i) = i$ we do not draw any arc.
\begin{figure}[t]
\centering
\vspace*{4cm}
\begin{picture}(60,20)(120, -65)
\setlength{\unitlength}{2mm}
\linethickness{.5mm}
\put(-2,0){\line(1,0){54}}
\put(0,0){\circle*{1,3}}\put(0,0){\makebox(0,-6)[c]{\small 1}}
\put(5,0){\circle*{1,3}}\put(5,0){\makebox(0,-6)[c]{\small 2}}
\put(10,0){\circle*{1,3}}\put(10,0){\makebox(0,-6)[c]{\small 3}}
\put(15,0){\circle*{1,3}}\put(15,0){\makebox(0,-6)[c]{\small 4}}
\put(20,0){\circle*{1,3}}\put(20,0){\makebox(0,-6)[c]{\small 5}}
\put(25,0){\circle*{1,3}}\put(25,0){\makebox(0,-6)[c]{\small 6}}
\put(30,0){\circle*{1,3}}\put(30,0){\makebox(0,-6)[c]{\small 7}}
\put(35,0){ \circle*{1,3}}\put(36,0){\makebox(0,-6)[c]{\small 8}}
\put(40,0){\circle*{1,3}}\put(40,0){\makebox(0,-6)[c]{\small 9}}
\put(45,0){\circle*{1,3}}\put(45,0){\makebox(0,-6)[c]{\small 10}}
\put(50,0){\circle*{1,3}}\put(50,0){\makebox(0,-6)[c]{\small 11}}
\green{\qbezier(0,0)(20,14)(40,0)
\qbezier(40,0)(42.5,6)(45,0)}
\red{\qbezier(4,0)(6.5,5)(9,0)
\qbezier(9,0)(18,10)(29,0)}
\blue{\qbezier(18,0)(20.5,5)(23.5,0)
\qbezier(23.2,0)(36,12)(48.5,0)
\qbezier(18,0)(34,-12)(48.5,0)}
%%%%%%%%%%%%%%%%%%%%%%%%%%
\red{\qbezier(2.5,0)(17,-14)(27.5,0)}
\green{\qbezier(-3,0)(22,-20)(42,0)}
\end{picture}
\caption{
   Pictorial representation of the permutation
   $\sigma = 9\,3\,7\,4\,6\,11\,2\,8\,10\,1\,5
           = (1,9,10)\,(2,3,7)\,(4)\,(5,6,11)\,(8) \in \Sym_{11}$.
 \label{fig.pictorial}
 \vspace*{8mm}
}
\end{figure}
Each vertex thus has either
out-degree = in-degree = 1 (if it is not a fixed point) or
out-degree = in-degree = 0 (if it is a fixed point).
Of course, the arrows on the arcs are redundant,
because the arrow on an arc above (resp.\ below) the axis
always points to the right (resp.\ left);
we therefore omit the arrows for simplicity.

We then say that a quadruplet $i < j < k < l$ forms an
\begin{itemize}
   \item {\em upper crossing}\/ (ucross) if $k = \sigma(i)$ and $l = \sigma(j)$;
   \item {\em lower crossing}\/ (lcross) if $i = \sigma(k)$ and $j = \sigma(l)$;
   \item {\em upper nesting}\/  (unest)  if $l = \sigma(i)$ and $k = \sigma(j)$;
   \item {\em lower nesting}\/  (lnest)  if $i = \sigma(l)$ and $j = \sigma(k)$.
\end{itemize}
We also consider some ``degenerate'' cases with $j=k$,
by saying that a triplet $i < j < l$ forms an
\begin{itemize}
   \item {\em upper joining}\/ (ujoin) if $j = \sigma(i)$ and $l = \sigma(j)$
      [i.e.\ the index $j$ is a cycle double rise];
   \item {\em lower joining}\/ (ljoin) if $i = \sigma(j)$ and $j = \sigma(l)$
      [i.e.\ the index $j$ is a cycle double fall];
   \item {\em upper pseudo-nesting}\/ (upsnest)
      if $l = \sigma(i)$ and $j = \sigma(j)$;
   \item {\em lower pseudo-nesting}\/ (lpsnest)
      if $i = \sigma(l)$ and $j = \sigma(j)$.
\end{itemize}
These are clearly degenerate cases of crossings and nestings, respectively.
See Figure~\ref{fig.crossnest}.
Note that $\upsnest(\sigma) = \lpsnest(\sigma)$ for all $\sigma$,
since for each fixed point~$j$,
the number of pairs $(i,l)$ with $i < j < l$ such that $l = \sigma(i)$
has to equal the number of such pairs with $i = \sigma(l)$;
we therefore write these two statistics simply as
\be
   \psnest(\sigma) \;\eqdef\; \upsnest(\sigma) \;=\;  \lpsnest(\sigma)
   \;.
\ee
And of course $\ujoin = \cdrise$ and $\ljoin = \cdfall$.

\begin{figure}[p]
\centering
\begin{picture}(30,15)(145, 10)
\setlength{\unitlength}{1.5mm}
\linethickness{.5mm}
\put(2,0){\line(1,0){28}}
\put(5,0){\circle*{1,3}}\put(5,0){\makebox(0,-6)[c]{\small $i$}}
\put(12,0){\circle*{1,3}}\put(12,0){\makebox(0,-6)[c]{\small $j$}}
\put(19,0){\circle*{1,3}}\put(19,0){\makebox(0,-6)[c]{\small $k$}}
\put(26,0){\circle*{1,3}}\put(26,0){\makebox(0,-6)[c]{\small $l$}}
\red{\qbezier(5,0)(12,10)(19,0)}
\blue{\qbezier(11,0)(18,10)(25,0)}
\put(15,-6){\makebox(0,-6)[c]{\small Upper crossing}}
%%%%%%%%%%%%%%%%%%%%
\put(43,0){\line(1,0){28}}
\put(47,0){\circle*{1,3}}\put(47,0){\makebox(0,-6)[c]{\small $i$}}
\put(54,0){\circle*{1,3}}\put(54,0){\makebox(0,-6)[c]{\small $j$}}
\put(61,0){\circle*{1,3}}\put(61,0){\makebox(0,-6)[c]{\small $k$}}
\put(68,0){\circle*{1,3}}\put(68,0){\makebox(0,-6)[c]{\small $l$}}
\red{\qbezier(47,0)(54,-10)(61,0)}
\blue{\qbezier(53,0)(60,-10)(67,0)}
\put(57,-6){\makebox(0,-6)[c]{\small Lower crossing}}
\end{picture}
\\[3.5cm]
\begin{picture}(30,15)(145, 10)
\setlength{\unitlength}{1.5mm}
\linethickness{.5mm}
\put(2,0){\line(1,0){28}}
\put(5,0){\circle*{1,3}}\put(5,0){\makebox(0,-6)[c]{\small $i$}}
\put(12,0){\circle*{1,3}}\put(12,0){\makebox(0,-6)[c]{\small $j$}}
\put(19,0){\circle*{1,3}}\put(19,0){\makebox(0,-6)[c]{\small $k$}}
\put(26,0){\circle*{1,3}}\put(26,0){\makebox(0,-6)[c]{\small $l$}}
\red{\qbezier(5,0)(15.5,10)(26,0)}
\blue{\qbezier(11,0)(14.5,5)(18,0)}
\put(15,-6){\makebox(0,-6)[c]{\small Upper nesting}}
%%%%%%%%%%%%%%%%%%%%
\put(43,0){\line(1,0){28}}
\put(47,0){\circle*{1,3}}\put(47,0){\makebox(0,-6)[c]{\small $i$}}
\put(54,0){\circle*{1,3}}\put(54,0){\makebox(0,-6)[c]{\small $j$}}
\put(61,0){\circle*{1,3}}\put(61,0){\makebox(0,-6)[c]{\small $k$}}
\put(68,0){\circle*{1,3}}\put(68,0){\makebox(0,-6)[c]{\small $l$}}
\red{\qbezier(47,0)(57,-10)(68,0)}
\blue{\qbezier(53,0)(56.5,-5)(60.5,0)}
\put(57,-6){\makebox(0,-6)[c]{\small Lower nesting}}
\end{picture}
\\[3.5cm]
\begin{picture}(30,15)(145, 10)
\setlength{\unitlength}{1.5mm}
\linethickness{.5mm}
\put(2,0){\line(1,0){28}}
\put(5,0){\circle*{1,3}}\put(5,0){\makebox(0,-6)[c]{\small $i$}}
\put(15.5,0){\circle*{1,3}}\put(15.5,0){\makebox(0,-6)[c]{\small $j$}}
%\put(19,0){\circle*{1,3}}\put(19,0){\makebox(0,-6)[c]%{\small $k$}}
\put(26,0){\circle*{1,3}}\put(26,0){\makebox(0,-6)[c]{\small $l$}}
\red{\qbezier(5,0)(10.5,10)(15.5,0)}
\blue{\qbezier(14.75,0)(19.5,10)(25.25,0)}
\put(15,-6){\makebox(0,-6)[c]{\small Upper joining}}
%%%%%%%%%%%%%%%%%%%%
\put(43,0){\line(1,0){28}}
\put(47,0){\circle*{1,3}}\put(47,0){\makebox(0,-6)[c]{\small $i$}}
\put(57.5,0){\circle*{1,3}}\put(57.5,0){\makebox(0,-6)[c]{\small $j$}}
%\put(51,0){\circle*{1,3}}\put(51,0){\makebox(0,-6)[c]{\small $k$}}
\put(68,0){\circle*{1,3}}\put(68,0){\makebox(0,-6)[c]{\small $l$}}
\red{\qbezier(47,0)(52,-10)(57.25,0)}
\blue{\qbezier(56.75,0)(62,-10)(67,0)}
%\blue{\qbezier(43,0)(46.5,-5)(50.5,0)}
\put(57,-6){\makebox(0,-6)[c]{\small Lower joining}}
\end{picture}
\\[3.5cm]
\begin{picture}(30,15)(145, 10)
\setlength{\unitlength}{1.5mm}
\linethickness{.5mm}
\put(2,0){\line(1,0){28}}
\put(5,0){\circle*{1,3}}\put(5,0){\makebox(0,-6)[c]{\small $i$}}
\put(15.5,0){\circle*{1,3}}\put(15.5,0){\makebox(0,-6)[c]{\small $j$}}
%\put(19,0){\circle*{1,3}}\put(19,0){\makebox(0,-6)[c]%{\small $k$}}
\put(26,0){\circle*{1,3}}\put(26,0){\makebox(0,-6)[c]{\small $l$}}
\red{\qbezier(5,0)(15.5,10)(26,0)}
%%\blue{\qbezier(14.25,0)(15,3)(14.75,0)}
\put(15,-6){\makebox(0,-6)[c]{\small Upper pseudo-nesting}}
%%%%%%%%%%%%%%%%%%%%
\put(43,0){\line(1,0){28}}
\put(47,0){\circle*{1,3}}\put(47,0){\makebox(0,-6)[c]{\small $i$}}
\put(57.5,0){\circle*{1,3}}\put(57.5,0){\makebox(0,-6)[c]{\small $j$}}
%\put(51,0){\circle*{1,3}}\put(51,0){\makebox(0,-6)[c]{\small $k$}}
\put(68,0){\circle*{1,3}}\put(68,0){\makebox(0,-6)[c]{\small $l$}}
\red{\qbezier(47,0)(57.5,-10)(68,0)}
%%\blue{\qbezier(56.75,0)(56.8,-3)(56.9,0)}
%\blue{\qbezier(53,0)(56.5,-5)(60.5,0)}
\put(57,-6){\makebox(0,-6)[c]{\small Lower pseudo-nesting}}
\end{picture}
\vspace*{3cm}
\caption{
   Crossing, nesting, joining and pseudo-nesting.
 \label{fig.crossnest}
}
\end{figure}

We can further refine the four crossing/nesting categories
by examining more closely the status of the inner index ($j$ or $k$)
whose {\em outgoing}\/ arc belongs to the crossing or nesting:
that is, $j$ for an upper crossing or nesting,
and $k$ for a lower crossing or nesting:
%we say that a quadruplet $i < j < k < l$ forms an

\medskip
\begin{center}
\begin{tabular}{c|c|c|c|c|}
 & ucross & unest & lcross & lnest \\
\hline
$j\in \Cval$ & ucrosscval & unestcval & &\\
$j\in \Cdrise$ & ucrosscdrise  & unestcdrise & &\\
$k\in \Cpeak$ & & & lcrosscpeak & lnestcpeak\\
$k\in \Cdfall$ & & & lcrosscdfall & lnestcdfall\\
\hline
\end{tabular}
\end{center}
\medskip
%\begin{itemize}
%   \item {\em upper crossing of type cval}\/ (ucrosscval)
%       if $k = \sigma(i)$ and $l = \sigma(j)$ and ${\sigma^{-1}(j) > j}$;
%   \item {\em upper crossing of type cdrise}\/ (ucrosscdrise)
%       \hbox{if $k = \sigma(i)$ and $l = \sigma(j)$ and ${\sigma^{-1}(j) < j}$;}
%   \item {\em lower crossing of type cpeak}\/ (lcrosscpeak)
%       if $i = \sigma(k)$ and $j = \sigma(l)$ and ${\sigma^{-1}(k) < k}$;
%   \item {\em lower crossing of type cdfall}\/ (lcrosscdfall)
%       if $i = \sigma(k)$ and $j = \sigma(l)$ and ${\sigma^{-1}(k) > k}$;
%   \item {\em upper nesting of type cval}\/  (unestcval)
%       if $l = \sigma(i)$ and $k = \sigma(j)$ and ${\sigma^{-1}(j) > j}$;
%   \item {\em upper nesting of type cdrise}\/  (unestcdrise)
%       if $l = \sigma(i)$ and $k = \sigma(j)$ and ${\sigma^{-1}(j) < j}$;
%   \item {\em lower nesting of type cpeak}\/  (lnestcpeak)
%       if $i = \sigma(l)$ and $j = \sigma(k)$ and ${\sigma^{-1}(k) < k}$;
%   \item {\em lower nesting of type cdfall}\/  (lnestcdfall)
%       if $i = \sigma(l)$ and $j = \sigma(k)$ and ${\sigma^{-1}(k) > k}$.
%\end{itemize}
See Figure~\ref{fig.refined_crossnest}.
Please note that for the ``upper'' quantities
the distinguished index
(i.e.\ the one for which we examine both $\sigma$ and $\sigma^{-1}$)
is in second position ($j$),
while for the ``lower'' quantities
the distinguished index is in third position ($k$).

\begin{figure}[p]
\centering
\begin{picture}(30,15)(145, 10)
\setlength{\unitlength}{1.5mm}
\linethickness{.5mm}
\put(2,0){\line(1,0){28}}
\put(5,0){\circle*{1,3}}\put(5,0){\makebox(0,-6)[c]{\small $i$}}
\put(12,0){\circle*{1,3}}\put(12,0){\makebox(0,-6)[c]{\small $j$}}
\put(19,0){\circle*{1,3}}\put(19,0){\makebox(0,-6)[c]{\small $k$}}
\put(26,0){\circle*{1,3}}\put(26,0){\makebox(0,-6)[c]{\small $l$}}
\red{\qbezier(5,0)(12,10)(19,0)}
\blue{\qbezier(11,0)(18,10)(25,0)}
\blue{\qbezier(10.5,0)(11.75,-1)(13,-2)}
\put(15,-6){\makebox(0,-6)[c]{\small Upper crossing of type cval}}
%%%%%%%%%%%%%%%%%%%%
\put(43,0){\line(1,0){28}}
\put(47,0){\circle*{1,3}}\put(47,0){\makebox(0,-6)[c]{\small $i$}}
\put(54,0){\circle*{1,3}}\put(54,0){\makebox(0,-6)[c]{\small $j$}}
\put(61,0){\circle*{1,3}}\put(61,0){\makebox(0,-6)[c]{\small $k$}}
\put(68,0){\circle*{1,3}}\put(68,0){\makebox(0,-6)[c]{\small $l$}}
\red{\qbezier(47,0)(54,10)(61,0)}
\blue{\qbezier(53,0)(60,10)(67,0)}
\blue{\qbezier(50,2)(51.25,1)(52.5,0)}
\put(57,-6){\makebox(0,-6)[c]{\small Upper crossing of type cdrise}}
\end{picture}
\\[3.5cm]
\begin{picture}(30,15)(145, 10)
\setlength{\unitlength}{1.5mm}
\linethickness{.5mm}
\put(2,0){\line(1,0){28}}
\put(5,0){\circle*{1,3}}\put(5,0){\makebox(0,-6)[c]{\small $i$}}
\put(12,0){\circle*{1,3}}\put(12,0){\makebox(0,-6)[c]{\small $j$}}
\put(19,0){\circle*{1,3}}\put(19,0){\makebox(0,-6)[c]{\small $k$}}
\put(26,0){\circle*{1,3}}\put(26,0){\makebox(0,-6)[c]{\small $l$}}
\red{\qbezier(5,0)(12,-10)(19,0)}
\red{\qbezier(16,2)(17.25, 1)(18.5,0)}
\blue{\qbezier(10.5,0)(18,-10)(25,0)}
\put(15,-6){\makebox(0,-6)[c]{\small Lower crossing of type cpeak}}
%%%%%%%%%%%%%%%%%%%%
\put(43,0){\line(1,0){28}}
\put(47,0){\circle*{1,3}}\put(47,0){\makebox(0,-6)[c]{\small $i$}}
\put(54,0){\circle*{1,3}}\put(54,0){\makebox(0,-6)[c]{\small $j$}}
\put(61,0){\circle*{1,3}}\put(61,0){\makebox(0,-6)[c]{\small $k$}}
\put(68,0){\circle*{1,3}}\put(68,0){\makebox(0,-6)[c]{\small $l$}}
\red{\qbezier(47,0)(54,-10)(61,0)}
\blue{\qbezier(53,0)(60,-10)(67,0)}
\red{\qbezier(59.5,0)(60.75,-1)(62,-2)}
\put(57,-6){\makebox(0,-6)[c]{\small Lower crossing of type cdfall}}
\end{picture}
\\[3.5cm]
\begin{picture}(30,15)(145, 10)
\setlength{\unitlength}{1.5mm}
\linethickness{.5mm}
\put(2,0){\line(1,0){28}}
\put(5,0){\circle*{1,3}}\put(5,0){\makebox(0,-5)[c]{\small $i$}}
\put(12,0){\circle*{1,3}}\put(12,0){\makebox(0,-5)[c]{\small $j$}}
\put(19,0){\circle*{1,3}}\put(19,0){\makebox(0,-5)[c]{\small $k$}}
\put(26,0){\circle*{1,3}}\put(26,0){\makebox(0,-5)[c]{\small $l$}}
\red{\qbezier(5,0)(15.5,10)(26,0)}
\blue{\qbezier(11,0)(14.5,5)(18,0)}
\blue{\qbezier(10.25,0)(11.5,-1)(12.75,-2)}
\put(15,-6){\makebox(0,-6)[c]{\small Upper nesting of type cval}}
%%%%%%%%%%%%%%%%%%%%
\put(43,0){\line(1,0){28}}
\put(47,0){\circle*{1,3}}\put(47,0){\makebox(0,-6)[c]{\small $i$}}
\put(54,0){\circle*{1,3}}\put(54,0){\makebox(0,-6)[c]{\small $j$}}
\put(61,0){\circle*{1,3}}\put(61,0){\makebox(0,-6)[c]{\small $k$}}
\put(68,0){\circle*{1,3}}\put(68,0){\makebox(0,-6)[c]{\small $l$}}
\red{\qbezier(47,0)(57,10)(68,0)}
\blue{\qbezier(53,0)(56.5,5)(60.5,0)}
\blue{\qbezier(50,2)(51.25,1)(52.5,0)}
\put(57,-6){\makebox(0,-6)[c]{\small Upper nesting of type cdrise}}
\end{picture}
\\[3.5cm]
\begin{picture}(30,15)(145, 10)
\setlength{\unitlength}{1.5mm}
\linethickness{.5mm}
\put(2,0){\line(1,0){28}}
\put(5,0){\circle*{1,3}}\put(5,0){\makebox(0,-5)[c]{\small $i$}}
\put(12,0){\circle*{1,3}}\put(12,0){\makebox(0,-5)[c]{\small $j$}}
\put(19,0){\circle*{1,3}}\put(19,0){\makebox(0,-5)[c]{\small $k$}}
\put(26,0){\circle*{1,3}}\put(26,0){\makebox(0,-5)[c]{\small $l$}}
\red{\qbezier(5,0)(15.5,-10)(26,0)}
\blue{\qbezier(11,0)(14.5,-5)(18,0)}
\blue{\qbezier(14.75,2)(16,1)(17.25,0)}
\put(15,-6){\makebox(0,-6)[c]{\small Lower nesting of type cpeak}}
%%%%%%%%%%%%%%%%%%%%
\put(43,0){\line(1,0){28}}
\put(47,0){\circle*{1,3}}\put(47,0){\makebox(0,-5)[c]{\small $i$}}
\put(54,0){\circle*{1,3}}\put(54,0){\makebox(0,-5)[c]{\small $j$}}
\put(61,0){\circle*{1,3}}\put(61,0){\makebox(0,-5)[c]{\small $k$}}
\put(68,0){\circle*{1,3}}\put(68,0){\makebox(0,-5)[c]{\small $l$}}
\red{\qbezier(47,0)(57,-10)(68,0)}
\blue{\qbezier(53,0)(56.5,-5)(60.5,0)}
\blue{\qbezier(59.5,0)(60.75, -1)(62,-2)}
\put(57,-6){\makebox(0,-6)[c]{\small Lower nesting of type cdfall}}
\end{picture}
\vspace*{3cm}
\caption{
   Refined categories of crossing and nesting.
 \label{fig.refined_crossnest}
}
\end{figure}

In fact, a central role in our work will be played
(just as in \cite{Sokal-Zeng_masterpoly,Deb-Sokal_genocchi})
by a refinement of these statistics:
rather than counting the {\em total}\/ numbers of quadruplets
$i < j < k < l$ that form upper (resp.~lower) crossings or nestings,
we will count the number of upper (resp.~lower) crossings or nestings
that use a particular vertex $j$ (resp.~$k$)
in second (resp.~third) position.
More precisely, we define the
\textbfit{index-refined crossing and nesting statistics}
\begin{subeqnarray}
   \ucross(j,\sigma)
   & = &
   \#\{ i<j<k<l \colon\: k = \sigma(i) \hbox{ and } l = \sigma(j) \}
         \\[2mm]
   \unest(j,\sigma)
   & = &
   \#\{ i<j<k<l \colon\: k = \sigma(j) \hbox{ and } l = \sigma(i) \}
      \\[2mm]
   \lcross(k,\sigma)
   & = &
   \#\{ i<j<k<l \colon\: i = \sigma(k) \hbox{ and } j = \sigma(l) \}
         \\[2mm]
   \lnest(k,\sigma)
   & = &
   \#\{ i<j<k<l \colon\: i = \sigma(l) \hbox{ and } j = \sigma(k) \}
 \label{def.ucrossnestjk}
\end{subeqnarray}
%
% {\bf But maybe we want to define these with second and third positions
%    reversed, i.e.
% \begin{subeqnarray}
%    \ucross^\star(k,\sigma)
%    & = &
%    \#\{ i<j<k<l \colon\: k = \sigma(i) \hbox{ and } l = \sigma(j) \}
%          \\[2mm]
%    \unest^\star(k,\sigma)
%    & = &
%    \#\{ i<j<k<l \colon\: k = \sigma(j) \hbox{ and } l = \sigma(i) \}
%       \\[2mm]
%    \lcross^\star(j,\sigma)
%    & = &
%    \#\{ i<j<k<l \colon\: i = \sigma(k) \hbox{ and } j = \sigma(l) \}
%          \\[2mm]
%    \lnest^\star(j,\sigma)
%    & = &
%    \#\{ i<j<k<l \colon\: i = \sigma(l) \hbox{ and } j = \sigma(k) \}
%  \label{def.ucrossnestjk.star}
% \end{subeqnarray}
% This would not only remove the $\sigma^{-1}$ in \reff{eq.xii.nestings},
% but even more importantly it would allow the subsequent equations
% to be interpreted in terms of $\ucross^\star$ and $\lcross^\star$.}
%
Note that $\ucross(j,\sigma)$ and $\unest(j,\sigma)$ can be nonzero
only when $j$ is an excedance
(that is, a cycle valley or a cycle double rise),
while $\lcross(k,\sigma)$ and $\lnest(k,\sigma)$ can be nonzero
only when $k$ is an anti-excedance
(that is, a cycle peak or a cycle double fall).

When $j$ is a fixed point, we also define the analogous quantity
for pseudo-nestings:
\be
   \psnest(j,\sigma)
   \;\eqdef\;
   \# \{i < j \colon\:  \sigma(i) > j \}
   \;=\;
   \# \{i > j \colon\:  \sigma(i) < j \}
   \;.
 \label{def.psnestj}
\ee
(Here the two expressions are equal because $\sigma$ is a bijection
 from $[1,j) \cup (j,n]$ to itself.)
In \cite[eq.~(2.20)]{Sokal-Zeng_masterpoly}
this quantity was called the {\em level}\/ of the fixed point $j$
and was denoted $\lev(j,\sigma)$.

\section{First S-fraction for cycle-alternating permutations} \label{sec.S-fraction1}

We recall that a permutation $\sigma$ is called {\em cycle-alternating}\/
if it has no cycle double rises, cycle double falls, or fixed points;
thus, each cycle of $\sigma$ is of even length (call it~$2m$)
and consists of $m$ cycle valleys and $m$ cycle peaks in alternation.
We write $\Sym^{\rm ca}_{2n}$ for the set of cycle-alternating permutations
of $[2n]$.

\subsection{First S-fraction (8 variables)}

In \cite[Section~2.15]{Sokal-Zeng_masterpoly}
the following 4-variable polynomials were introduced:
\be
   Q_n(x,y,u,v)
   \;=\;
   \sum_{\sigma \in \Sym^{\rm ca}_{2n}}
   x^{\eareccpeak(\sigma)}
   y^{\ereccval(\sigma)}
   u^{\nrcpeak(\sigma)}
   v^{\nrcval(\sigma)}
   \;,
 \label{def.Qn.cycle-alternating}
\ee
where
\be
   \eareccpeak(\sigma)
   \;=\;
   |\Eareccpeak(\sigma)|
   \;=\;
   |\Arec(\sigma) \cap \Cpeak(\sigma)|
\ee
and likewise for the others.
It was then shown \cite[Theorem~2.18]{Sokal-Zeng_masterpoly}
that the ordinary generating function of the polynomials $Q_n$
has the S-type continued fraction
\be
   \sum_{n=0}^\infty Q_n(x,y,u,v) \: t^n
   \;=\;
   \cfrac{1}{1 - \cfrac{x y t}{1 -  \cfrac{(x\!+\!u)(y\!+\!v) t}{1 - \cfrac{(x\!+\!2u)(y\!+\!2v) t}{1 - \cdots}}}}
   \label{eq.thm.perm.Stype.cycle-alternating}
\ee
with coefficients
\be
   \alpha_n   \;=\;   [x + (n-1)u] \: [y + (n-1)v]
   \;.
 \label{def.weights.perm.Stype.cycle-alternating}
\ee

We now generalize this by treating even and odd indices separately.
That is, we define the polynomials
\begin{eqnarray}
   & & \hspace*{-6mm}
   Q_n(\xe,\ye,\ue,\ve,\xo,\yo,\uo,\vo)
   \;=\;
        \nonumber \\[4mm]
   & & \qquad
   \sum_{\sigma \in \Sym^{\rm ca}_{2n}}
   \xe^{\eareccpeakeven(\sigma)}
   \ye^{\ereccvaleven(\sigma)}
   \ue^{\nrcpeakeven(\sigma)}
   \ve^{\nrcvaleven(\sigma)}
          \:\times
       \qquad\qquad
       \nonumber \\[0mm]
   & & \qquad\qquad\;\,
   \xo^{\eareccpeakodd(\sigma)}
   \yo^{\ereccvalodd(\sigma)}
   \uo^{\nrcpeakodd(\sigma)}
   \vo^{\nrcvalodd(\sigma)}
   \;,
 \label{def.Qn.cycle-alternating.evenodd}
\end{eqnarray}
where
\be
   \eareccpeakeven(\sigma)
   \;=\;
   |\Eareccpeakeven(\sigma)|
   \;=\;
   |\Arec(\sigma) \cap \Cpeak(\sigma) \cap \Even|
\ee
and likewise for the others.
We then have:

\begin{theorem}[First S-fraction for cycle-alternating permutations]
   \label{thm.first.Sfrac}
The ordinary generating function of the polynomials
$Q_n(\xe,\ye,\ue,\ve,\xo,\yo,\uo,\vo)$ has the S-type continued fraction
\be
   \sum_{n=0}^\infty Q_n(\xe,\ye,\ue,\ve,\xo,\yo,\uo,\vo) \: t^n
   \;=\;
   \cfrac{1}{1 - \cfrac{\xe \yo t}{1 -  \cfrac{(\xo\!+\!\uo)(\ye\!+\!\ve) t}{1 - \cfrac{(\xe\!+\!2\ue)(\yo\!+\!2\vo) t}{1 - \cdots}}}}
   \label{eq.thm.perm.Stype.cycle-alternating.evenodd}
\ee
with coefficients
\begin{subeqnarray}
   \alpha_{2k-1}  & = &  [\xe + (2k-2) \ue] \: [\yo + (2k-2) \vo] \\[1mm]
   \alpha_{2k}    & = &  [\xo + (2k-1) \uo] \: [\ye + (2k-1) \ve]
 \label{def.weights.perm.Stype.cycle-alternating.evenodd}
\end{subeqnarray}
\end{theorem}

The proof of this theorem is based on
the ``master'' S-fraction for cycle-alternating permutations
\cite[Theorem~2.20]{Sokal-Zeng_masterpoly}
together with the following lemma:

\begin{lemma}[Key lemma]
   \label{lemma.cycle-alt}
If $\sigma$ is a cycle-alternating permutation of $[2n]$, then
\begin{subeqnarray}
   \hbox{\rm cycle valleys:}
   & \!\!  &
   \ucross(i,\sigma) \,+\, \unest(i,\sigma)
   \;\equiv\;
   i-1 \: \pmod{2}
         \\[1mm]
   \hbox{\rm cycle peaks:}
   & \!\!  &
   \lcross(i,\sigma) \,+\, \lnest(i,\sigma)
   \;\equiv\;
   i \: \pmod{2}
\end{subeqnarray}
for all $i \in [2n]$.
% {\bf Is this right??????
%    This is what I get from \cite[eqns.~(6.6)--(6.9)]{Sokal-Zeng_masterpoly}
%    and the definition of the steps $s_i = h_i - h_{i-1}$.}
\end{lemma}

%% \firstproof
%% {\bf Do it directly on the permutation!!!}
%% \qed
%% 
%% \secondproof
\proof
In the Foata--Zeilberger bijection for permutations
as presented in \cite[Section~6.1]{Sokal-Zeng_masterpoly},
permutations $\sigma \in \Sym_N$ are mapped bijectively onto
pairs $(\omega,\xi)$,
where $\omega = (\omega_0,\ldots,\omega_N)$ is a Motzkin path of length $N$,
and $\xi$ is a suitably defined set of integer labels $(\xi_1,\ldots,\xi_N)$.
We write $h_i$ for the height of the Motzkin path after step~$i$,
i.e.\ $\omega_i = (i, h_i)$.
{}From \cite[eqns.~(6.6)--(6.9)]{Sokal-Zeng_masterpoly}
and the definition of the steps $s_i = h_i - h_{i-1}$,
we get
\begin{subeqnarray}
   \textrm{cval:}
   & &
   \ucross(i,\sigma) \,+\, \unest(i,\sigma)
   \;=\;
   h_{i-1}
   \;=\;
   h_i - 1
      \\[1mm]
   \textrm{cdrise:}
   & &
   \ucross(i,\sigma) \,+\, \unest(i,\sigma)
   \;=\;
   h_{i-1} - 1
   \;=\;
   h_i - 1
      \\[1mm]
   \textrm{cpeak:}
   & &
   \lcross(i,\sigma) \,+\, \lnest(i,\sigma)
   \;=\;
   h_{i-1} - 1
   \;=\;
   h_i
      \\[1mm]
   \textrm{cdfall:}
   & &
   \lcross(i,\sigma) \,+\, \lnest(i,\sigma)
   \;=\;
   h_{i-1} - 1
   \;=\;
   h_i - 1
\end{subeqnarray}
On the other hand, for a cycle-alternating permutation,
the path $\omega$ is a Dyck path (i.e.~there are no level steps),
so that $h_i \equiv i \pmod{2}$.
\qed

We now apply the first master S-fraction for cycle-alternating permutations
\cite[Theorem~2.20]{Sokal-Zeng_masterpoly},
which states that the polynomials
\be
   Q_n(\bsfa,\bsfb)
   \;=\;
   \sum_{\sigma \in \Sym^{\rm ca}_{2n}}
   \;\:
   \prod\limits_{i \in \Cval}  \! \sfa_{\ucross(i,\sigma),\,\unest(i,\sigma)}
   \prod\limits_{i \in \Cpeak} \!\! \sfb_{\lcross(i,\sigma),\,\lnest(i,\sigma)}
 \label{def.Qn.firstmaster.cycle-alternating}
\ee
in indeterminates
$\bsfa = (\sfa_{\ell,\ell'})_{\ell,\ell' \ge 0}$
and $\bsfb = (\sfb_{\ell,\ell'})_{\ell,\ell' \ge 0}$
have the S-type continued fraction
\be
   \sum_{n=0}^\infty Q_n(\bsfa,\bsfb)  \: t^n
   \;=\;
   \cfrac{1}{1 - \cfrac{\sfa_{00} \sfb_{00} t}{1 - \cfrac{(\sfa_{01} + \sfa_{10})(\sfb_{01} + \sfb_{10}) t}{1 - \cfrac{(\sfa_{02} + \sfa_{11} + \sfa_{20})(\sfb_{02} + \sfb_{11} + \sfb_{20}) t}{1 - \cdots}}}}
   \label{eq.thm.permutations.Stype.final1.cycle-alternating}
\ee
with coefficients
\be
   \alpha_n   \;=\;
   \biggl( \sum_{\ell=0}^{n-1} \sfa_{\ell,n-1-\ell} \biggr)
   \biggl( \sum_{\ell=0}^{n-1} \sfb_{\ell,n-1-\ell} \biggr)
   \;.
   \label{eq.thm.permutations.Stype.final1.cycle-alternating.weights}
\ee
We also use the following general fact about permutations
%\cite[Lemma~2.10]{Sokal-Zeng_masterpoly}:
\begin{lemma} {\hspace*{-2mm}\bf \cite[Lemma~2.10]{Sokal-Zeng_masterpoly}\ }
   \label{lemma.SZ.nesting}
Let $\sigma$ be a permutation.
\vspace*{-2mm}
\begin{itemize}
   \item[(a)]  If $i$ is a cycle valley or cycle double rise,
       then $i$ is a record if and only if $\unest(i,\sigma) = 0$;
       and in this case it is an exclusive record.
   \item[(b)]  If $i$ is a cycle peak or cycle double fall,
       then $i$ is an antirecord if and only if $\lnest(i,\sigma) = 0$;
       and in this case it is an exclusive antirecord.
\end{itemize}
\end{lemma}
%\begin{itemize}
%   \item[(a)]  If $i$ is a cycle valley or cycle double rise,
%       then $i$ is a record if and only if $\unest(i,\sigma) = 0$;
%       and in this case it is an exclusive record.
%   \item[(b)]  If $i$ is a cycle peak or cycle double fall,
%       then $i$ is an antirecord if and only if $\lnest(i,\sigma) = 0$;
%       and in this case it is an exclusive antirecord.
%\end{itemize}
Applying Lemma~\ref{lemma.SZ.nesting} to cycle valleys 
and peaks in a cycle-alternating permutation,
it follows that the polynomials $Q_n(\xe,\ye,\ue,\ve,\xo,\yo,\uo,\vo)$
are obtained by specializing $Q_n(\bsfa,\bsfb)$ to
\begin{subeqnarray}
   \sfa_{\ell,\ell'}
   & = &
   \begin{cases}
       \yo    &  \textrm{if $\ell' = 0$ and $\ell+\ell'$ is even}  \\
       \vo    &  \textrm{if $\ell' \ge 1$ and $\ell+\ell'$ is even} \\
       \ye    &  \textrm{if $\ell' = 0$ and $\ell+\ell'$ is odd}  \\
       \ve    &  \textrm{if $\ell' \ge 1$ and $\ell+\ell'$ is odd}
   \end{cases}
       \\[2mm]
   \sfb_{\ell,\ell'}
   & = &
   \begin{cases}
       \xe    &  \textrm{if $\ell' = 0$ and $\ell+\ell'$ is even}  \\
       \ue    &  \textrm{if $\ell' \ge 1$ and $\ell+\ell'$ is even} \\
       \xo    &  \textrm{if $\ell' = 0$ and $\ell+\ell'$ is odd}  \\
       \uo    &  \textrm{if $\ell' \ge 1$ and $\ell+\ell'$ is odd}
   \end{cases}
\end{subeqnarray}
Inserting these into
\reff{eq.thm.permutations.Stype.final1.cycle-alternating.weights}
yields
\reff{def.weights.perm.Stype.cycle-alternating.evenodd}.
This completes the proof of Theorem~\ref{thm.first.Sfrac}.

\bigskip

{\bf Remarks.}
1. If we specialize Theorem~\ref{thm.first.Sfrac}
to $\ue = \xe$, $\uo = \xo$, $\ye = \ve = \yo = \vo = 1$,
we obtain an S-fraction with coefficients
$\alpha_{2k-1} = (2k-1)^2 \xe$, $\alpha_{2k} = (2k)^2 \xo$.
This is the S-fraction for the polynomials $\scre_{2n}(\xe,\xo)$,
which generalize the secant numbers $E_{2n}$ as follows:
Start from the Jacobian elliptic functions $\sn$, $\cn$, $\dn$
in Dumont's \cite{Dumont_81a} symmetric parametrization, which are
defined by the system of differential equations
\begin{subeqnarray}
   {d \over du} \, \sn(u;a,b) & = & \cn(u;a,b) \: \dn(u;a,b)   \\[2mm]
   {d \over du} \, \cn(u;a,b) & = & a^2 \, \dn(u;a,b) \: \sn(u;a,b)   \\[2mm]
   {d \over du} \, \dn(u;a,b) & = & b^2 \, \sn(u;a,b) \: \cn(u;a,b)
 \label{eq.system.elliptic}
\end{subeqnarray}
with the initial conditions
$\sn(0;a,b) = 0$, $\cn(0;a,b) = 1$, $\dn(0;a,b) = 1$.
Note the special cases
\begin{subeqnarray}
   \sn(u;a,0)  & = &  {\sinh au \over a}  \\[2mm]
   \cn(u;a,0)  & = &  \cosh au            \\[2mm]
   \dn(u;a,0)  & = &  1
\end{subeqnarray}
and
\begin{subeqnarray}
   \sn(u;a,a)  & = &  {\tan au \over a}  \\[2mm]
   \cn(u;a,a)  & = &  \sec au  \\[2mm]
   \dn(u;a,a)  & = &  \sec au
  \label{eq.jacobian.tansec}
\end{subeqnarray}
Thus, $\sn$ is a simultaneous generalization of tangent and hyperbolic sine,
while $\cn$ is a simultaneous generalization of secant and hyperbolic cosine.
The Taylor expansions around $u=0$ have the form
\begin{eqnarray}
   \cn(u;a,b)
   & = &
   \sum_{n=0}^\infty \scre_{2n}(a^2,b^2) \, {u^{2n} \over (2n)!}
       \label{eq.cn.taylor}  \\[2mm]
   \sn(u;a,b)
   & = &
   \sum_{n=0}^\infty \scre_{2n+1}(a^2,b^2) \, {u^{2n+1} \over (2n+1)!}
       \label{eq.sn.taylor}
\end{eqnarray}
where $\scre_{2n}(\alpha,\beta)$ and $\scre_{2n+1}(\alpha,\beta)$
are homogeneous polynomials of degree $n$
in the indeterminates $\alpha$ and $\beta$,
with nonnegative integer coefficients.
Then the ordinary generating function of the polynomials $\scre_{2n}$
can be expressed as an S-fraction
\be
   \sum_{n=0}^\infty t^n \, \scre_{2n}(\alpha,\beta)
   \;=\;
   \cfrac{1}{1 - \cfrac{1^2 \alpha \, t}{1 - \cfrac{2^2 \beta \, t}{1 -  \cfrac{3^2 \alpha \, t}{1- \cdots}}}}
 \label{eq.scre2n.contfrac}
\ee
with coefficients
\begin{subeqnarray}
   \alpha_{2k-1}  & = &  (2k-1)^2 \, \alpha \\[2mm]
   \alpha_{2k}    & = &  (2k)^2 \, \beta
 \label{def.weights.scre2n}
\end{subeqnarray}
This S-fraction goes back to
Stieltjes \cite[p.~H17]{Stieltjes_1889} and Rogers \cite[p.~77]{Rogers_07};
modern presentations of Rogers' elegant proof
can be found in Flajolet and Fran\c{c}on \cite[Theorem~1]{Flajolet_89}
and Milne \cite[Theorems~3.2 and 3.11]{Milne_02}.
The polynomials $\scre_{2n}(\alpha,\beta)$
are also specializations of the three-variable Schett polynomials
\cite{Schett_76,Schett_77,Dumont_79,Dumont_81a,Dumont_86}
$X_{2n}(x,y,z)$,
for which Dumont \cite{Dumont_79,Dumont_81a,Dumont_86}
has provided a combinatorial interpretation in terms of
even and odd cycle peaks in permutations of $[2n]$.
When specialized to $z=0$,
one obtains cycle-alternating permutations,
and Dumont's interpretation agrees with ours.
The continued fraction \reff{eq.scre2n.contfrac}
can be found on \cite[p.~38]{Dumont_86}.

2. The continued fraction of Theorem 4.1, specialized to
$x_e = u_e = v_e = x_o = u_o = v_0 = 1$,
enumerates cycle-alternating permutations of $[2n]$ with respect to
odd and even records, and has coefficients
$\alpha_{2k-1} = (2k-1)(2k-2 + \yo)$, $\alpha_{2k} = (2k)(2k-1+\ye)$.
This same continued fraction was found by
Randrianarivony and Zeng \cite[Th\'eor\`eme~3]{Randrianarivony_94b}
as enumerating even and odd record {\em values}\/
in alternating (not cycle-alternating) permutations.
This is a ``linear statistics" analogue of our ``cyclic statistics" results.

\subsection[$p,q$-generalization of first S-fraction (16 variables)]{$\bm{p,q}$-generalization of first S-fraction (16 variables)}

We can further extend Theorem~\ref{thm.first.Sfrac}
by introducing a $p,q$-generalization.
The statistics on permutations corresponding to the variables $p$ and $q$
will be crossings and nestings,
as defined in Section~\ref{subsec.statistics.2}.

Recall first that for an integer $n \ge 0$ we define
\be
   [n]_{p,q}
   \;=\;
   {p^n - q^n \over p-q}
   \;=\;
   \sum\limits_{j=0}^{n-1} p^j q^{n-1-j}
\ee
where $p$ and $q$ are indeterminates;
it is a homogeneous polynomial of degree $n-1$ in $p$ and $q$,
which is symmetric in $p$ and $q$.
In particular, $[0]_{p,q} = 0$ and $[1]_{p,q} = 1$;
and for $n \ge 1$ we have the recurrence
\be
   [n]_{p,q}
   \;=\;
   p \, [n-1]_{p,q} \,+\, q^{n-1}
   \;=\;
   q \, [n-1]_{p,q} \,+\, p^{n-1}
   \;.
 \label{eq.recurrence.npq}
\ee
If $p=1$, then $[n]_{1,q}$ is the well-known $q$-integer
\be
   [n]_q
   \;=\;  [n]_{1,q}
   \;=\; {1 - q^n \over 1-q}
   \;=\;  \begin{cases}
               0  & \textrm{if $n=0$}  \\
               1+q+q^2+\ldots+q^{n-1}  & \textrm{if $n \ge 1$}
          \end{cases}
\ee
If $p=0$, then
\be
   [n]_{0,q}
   \;=\;  \begin{cases}
               0  & \textrm{if $n=0$}  \\
               q^{n-1}  & \textrm{if $n \ge 1$}
          \end{cases}
\ee

In \cite[Section~2.15]{Sokal-Zeng_masterpoly}
the following 8-variable polynomials were introduced:
\begin{eqnarray}
   & &
   Q_{n}(x,y,u,v,p_{+},p_{-},q_{+},q_{-})
   \;=\;
       \nonumber \\[4mm]
   & & \qquad\qquad
   \sum_{\sigma \in \Sym^{\rm ca}_{2n}}
   x^{\eareccpeak(\sigma)}
   y^{\ereccval(\sigma)}
   u^{\nrcpeak(\sigma)}
   v^{\nrcval(\sigma)}
       \:\times
       \qquad\qquad
       \nonumber \\[0mm]
   & & \qquad\qquad\qquad\;\:
   p_{-}^{\lcrosscpeak(\sigma)}
   p_{+}^{\ucrosscval(\sigma)}
   q_{-}^{\lnestcpeak(\sigma)}
   q_{+}^{\unestcval(\sigma)}
 \label{def.Qn.cycle-alternating.pqgen}
\end{eqnarray}
%\be
%   \eareccpeak(\sigma)
%   \;=\;
%   |\Eareccpeak(\sigma)|
%   \;=\;
%   |\Arec(\sigma) \cap \Cpeak(\sigma)|
%\ee
It was then shown \cite[Theorem~2.19]{Sokal-Zeng_masterpoly}
that the ordinary generating function of the polynomials $Q_n$
has the S-type continued fraction
\begin{eqnarray}
   & & \hspace*{-2cm}
   \sum_{n=0}^\infty
   Q_{n}(x,y,u,v,p_{+},p_{-},q_{+},q_{-}) \; t^n
   \;=\;
       \nonumber \\
   & & \!\!\!\!
   \cfrac{1}{1 - \cfrac{x y t}{1 - \cfrac{(p_{-} x\!+\! q_{-} u)(p_{+} y\!+\! q_{+} v) t}{1 - \cfrac{(p_{-}^2 x\!+\! q_{-} [2]_{p_{-},q_{-}} u)(p_{+}^2 y\!+\! q_{+} [2]_{p_{+},q_{+}} v) t}{1 - \cdots}}}}
   \label{eq.thm.perm.Stype.cycle-alternating.pqgen}
\end{eqnarray}
with coefficients
\be
   \alpha_n   \;=\;   (p_{-}^{n-1} x + q_{-} \, [n-1]_{p_{-},q_{-}} u)
                  \: (p_{+}^{n-1} y + q_{+} \, [n-1]_{p_{+},q_{+}} v)
   \;.
 \label{def.weights.perm.Stype.cycle-alternating.pqgen}
\ee

We now generalize this by treating even and odd indices separately.
That is, we define the polynomials
\begin{eqnarray}
   & & \hspace*{-6mm}
	Q_n(\xe,\ye,\ue,\ve,\xo,\yo,\uo,\vo,p_{-1},p_{-2},p_{+1},p_{+2},q_{-1},q_{-2},q_{+1},q_{+2})
   \;=\;
        \nonumber \\[4mm]
   & & \qquad
   \sum_{\sigma \in \Sym^{\rm ca}_{2n}}
   \xe^{\eareccpeakeven(\sigma)}
   \ye^{\ereccvaleven(\sigma)}
   \ue^{\nrcpeakeven(\sigma)}
   \ve^{\nrcvaleven(\sigma)}
          \:\times
       \qquad\qquad
       \nonumber \\[-2mm]
   & & \qquad\qquad\;\,
   \xo^{\eareccpeakodd(\sigma)}
   \yo^{\ereccvalodd(\sigma)}
   \uo^{\nrcpeakodd(\sigma)}
   \vo^{\nrcvalodd(\sigma)}
	\:\times
       \qquad\qquad
       \nonumber \\[2mm]
   & & \qquad\qquad\;\,
   p_{-1}^{{\rm lcrosscpeakeven}(\sigma)}
   p_{-2}^{{\rm lcrosscpeakodd}(\sigma)}
   p_{+1}^{{\rm ucrosscvalodd}(\sigma)}
   p_{+2}^{{\rm ucrosscvaleven}(\sigma)}
	 \:\times
       \qquad\qquad
       \nonumber \\[2mm]
   & & \qquad\qquad\;\,
   q_{-1}^{{\rm lnestcpeakeven}(\sigma)}
   q_{-2}^{{\rm lnestcpeakodd}(\sigma)}
   q_{+1}^{{\rm unestcvalodd}(\sigma)}
   q_{+2}^{{\rm unestcvaleven}(\sigma)}
   \;,
\label{def.Qn.cycle-alternating.evenodd.pqgen}
\end{eqnarray}
where
\be
   {\rm lcrosscpeakeven}(\sigma)
   \;=\;
   \sum_{j\in \Cpeak(\sigma) \cap \Even} \lcross(j,\sigma)
\ee
and likewise for the others.
We then have:

\begin{theorem}[First S-fraction for cycle-alternating permutations, $p,q$-generalization]
   \label{thm.first.Sfrac.pqgen}
The ordinary generating function of the
polynomials~\reff{def.Qn.cycle-alternating.evenodd.pqgen}
%{\hbox{$Q_n(\xe,\ye,\ue,\ve,\xo,\yo,\uo,\vo,p_{-1},p_{-2},p_{+1},p_{+2},q_{-1},q_{-2},q_{+1},q_{+2})$}}
has the S-type continued fraction
\begin{subeqnarray}
   \sum_{n=0}^\infty Q_n(\xe,\ye,\ue,\ve,\xo,\yo,\uo,\vo,p_{-1},p_{-2},p_{+1},p_{+2},q_{-1},q_{-2},q_{+1},q_{+2}) \: t^n\nonumber\\
   \;=\;
	\cfrac{1}{1 - \cfrac{\xe \yo t}{1 -  \cfrac{(p_{-2}\xo\!+\!q_{-2}\uo)(p_{+2}\ye\!+\!q_{+2}\ve) t}{1 - \cfrac{(p_{-1}^2\xe\!+\!q_{-1}[2]_{p_{-1},q_{-1}}\ue)(p_{+1}^2\yo\!+\!q_{+1}[2]_{p_{+1},q_{+1}}\vo) t}{1 - \cdots}}}}
   \label{eq.thm.perm.Stype.cycle-alternating.evenodd.pqgen}
\end{subeqnarray}
with coefficients
\begin{subeqnarray}
   \alpha_{2k-1}  & = &  (p_{-1}^{2k-2}\xe + q_{-1}[2k-2]_{p_{-1},q_{-1}} \ue) \: (p_{+1}^{2k-2}\yo + q_{+1}[2k-2]_{p_{+1},q_{+1}} \vo) \nonumber\\
	\mbox{}\\
   \alpha_{2k}    & = &  (p_{-2}^{2k-1}\xo + q_{-2}[2k-1]_{p_{-2},q_{-2}} \uo) \: (p_{+2}^{2k-1}\ye + q_{+2}[2k-1]_{p_{+2},q_{+2}} \ve)\nonumber\\
	\mbox{}
 \label{def.weights.perm.Stype.cycle-alternating.evenodd.pqgen}
\end{subeqnarray}
\end{theorem}

\proof Notice that we obtain the polynomials
\reff{def.Qn.cycle-alternating.evenodd.pqgen}
%% $Q_n(\xe,\ye,\ue,\ve,\xo,\yo,\uo,\vo,p_{-1},p_{-2},p_{+1},p_{+2},q_{-1},q_{-2},q_{+1},q_{+2})$
by specializing the polynomials 
$Q_n(\bsfa,\bsfb)$ defined in
\reff{def.Qn.firstmaster.cycle-alternating}
to 
\begin{subeqnarray}
   \sfa_{\ell,\ell'}
   & = &
   \begin{cases}
       p_{+1}^{\ell} \yo    &  \textrm{if $\ell' = 0$ and $\ell+\ell'$ is even}  \\
       p_{+1}^{\ell} q_{+1}^{\ell'}\vo    &  \textrm{if $\ell' \ge 1$ and $\ell+\ell'$ is even} \\
       p_{+2}^{\ell} \ye    &  \textrm{if $\ell' = 0$ and $\ell+\ell'$ is odd}  \\
       p_{+2}^{\ell} q_{+2}^{\ell'}\ve    &  \textrm{if $\ell' \ge 1$ and $\ell+\ell'$ is odd}
   \end{cases}
       \\[2mm]
   \sfb_{\ell,\ell'}
   & = &
   \begin{cases}
       p_{-1}^{\ell} \xe    &  \textrm{if $\ell' = 0$ and $\ell+\ell'$ is even}  \\
       p_{-1}^{\ell} q_{-1}^{\ell'}\ue    &  \textrm{if $\ell' \ge 1$ and $\ell+\ell'$ is even} \\
       p_{-2}^{\ell} \xo    &  \textrm{if $\ell' = 0$ and $\ell+\ell'$ is odd}  \\
       p_{-2}^{\ell} q_{-2}^{\ell'}\uo    &  \textrm{if $\ell' \ge 1$ and $\ell+\ell'$ is odd}
   \end{cases}
\end{subeqnarray}
Inserting these into
\reff{eq.thm.permutations.Stype.final1.cycle-alternating.weights}
yields
\reff{def.weights.perm.Stype.cycle-alternating.evenodd.pqgen}.
This along with Lemma~\ref{lemma.cycle-alt}
completes the proof of Theorem~\ref{thm.first.Sfrac.pqgen}.
\qed

{\bf Remark.}  We can also obtain Biane's \cite[section~6]{Biane_93}
S-fraction for the $q$-secant numbers $E_{2n}(q)$, defined as counting
cycle-alternating permutations by number of inversions.
We recall \cite[Proposition~2.24]{Sokal-Zeng_masterpoly}
that the number of inversions in a permutation $\sigma$ satisfies
\be
   \inv
   \;=\;
   \cval + \cdrise + \cdfall + \ucross + \lcross
                    + 2(\unest + \lnest + \psnest)
   \;.
\ee
For a cycle-alternating permutation we have
$\cdrise = \cdfall = \psnest = 0$.
Evaluating \reff{def.Qn.cycle-alternating.pqgen}
with $x=u=1$, $y=v=q$, $p_- = p_+ = q$ and $q_- = q_+ = q^2$
and using
\reff{eq.thm.perm.Stype.cycle-alternating.pqgen}/%
\reff{def.weights.perm.Stype.cycle-alternating.pqgen},
we obtain after a bit of algebra $\alpha_n = q^{2n-1} [n]_q^2$,
which agrees with Biane's formula.
\myendremark

\section{Second S-fraction for cycle-alternating permutations} \label{sec.S-fraction2}

\subsection{Second S-fraction (counting of cycles)}

It is now natural to want to generalize the foregoing polynomials
by introducing also the counting of cycles.
That is, we define the polynomials
\begin{eqnarray}
   & & \hspace*{-6mm}
   \widehat{Q}_n(\xe,\ye,\ue,\ve,\xo,\yo,\uo,\vo,\lambda)
   \;=\;
        \nonumber \\[3mm]
   & & \qquad
   \sum_{\sigma \in \Sym^{\rm ca}_{2n}}
   \xe^{\eareccpeakeven(\sigma)}
   \ye^{\ereccvaleven(\sigma)}
   \ue^{\nrcpeakeven(\sigma)}
   \ve^{\nrcvaleven(\sigma)}
          \:\times
       \qquad\qquad
       \nonumber \\[0mm]
   & & \qquad\qquad\;\,
   \xo^{\eareccpeakodd(\sigma)}
   \yo^{\ereccvalodd(\sigma)}
   \uo^{\nrcpeakodd(\sigma)}
   \vo^{\nrcvalodd(\sigma)}
   \lambda^{\cyc(\sigma)}
   \;,
 \label{def.Qnhat.cycle-alternating.evenodd}
\end{eqnarray}
which generalize the polynomials \reff{def.Qn.cycle-alternating.evenodd}
by including a factor $\lambda^{\cyc(\sigma)}$.

Unfortunately, it is not possible to obtain an S-fraction
(or even a J-fraction) with polynomial coefficients
for \reff{def.Qnhat.cycle-alternating.evenodd}.
Indeed, it is not possible to obtain a J-fraction
with polynomial coefficients even for the 3-variable specialization
\be
   P_n(x,y,\lambda)
   \;=\; 
   \sum_{\sigma \in \Sym^{\rm ca}_{2n}}
       x^{\eareccpeak(\sigma)}
       y^{\ereccval(\sigma)}
        \lambda^{\cyc(\sigma)}
   \;.
 \label{eq.Pn.xylam}
\ee  
We show this in the Appendix.

However, it is possible to obtain a nice S-fraction
for the polynomials \reff{def.Qnhat.cycle-alternating.evenodd}
if we make the specializations $\ye=\ve$ and $\yo=\vo$:

% In \cite[Section~2.15]{Sokal-Zeng_masterpoly}
% another family 4-variable polynomials were introduced:
% \be
%    \widehat{Q}_n(x,y,u,\lambda)
%    \;=\;
%    \sum_{\sigma \in \Sym^{\rm ca}_{2n}}
%    x^{\eareccpeak(\sigma)}
%    u^{\nrcpeak(\sigma)}
%    y^{\cval(\sigma)}
%    \lambda^{\cyc(\sigma)}
%    \;.
%  \label{def.Qnhat.cycle-alternating}
% \ee
% It was then shown \cite[Theorem~2.21]{Sokal-Zeng_masterpoly}
% that the ordinary generating function of the polynomials $\widehat{Q}_n$
% has the S-type continued fraction
% \be
%    \sum_{n=0}^\infty \widehat{Q}_n(x,y,u,\lambda) \: t^n
%    \;=\;
%    \cfrac{1}{1 - \cfrac{\lambda x y t}{1 -  \cfrac{(\lambda\!+\! 1)(x\!+\! u)y t}{1 - \cfrac{(\lambda\!+\! 2)(x\!+\!2u)y t}{1 - \cdots}}}}
%    \label{eq.thm.perm.Stype.cycle-alternating.second}
% \ee
% with coefficients
% \be
%    \alpha_n   \;=\;   (\lambda+n-1) \: [x + (n-1)u]\: y
%    \;.
%  \label{def.weights.perm.Stype.cycle-alternating.second}
% \ee
% 
% We now generalize this by treating even and odd indices separately.
% We then have:

\begin{theorem}[Second S-fraction for cycle-alternating permutations]
   \label{thm.second.Sfrac}
The ordinary generating function of the polynomials
$\widehat{Q}_n(\xe,\ye,\ue,\ve,\xo,\yo,\uo,\vo,\lambda)$
with the specializations $\ye=\ve$ and $\yo=\vo$
has the S-type continued fraction
\be
   \sum_{n=0}^\infty \widehat{Q}_n(\xe,\ye,\ue,\ye,\xo,\yo,\uo,\yo,\lambda) \: t^n
   \;=\;
   \cfrac{1}{1 - \cfrac{\lambda\xe \yo t}{1 -  \cfrac{(\lambda\!+\! 1)(\xo\!+\!\uo)\ye t}{1 - \cfrac{(\lambda\!+\! 2)(\xe\!+\!2\ue)\yo t}{1 - \cdots}}}}
   \label{eq.thm.perm.Stype.cycle-alternating.evenodd.second}
\ee
with coefficients
\begin{subeqnarray}
   \alpha_{2k-1}  & = &  (\lambda+2k-2) \: [\xe + (2k-2) \ue] \: \yo \\[1mm]
   \alpha_{2k}    & = &  (\lambda+2k-1) \:[\xo + (2k-1) \uo] \: \ye
 \label{def.weights.perm.Stype.cycle-alternating.evenodd.second}
\end{subeqnarray}
\end{theorem}

\noindent
This generalizes \cite[Theorem~2.21]{Sokal-Zeng_masterpoly}
by treating even and odd indices separately.

The proof of Theorem~\ref{thm.second.Sfrac} will be obtained 
as a specialization of Theorem~\ref{thm.second.Sfrac.pqgen}.

\bigskip

When $\xe=\ye=\ue=\ve=\xo=\yo=\uo=\vo=1$, the polynomials $\widehat{Q}_n$
reduce to the secant power polynomials $E_{2n}(\lambda)$,
as will be shown in Section~\ref{subsec.contfrac.secpower}.

\subsection[$p,q$-generalization of second S-fraction]{$\bm{p,q}$-generalization of second S-fraction}

% In \cite[Section~2.15]{Sokal-Zeng_masterpoly},
% a 7-variable $p,q$-generalization of
% \reff{eq.thm.perm.Stype.cycle-alternating.second} was introduced:
% \begin{eqnarray}
%    & &
% 	\widehat{Q}_n(x,y,u,p_{+}, p_{-},q_{-} ,\lambda)
% 	\;=\;
% 	\sum_{\sigma \in \Sym^{\rm ca}_{2n}}
%    x^{\eareccpeak(\sigma)}   
%    u^{\nrcpeak(\sigma)}   y^{\cval(\sigma)}\times \qquad\qquad
%         \nonumber\\
%    & & \qquad\qquad\qquad
% 	 p_{-}^{\lcrosscpeak(\sigma)}
% 	 q_{-}^{\lnestcpeak(\sigma)}
%     p_{+}^{\ucrosscval(\sigma)+\unestcval(\sigma)}  \lambda^{\cyc(\sigma)}.
%  \label{def.Qnhat.cycle-alternating.pqgen}
% \end{eqnarray}
% %\be
% %   \eareccpeak(\sigma)
% %   \;=\;
% %   |\Eareccpeak(\sigma)|
% %   \;=\;
% %   |\Arec(\sigma) \cap \Cpeak(\sigma)|
% %\ee
% It was then shown \cite[Theorem~2.22]{Sokal-Zeng_masterpoly}
% that the ordinary generating function of the polynomials
% \reff{def.Qnhat.cycle-alternating.pqgen}
% has the S-type continued fraction
% 
% \be
%    \sum_{n=0}^\infty  \widehat{Q}_n(x,y,u,p_{+}, p_{-},q_{-} ,\lambda) \: t^n
%    \;=\;
%    \cfrac{1}{1 - \cfrac{\lambda x y t}{1 -  \cfrac{(\lambda\!+\! 1)(p_{-}x\!+\!q_{-}u)p_{+}y t}{1 - \cfrac{(\lambda\!+\! 2)(p_{-}^2 x\!+\!q_{-}[2]_{p_{-},q_{-}}u)p_{+}^2 y t}{1 - \cdots}}}}
%   \qquad
%    \label{eq.thm.perm.Stype.cycle-alternating.pqgen.second}
% \ee
% with coefficients
% \be
%    \alpha_n   \;=\;  (\lambda + n-1) \: (p_{-}^{n-1}x + q_{-}[n-1]_{p_{-},q_{-}}u) \: p_{+}^{n-1}y 
%    \;.
%  \label{def.weights.perm.Stype.cycle-alternating.pqgen.second}
% \ee
% 
% We now generalize this by treating even and odd indices separately.

We now introduce a $p,q$-generalization of
\reff{def.Qnhat.cycle-alternating.evenodd}:
\begin{eqnarray}
   & & \hspace*{-6mm}
   \widehat{Q}_n(\xe,\ye,\ue,\ve,\xo,\yo,\uo,\vo,p_{-1},p_{-2},p_{+1},p_{+2},q_{-1},q_{-2},q_{+1}, q_{+2},\lambda)
	\;=\;
        \nonumber \\[4mm]
   & & \qquad
   \sum_{\sigma \in \Sym^{\rm ca}_{2n}}
   \xe^{\eareccpeakeven(\sigma)}
   \ye^{\ereccvaleven(\sigma)}
   \ue^{\nrcpeakeven(\sigma)}
   \ve^{\nrcvaleven(\sigma)}
          \:\times
       \qquad\qquad
       \nonumber \\[-2mm]
   & & \qquad\qquad\;\,
   \xo^{\eareccpeakodd(\sigma)}
   \yo^{\ereccvalodd(\sigma)}
   \uo^{\nrcpeakodd(\sigma)}
   \vo^{\nrcvalodd(\sigma)}
	\:\times
       \qquad\qquad
       \nonumber \\[2mm]
   & & \qquad\qquad\;\,
   p_{-1}^{{\rm lcrosscpeakeven}(\sigma)}
   p_{-2}^{{\rm lcrosscpeakodd}(\sigma)}
   p_{+1}^{{\rm ucrosscvalodd}(\sigma)}
   p_{+2}^{{\rm ucrosscvaleven}(\sigma)}
	 \:\times
       \qquad\qquad
       \nonumber \\[2mm]
   & & \qquad\qquad\;\,
   q_{-1}^{{\rm lnestcpeakeven}(\sigma)}
   q_{-2}^{{\rm lnestcpeakodd}(\sigma)}
   q_{+1}^{{\rm unestcvalodd}(\sigma)}
   q_{+2}^{{\rm unestcvaleven}(\sigma)}
   \lambda^{\cyc(\sigma)}
   \;.
\label{def.Qnhat.cycle-alternating.evenodd.pqgen}
\end{eqnarray}
We then have:

\begin{theorem}[Second S-fraction for cycle-alternating permutations, $p,q$-generalization]
   \label{thm.second.Sfrac.pqgen}
The ordinary generating function of the polynomials~\reff{def.Qnhat.cycle-alternating.evenodd.pqgen}
%{\hbox{$Q_n(\xe,\ye,\ue,\ve,\xo,\yo,\uo,\vo,p_{-1},p_{-2},p_{+1},p_{+2},q_{-1},q_{-2},q_{+1},q_{+2})$}} 
with the specializations
$\ye=\ve$, $\yo=\vo$, $p_{+1} = q_{+1}$ and $p_{+2}=q_{+2}$
has the S-type continued fraction
\begin{eqnarray}
   & & \hspace*{-1cm}
   \sum_{n=0}^\infty \widehat{Q}_n(\xe,\ye,\ue,\ye,\xo,\yo,\uo,\yo,p_{-1},p_{-2},p_{+1},p_{+2},q_{-1},q_{-2},p_{+1},p_{+2},\lambda) \: t^n  \qquad
      \nonumber \\
   & & \hspace*{1cm}
   =\;
	\cfrac{1}{1 - \cfrac{\lambda\xe \yo t}{1 -  \cfrac{(\lambda+1)(p_{-2}\xo\!+\!q_{-2}\uo)p_{+2}\ye t}{1 - \cfrac{(\lambda+2)(p_{-1}^2\xe\!+\!q_{-1}[2]_{p_{-1},q_{-1}}\ue)p_{+1}^2\yo t}{1 - \cdots}}}}
   \label{eq.thm.perm.Stype.cycle-alternating.evenodd.pqgen.second}
\end{eqnarray}
with coefficients
\begin{subeqnarray}
   \alpha_{2k-1}  & = &  (\lambda+2k-2)\:(p_{-1}^{2k-2}\xe + q_{-1}[2k-2]_{p_{-1},q_{-1}} \ue) \: p_{+1}^{2k-2}\yo \nonumber\\
	\mbox{}\\
   \alpha_{2k}    & = &  (\lambda+2k-1)\:(p_{-2}^{2k-1}\xo + q_{-2}[2k-1]_{p_{-2},q_{-2}} \uo) \: p_{+2}^{2k-1}\ye \nonumber\\
	\mbox{}
 \label{def.weights.perm.Stype.cycle-alternating.evenodd.pqgen.second}
\end{subeqnarray}
\end{theorem}

\noindent
This generalizes \cite[Theorem~2.22]{Sokal-Zeng_masterpoly}
by treating even and odd indices separately.

The proof of this theorem
(and thus also of Theorem~\ref{thm.second.Sfrac})
will be based on the ``second master''
S-fraction for cycle-alternating permutations
\cite[Theorem~2.23]{Sokal-Zeng_masterpoly}
together with Lemma~\ref{lemma.cycle-alt}.
We first recall that second master S-fraction:
it states that the polynomials
\be
   \widehat{Q}_n(\bsfa,\bsfb,\lambda)
   \;=\;
   \sum_{\sigma \in \Sym^{\rm ca}_{2n}}
   \;\:\lambda^{\cyc(\sigma)}
   \prod\limits_{i \in \Cval}  \! \sfa_{\ucross(i,\sigma)\,+\,\unest(i,\sigma)}
   \prod\limits_{i \in \Cpeak} \!\! \sfb_{\lcross(i,\sigma),\,\lnest(i,\sigma)}
 \label{def.Qnhat.secondmaster.cycle-alternating}
\ee
in indeterminates
$\bsfa = (\sfa_{\ell})_{\ell \ge 0}$
and $\bsfb = (\sfb_{\ell,\ell'})_{\ell,\ell' \ge 0}$
have the S-type continued fraction
\be
   \sum_{n=0}^\infty \widehat{Q}_n(\bsfa,\bsfb,\lambda)  \: t^n
   \;=\;
   \cfrac{1}{1 - \cfrac{\lambda\sfa_{0} \sfb_{00} t}{1 - \cfrac{(\lambda + 1)\sfa_{1} (\sfb_{01} + \sfb_{10}) t}{1 - \cfrac{(\lambda + 2)\sfa_{2}(\sfb_{02} + \sfb_{11} + \sfb_{20}) t}{1 - \cdots}}}}
   \label{eq.thm.permutations.Stype.final1.cycle-alternating.second}
\ee
with coefficients
\be
   \alpha_n   \;=\;
   (\lambda + n-1) \sfa_{n-1}
   \biggl( \sum_{\ell=0}^{n-1} \sfb_{\ell,n-1-\ell} \biggr)
   \;.
   \label{eq.thm.permutations.Stype.final1.cycle-alternating.weights.second}
\ee
Notice that, in contrast to \reff{def.Qn.firstmaster.cycle-alternating},
the indeterminates $\bsfa$ now have only one index:
we are unable to count $\ucross$ and $\unest$ separately,
but can only count their sum.
This is the price we pay for being able to count cycles.

\proofof{Theorem~\ref{thm.second.Sfrac.pqgen}}
Notice that we obtain the polynomials
\reff{def.Qnhat.cycle-alternating.evenodd.pqgen}
%% $\widehat{Q}_n(\xe,\ye,\ue,\ye,\xo,\yo,\uo,\yo,p_{-1},p_{-2},p_{+1},p_{+2},q_{-1},q_{-2},p_{+1},p_{+2},\lambda)$\\
by specializing the polynomials 
$\widehat{Q}_n(\bsfa,\bsfb,\lambda)$ defined in
\reff{def.Qnhat.secondmaster.cycle-alternating}
to 
\begin{subeqnarray}
   \sfa_{\ell}
   & = &
   \begin{cases}
       p_{+1}^{\ell} \yo    &  \textrm{if $\ell$ is even}  \\
       
       p_{+2}^{\ell} \ye    &  \textrm{if $\ell$ is odd}  
   \end{cases}
       \\[2mm]
   \sfb_{\ell,\ell'}
   & = &
   \begin{cases}
       p_{-1}^{\ell} \xe    &  \textrm{if $\ell' = 0$ and $\ell+\ell'$ is even}  \\
       p_{-1}^{\ell} q_{-1}^{\ell'}\ue    &  \textrm{if $\ell' \ge 1$ and $\ell+\ell'$ is even} \\
       p_{-2}^{\ell} \xo    &  \textrm{if $\ell' = 0$ and $\ell+\ell'$ is odd}  \\
       p_{-2}^{\ell} q_{-2}^{\ell'}\uo    &  \textrm{if $\ell' \ge 1$ and $\ell+\ell'$ is odd}
   \end{cases}
 \label{eq.sfasfb.pq}
\end{subeqnarray}
Inserting these into
\reff{eq.thm.permutations.Stype.final1.cycle-alternating.weights.second}
yields
\reff{def.weights.perm.Stype.cycle-alternating.evenodd.pqgen.second}.
This along with Lemma~\ref{lemma.cycle-alt} completes the proof of Theorem~\ref{thm.second.Sfrac.pqgen}.
\qed

\proofof{Theorem~\ref{thm.second.Sfrac}}
Specialize Theorem~\ref{thm.second.Sfrac.pqgen}
to $p_{+1}=p_{+2}=p_{-1}=p_{-2}=q_{-1}=q_{-2} = 1$.
\qed

\subsection{Alternating cycles}

Let $\Cyc_N \subseteq \Sym_N$ denote the set of permutations of $[N]$
consisting of exactly one cycle.
Then let $\Altcyc_{2n} = \Cyc_{2n} \cap \Sym_{2n}^{\rm ca}$
denote the set of cycle-alternating permutations of $[2n]$
consisting of exactly one cycle,
or \textbfit{alternating cycles} for short.
Dumont \cite[section~3 and p.~40]{Dumont_86} showed that
alternating cycles are enumerated by the
down-shifted tangent numbers, i.e.~$|\Altcyc_{2n}| = E_{2n-1}$.
Here we will recover this result,
and extend it to a more refined S-fraction.

We can enumerate alternating cycles by extracting the coefficient
of $\lambda^1$ in the preceding results
(namely, Theorems~\ref{thm.second.Sfrac} and \ref{thm.second.Sfrac.pqgen}
 and eq.~\reff{eq.thm.permutations.Stype.final1.cycle-alternating.second}).
We begin by defining the polynomials analogous to
\reff{def.Qn.cycle-alternating.evenodd},
\reff{def.Qn.cycle-alternating.evenodd.pqgen}
and \reff{def.Qnhat.secondmaster.cycle-alternating}
but restricted to alternating cycles:
\begin{eqnarray}
   & &
   Q_n^\Altcyc(\xe,\ye,\ue,\ve,\xo,\yo,\uo,\vo)
   \;=\;
        \nonumber \\[4mm]
   & & \qquad
   \sum_{\sigma \in \Altcyc_{2n}}
   \xe^{\eareccpeakeven(\sigma)}
   \ye^{\ereccvaleven(\sigma)}
   \ue^{\nrcpeakeven(\sigma)}
   \ve^{\nrcvaleven(\sigma)}
          \:\times
       \qquad\qquad
       \nonumber \\[0mm]
   & & \qquad\qquad\;\,
   \xo^{\eareccpeakodd(\sigma)}
   \yo^{\ereccvalodd(\sigma)}
   \uo^{\nrcpeakodd(\sigma)}
   \vo^{\nrcvalodd(\sigma)}
 \label{def.Qn.altcyc}
       \\[4mm]
   & &
        Q_n^\Altcyc(\xe,\ye,\ue,\ve,\xo,\yo,\uo,\vo,p_{-1},p_{-2},p_{+1},p_{+2},q_{-1},q_{-2},q_{+1},q_{+2})
   \;=\;
        \nonumber \\[4mm]
   & & \qquad
   \sum_{\sigma \in \Altcyc_{2n}}
   \xe^{\eareccpeakeven(\sigma)}
   \ye^{\ereccvaleven(\sigma)}
   \ue^{\nrcpeakeven(\sigma)}
   \ve^{\nrcvaleven(\sigma)}
          \:\times
       \qquad\qquad
       \nonumber \\[-2mm]
   & & \qquad\qquad\;\,
   \xo^{\eareccpeakodd(\sigma)}
   \yo^{\ereccvalodd(\sigma)}
   \uo^{\nrcpeakodd(\sigma)}
   \vo^{\nrcvalodd(\sigma)}
        \:\times
       \qquad\qquad
       \nonumber \\[2mm]
   & & \qquad\qquad\;\,
   p_{-1}^{{\rm lcrosscpeakeven}(\sigma)}
   p_{-2}^{{\rm lcrosscpeakodd}(\sigma)}
   p_{+1}^{{\rm ucrosscvalodd}(\sigma)}
   p_{+2}^{{\rm ucrosscvaleven}(\sigma)}
         \:\times
       \qquad\qquad
       \nonumber \\[2mm]
   & & \qquad\qquad\;\,
   q_{-1}^{{\rm lnestcpeakeven}(\sigma)}
   q_{-2}^{{\rm lnestcpeakodd}(\sigma)}
   q_{+1}^{{\rm unestcvalodd}(\sigma)}
   q_{+2}^{{\rm unestcvaleven}(\sigma)}
 \label{def.Qn.altcyc.pqgen}
       \\[4mm]
   & &
   Q^\Altcyc_n(\bsfa,\bsfb)
   \;=\;
   \sum_{\sigma \in \Altcyc_{2n}}
   \;\:
   \prod\limits_{i \in \Cval}  \! \sfa_{\ucross(i,\sigma)\,+\,\unest(i,\sigma)}
   \prod\limits_{i \in \Cpeak} \!\! \sfb_{\lcross(i,\sigma),\,\lnest(i,\sigma)}
 \label{def.Qn.altcyc.master}
\end{eqnarray}

Since the results for alternating cycles will arise by specialization
of the results that include counting of cycles,
the same specializations are needed to obtain nice continued fractions.
Extracting the coefficient of $\lambda^1$ in Theorem~\ref{thm.second.Sfrac},
we deduce the following:

\begin{corollary}[S-fraction for alternating cycles]
   \label{cor.second.Sfrac.altcyc}
The ordinary generating function of the polynomials \reff{def.Qn.altcyc}
%% $Q^\Altcyc_{n+1}(\xe,\ye,\ue,\ve,\xo,\yo,\uo,\vo)$
with the specializations $\ye=\ve$ and $\yo=\vo$
has the S-type continued fraction
\be
   \sum_{n=0}^\infty Q^\Altcyc_{n+1}(\xe,\ye,\ue,\ye,\xo,\yo,\uo,\yo) \: t^n
   \;=\;
   \cfrac{\xe \yo}{1 -  \cfrac{(\xo\!+\!\uo)\ye t}{1 - \cfrac{ 2(\xe\!+\!2\ue)\yo t}{1 - \cdots}}}
   \label{eq.cor.second.Sfrac.altcyc}
\ee
with coefficients
\begin{subeqnarray}
   \alpha_{2k-1}  & = &  (2k-1) \: [\xo + (2k-1) \uo] \: \ye \\[1mm]
   \alpha_{2k}    & = &  2k \:[\xe + 2k \ue] \: \yo
 \label{def.weights.cor.second.Sfrac.altcyc}
\end{subeqnarray}
\end{corollary}

Specializing to $\xe=\ye=\ue=\xo=\yo=\uo=1$,
we see that $|\Altcyc_{2n+2}|$ has an S-fraction
with coefficients $\alpha_n = n(n+1)$,
which is well known to be the S-fraction
for the tangent numbers $E_{2n+1}$ \cite[A000182]{OEIS}.
In other words, $|\Altcyc_{2n}| = E_{2n-1}$.

Similarly, extracting the coefficient of $\lambda^1$ in
Theorem~\ref{thm.second.Sfrac.pqgen}, we deduce:

\begin{corollary}[S-fraction for alternating cycles, $p,q$-generalization]
   \label{cor.second.Sfrac.pqgen}
The ordinary generating function of the
polynomials~\reff{def.Qn.altcyc.pqgen}
with the specializations
$\ye=\ve$, $\yo=\vo$, $p_{+1} = q_{+1}$ and $p_{+2}=q_{+2}$
has the S-type continued fraction
\begin{eqnarray}
   & & \hspace*{-1cm}
   \sum_{n=0}^\infty Q^\Altcyc_{n+1}(\xe,\ye,\ue,\ye,\xo,\yo,\uo,\yo,p_{-1},p_{-2},p_{+1},p_{+2},q_{-1},q_{-2},p_{+1},p_{+2}) \: t^n  \qquad
     \nonumber\\
   & & \hspace*{1cm}
   =\;
	\cfrac{\xe \yo}{1 -  \cfrac{(p_{-2}\xo\!+\!q_{-2}\uo)p_{+2}\ye t}{1 - \cfrac{2(p_{-1}^2\xe\!+\!q_{-1}[2]_{p_{-1},q_{-1}}\ue)p_{+1}^2\yo t}{1 - \cdots}}}
   \label{eq.cor.second.Sfrac.pqgen}
\end{eqnarray}
with coefficients
\begin{subeqnarray}
   \alpha_{2k-1}    & = &  (2k-1)\:(p_{-2}^{2k-1}\xo + q_{-2}[2k-1]_{p_{-2},q_{-2}} \uo) \: p_{+2}^{2k-1}\ye \\
	\mbox{}
   \alpha_{2k}  & = &  2k \:(p_{-1}^{2k}\xe + q_{-1}[2k]_{p_{-1},q_{-1}} \ue) \: p_{+1}^{2k}\yo 
 \label{def.weights.cor.second.Sfrac.pqgen}
\end{subeqnarray}
\end{corollary}

Finally, extracting the coefficient of $\lambda^1$ in
\reff{eq.thm.permutations.Stype.final1.cycle-alternating.second}/%
\reff{eq.thm.permutations.Stype.final1.cycle-alternating.weights.second},
we obtain:

\begin{corollary}[Master S-fraction for alternating cycles]
   \label{cor.second.Sfrac.master}
The ordinary generating function of the
polynomials~\reff{def.Qn.altcyc.master}
has the S-type continued fraction
\be
   \sum_{n=0}^\infty Q_{n+1}^\Altcyc(\bsfa,\bsfb)  \: t^n
   \;=\;
   \cfrac{\sfa_{0} \sfb_{00}}{1 - \cfrac{\sfa_{1} (\sfb_{01} + \sfb_{10}) t}{1 - \cfrac{2\sfa_{2}(\sfb_{02} + \sfb_{11} + \sfb_{20}) t}{1 - \cdots}}}
   \label{eq.cor.second.Sfrac.master}
\ee
with coefficients
\be
   \alpha_n   \;=\;
   n \, \sfa_{n}
   \biggl( \sum_{\ell=0}^{n} \sfb_{\ell,n-\ell} \biggr)
   \;.
   \label{def.weights.cor.second.Sfrac.master}
\ee
\end{corollary}

\subsection{Open problem: Combinatorial model for the tangent numbers}

Let us start by noting the inequalities
\be
   E_{2n-1} \;\le\; E_{2n} \;\le\; E_{2n+1} \;\le\; (2n)!
   \;.
\ee
We have just seen that
the down-shifted tangent numbers $E_{2n-1}$ enumerate the set $\Altcyc_{2n}$
of alternating cycles on $[2n]$,
while the secant numbers $E_{2n}$ enumerate the set $\Sym_{2n}^{\rm ca}$
of cycle-alternating permutations of $[2n]$.
The inequality $E_{2n} \le E_{2n+1} \le (2n)!$
suggests that also the tangent numbers $E_{2n+1}$
should enumerate some class of permutations of $[2n]$
that contains $\Sym_{2n}^{\rm ca}$.
What class?  We pose this as an open problem:

\begin{openproblem}
\hfill\quad\break
\vspace*{-7mm}
\begin{itemize}
   \item[(a)]
Find a natural class $\scrc$ satisfying
$\Sym_{2n}^{\rm ca} \subseteq \scrc \subseteq \Sym_{2n}$
and $|\scrc| = E_{2n+1}$.
   \item[(b)]
Find statistics on $\scrc$ whose generating polynomials
have a nice S-fraction that generalizes the coefficients
$\alpha_n = n(n+1)$ for the tangent numbers.
Ideally these statistics would also be such that,
by specializing some variables to zero,
we reobtain some or all of the continued fractions
of the present paper for cycle-alternating permutations.
\end{itemize}
\end{openproblem}

We remark that an analogous situation occurs
for the Genocchi numbers $g_n$ and the median Genocchi numbers $h_n$
(in the notation of \cite{Deb-Sokal_genocchi}),
which satisfy
\be
   g_{n-1} \;\le\; h_n \;\le\; g_n \;\le\; h_{n+1}  \;\le\; (2n)!
   \;.
\ee
In this case the relevant objects are a class of permutations of $[2n]$
called {\em D-permutations}\/ \cite{Lazar_22,Lazar_20}
and their subsets
\be
   \hbox{D-cycles}  \;\subseteq\;
   \hbox{D-derangements}  \;\subseteq\;
   \hbox{D-semiderangements}  \;\subseteq\;
   \hbox{D-permutations}
   \;,
\ee
and we have \cite{Dumont_74,Dumont_94,Lazar_22,Lazar_20,Deb-Sokal_genocchi}
\begin{subeqnarray}
   |\hbox{D-cycles}|        & = &  g_{n-1}   \\[1mm]
   |\hbox{D-derangements}|  & = &  h_n       \\[1mm]
   |\hbox{D-semiderangements}|  & = &  g_n       \\[1mm]
   |\hbox{D-permutations}|  & = &  h_{n+1}
\end{subeqnarray}
Furthermore, in \cite{Deb-Sokal_genocchi} we obtained continued fractions
for multivariate polynomials enumerating D-permutations,
which by specializing some variables to zero
gave continued fractions enumerating D-semiderangements and D-derangements;
furthermore, by extracting the coefficient of $\lambda^1$
in a cycle-counting result,
we obtained continued fractions enumerating D-cycles.
It would be nice to have analogous results in the present context.

\section{Alternating Laguerre digraphs}   \label{sec.alternating}

We now place permutations in a more general context:
that of Laguerre digraphs.
Recall from the Introduction that a \textbfit{Laguerre digraph}
\cite{Foata_84,Sokal_multiple_laguerre,latpath_laguerre}
is a digraph in which each vertex has out-degree 0 or~1 and in-degree 0 or~1.
%It follows that each weakly connected component of a Laguerre digraph
% \footnote{
%    By a ``weakly connected component'' of a digraph,
%    we mean a connected component of the underlying undirected graph.
%% THIS SHOULD BE "CORRESPONDING TO" a connected component of the
%%    underlying undirected graph.
% }
%is either a directed path of some length $\ell \ge 0$
%(where a path of length 0 is an isolated vertex)
%or else a directed cycle of some length $\ell \ge 1$.
%(where a cycle of length 1 is a loop).
For each integer $n \ge 0$,
let us write $\LD_n$ for the set of Laguerre digraphs
on the vertex set $[n] \eqdef \{1,\ldots,n\}$;
and for $n \ge k \ge 0$,
let us write $\LD_{n,k}$ for the set of Laguerre digraphs
on the vertex set $[n]$ with $k$ paths.
A Laguerre digraph with no paths is simply a collection of directed cycles,
i.e.~the digraph associated to a permutation $\sigma \in \Sym_n$
\cite[pp.~22--23]{Stanley_12}.
There is thus a natural bijection $\LD_{n,0} \simeq \Sym_n$.

The study of Laguerre digraphs $\LD_{n,k}$
thus generalizes the study of permutations $\LD_{n,0} \simeq \Sym_n$.
In a similar way, we will define a class of
{\em alternating Laguerre digraphs}\/ $\LD^{\rm alt}_{n,k}$
so as to generalize the study of cycle-alternating permutations
$\LD^{\rm alt}_{n,0} \simeq \Sym^{\rm ca}_n$.
Just as $\Sym^{\rm ca}_n \neq \emptyset$ only for $n$ even,
we will see that $\LD^{\rm alt}_{n,k} \neq \emptyset$ only for $n+k$ even.

We begin (Section~\ref{subsec.egf})
by deriving an exponential generating function
for Laguerre digraphs with a very general set of weights
\cite{latpath_laguerre}.
Then (Section~\ref{subsec.altlag})
we define alternating Laguerre digraphs
and specialize to obtain their exponential generating function.
Next (Section~\ref{subsec.contfrac.secpower})
we obtain the generalized Jacobi--Rogers and Stieltjes--Rogers polynomials
associated to the secant power polynomials,
using exponential generating functions.
Then (Section~\ref{subsec.contfrac.secondmaster})
we obtain the generalized Stieltjes--Rogers polynomials
associated to the second master S-fraction
for cycle-alternating permutations
[eqns.~\reff{eq.thm.permutations.Stype.final1.cycle-alternating.second}/%
\reff{eq.thm.permutations.Stype.final1.cycle-alternating.weights.second}],
using combinatorial methods.
Finally (Section~\ref{subsec.contfrac.specialization})
we specialize this latter result
to obtain the generalized Stieltjes-Rogers polynomials
corresponding to the S-fractions of
Theorems~\ref{thm.second.Sfrac} and \ref{thm.second.Sfrac.pqgen}.

\subsection{Exponential generating function for Laguerre digraphs}
   \label{subsec.egf}

Before defining our statistics on Laguerre digraphs,
we first need to make a convention about {\em boundary conditions}\/
at the two ends of a path.
We will here use 0--0 boundary conditions:
that is, we extend the Laguerre digraph $G$ on the vertex set $[n]$
to a digraph $\widehat{G}$ on the vertex set $[n] \cup \{0\}$
by decreeing that any vertex $i \in [n]$
that has in-degree (resp.~out-degree) 0 in $G$
will receive an incoming (resp.~outgoing) edge
from (resp.~to) the vertex 0.
In this way each vertex $i \in [n]$
will have a unique predecessor $p(i) \in [n] \cup \{0\}$
and a unique successor $s(i) \in [n] \cup \{0\}$.
We then say that a vertex $i \in [n]$ is a
\begin{itemize}
   \item {\em peak}\/ (p) if $p(i) < i > s(i)$;
   \item {\em valley}\/ (v) if $p(i) > i < s(i)$;
   \item {\em double ascent}\/ (da) if $p(i) < i < s(i)$;
   \item {\em double descent}\/ (dd) if $p(i) > i > s(i)$;
   \item {\em fixed point}\/ (fp) if $p(i) = i = s(i)$.
\end{itemize}
(Note that ``fixed point'' is a synonym of ``loop''.)
When these concepts are applied to the cycles of a Laguerre digraph,
we obtain the usual \textbfit{cycle classification}
of indices in a permutation
as cycle peaks, cycle valleys, cycle double rises, cycle double falls
and fixed points \cite{Zeng_93,Sokal-Zeng_masterpoly}.
When applied to the paths of a Laguerre digraph,
we obtain the usual \textbfit{linear classification}
of indices in a permutation (written in word form)
as peaks, valleys, double ascents or double descents
\cite[p.~45]{Stanley_12}.
Note that, because of the 0--0 boundary conditions,
an isolated vertex is always a peak,
the initial vertex of a path is always a peak or double ascent,
and the final vertex of a path is always a peak or double descent;
moreover, each path contains at least one peak.

We write $\pk(G), \val(G), \da(G), \dd(G), \fp(G)$
for the number of vertices $i \in [n]$ that are, respectively,
peaks, valleys, double ascents, double descents or fixed points.
We then introduce the
\textbfit{multivariate Laguerre coefficient matrix}
\be
   \sfLhat^{(\alpha)}(\yp,\yv,\yda,\ydd,\yfp)_{n,k}
   \;\eqdef\;
   \sum_{G \in \LD_{n,k}}
       y_{\rm p}^{\pk(G)} y_{\rm v}^{\val(G)}
       y_{\rm da}^{\da(G)} y_{\rm dd}^{\dd(G)} y_{\rm fp}^{\fp(G)}
%%            \, w_0^{\iv(G)} \, w_{\ge 1}^{\pa_{\ge 1}(G)}
            \, (1+\alpha)^{\cyc(G)}
   \:.
 \label{def.coeffmat.gen.second}
\ee
This polynomial is homogeneous of degree $n$ in $\yp,\yv,\yda,\ydd,\yfp$.
(The weight per cycle is here written as $1+\alpha$ instead of $\lambda$
 because these polynomials generalize the Laguerre polynomials
 with parameter $\alpha$.)

Indeed, we can go farther, by assigning different weights
for the vertices belonging to a cycle or to a path.
For a Laguerre digraph $G$,
let us write $\pkcyc(G)$, $\valcyc(G)$, $\dacyc(G)$, $\ddcyc(G)$, $\fp(G)$
for the number of peaks, valleys, double ascents, double descents
and fixed points that belong to a cycle of $G$,
and $\pkpa(G)$, $\valpa(G)$, $\dapa(G)$, $\ddpa(G)$
for the number of peaks, valleys, double ascents and double descents
that belong to a path of $G$
(of course fixed points can only belong to a cycle).
We then assign weights $\yp,\yv,\yda,\ydd,\yfp$
to the vertices belonging to a cycle,
and weights $\zp,\zv,\zda,\zdd$
to the vertices belonging to a path.
We therefore define the
\textbfit{generalized multivariate Laguerre coefficient matrix}
\begin{eqnarray}
   & &
   \!\!\!\!
   \sfLtilde^{(\alpha)}(\yp,\yv,\yda,\ydd,\yfp, \zp,\zv,\zda,\zdd)_{n,k}
         \nonumber \\[3mm]
   & &
   \quad\eqdef
   \sum_{G \in \LD_{n,k}}
       y_{\rm p}^{\pkcyc(G)} y_{\rm v}^{\valcyc(G)}
       y_{\rm da}^{\dacyc(G)} y_{\rm dd}^{\ddcyc(G)} y_{\rm fp}^{\fp(G)}
       \,
       z_{\rm p}^{\pkpa(G)} z_{\rm v}^{\valpa(G)}
       z_{\rm da}^{\dapa(G)} z_{\rm dd}^{\ddpa(G)}
         \:\times
         \nonumber \\[-1mm]
   & & \qquad\qquad\quad\;\,
%%            \, w_0^{\iv(G)} \, w_{\ge 1}^{\pa_{\ge 1}(G)}
            \, (1+\alpha)^{\cyc(G)}
   \:.
   \qquad
 \label{def.coeffmat.gen.second.gen}
\end{eqnarray}
This polynomial is homogeneous of degree $n$ in
$\yp,\yv,\yda,\ydd,\yfp,\zp,\zv,\zda,\zdd$.

We now proceed to compute the exponential generating functions
for the matrices \reff{def.coeffmat.gen.second.gen}.
We do this by combining the known exponential generating functions
for permutations with cyclic statistics \cite[Th\'eor\`eme~1]{Zeng_93}
and permutations with linear statistics \cite[Proposition~4]{Zeng_93}.
These formulae are as follows:

\medskip

1) Permutations with cyclic statistics are enumerated by the polynomials
\be
   P_n^{\rm cyc}(\yp,\yv,\yda,\ydd,\yfp,\lambda)
   \;\eqdef\;
   \sum_{\sigma \in \Sym_n}
      y_{\rm p}^{\pkcyc(\sigma)} y_{\rm v}^{\valcyc(\sigma)}
       y_{\rm da}^{\dacyc(\sigma)} y_{\rm dd}^{\ddcyc(\sigma)}
       y_{\rm fp}^{\fp(\sigma)}
      \, \lambda^{\cyc(\sigma)}
\ee
where $\pkcyc(\sigma)$, $\valcyc(\sigma)$, $\dacyc(\sigma)$, $\ddcyc(\sigma)$,
$\fp(\sigma)$
denote the number of cycle peaks, cycle valleys, cycle double rises,
cycle double falls and fixed points in $\sigma$,
and $\cyc(\sigma)$ denotes the number of cycles in $\sigma$.
By convention we set $P_0^{\rm cyc} = 1$.
The polynomial $P_n^{\rm cyc}$ is homogeneous of degree $n$ in
$\yp,\yv,\yda,\ydd,\yfp$.
We write
\be
   F(t;\yp,\yv,\yda,\ydd,\yfp,\lambda)
   \;\eqdef\;
   \sum_{n=0}^\infty P_n^{\rm cyc}(\yp,\yv,\yda,\ydd,\yfp,\lambda)
       \: {t^n \over n!}
\ee
for the corresponding exponential generating function.

\begin{lemma} {\hspace*{-2mm}\bf\cite[Th\'eor\`eme~1]{Zeng_93}\ }
   \label{lemma.zeng.cyclic}
We have
\begin{subeqnarray}
   F(t;\yp,\yv,\yda,\ydd,\yfp,\lambda)
   & = &
   e^{\lambda \yfp t}
   \,
   \biggl( {r_1 \,-\, r_2
            \over
            r_1 e^{r_2 t} \,-\, r_2 e^{r_1 t}
           }
   \biggr)^{\! \lambda}
       \\[2mm]
   & = &
   F(t;\yp,\yv,\yda,\ydd,\yfp,1)^\lambda
 \label{eq.lemma.zeng.cyclic.1}
\end{subeqnarray}
where $r_1 r_2 = \yp \yv$ and $r_1 + r_2 = \yda + \ydd$.
Otherwise put, $r_1$ and $r_2$ are the roots (in~either order)
of the quadratic equation $\rho^2 - (\yda + \ydd) \rho + \yp \yv = 0$.
Concretely,
\be
   r_{1,2}
   \;=\;
   {\yda + \ydd \pm \sqrt{(\yda + \ydd)^2 - 4\yp \yv}
    \over
    2
   }
   \;.
 \label{eq.lemma.zeng.cyclic.2}
\ee
\end{lemma}

\medskip

2) Permutations with linear statistics and 0--0 boundary conditions
are enumerated by the polynomials
\be
   P_n^{\rm lin(00)}(\zp,\zv,\zda,\zdd)
   \;\eqdef\;
   \sum_{\sigma \in \Sym_n}
      z_{\rm p}^{\pkpa(\sigma)} z_{\rm v}^{\valpa(\sigma)}
       z_{\rm da}^{\dapa(\sigma)} z_{\rm dd}^{\ddpa(\sigma)}
\ee
where $\pkpa(\sigma)$, $\valpa(\sigma)$, $\dapa(\sigma)$, $\ddpa(\sigma)$
denote the number of peaks, valleys, double ascents and double descents
in the permutation $\sigma$ written as a word $\sigma_1 \cdots \sigma_n$,
where we impose the boundary conditions $\sigma_0 = \sigma_{n+1} = 0$.
By convention we restrict attention to $n \ge 1$.
The polynomial $P_n^{\rm lin(00)}$ is homogeneous of degree $n$ in
$\zp,\zv,\zda,\zdd$.
We write
\be
   G(t;\zp,\zv,\zda,\zdd)
   \;\eqdef\;
   \sum_{n=1}^\infty P_n^{\rm lin(00)}(\zp,\zv,\zda,\zdd)
       \: {t^n \over n!}
\ee
for the corresponding exponential generating function
(note that the sum starts at $n=1$).

\begin{lemma}  {\hspace*{-2mm}\bf\cite[Proposition~4]{Zeng_93}\ }
   \label{lemma.zeng.linear}
We have
\be
   G(t;\zp,\zv,\zda,\zdd)
   \;=\;
   \zp
   \:
   \biggl( {e^{r_1 t} \,-\, e^{r_2 t}
            \over
            r_1 e^{r_2 t} \,-\, r_2 e^{r_1 t}
           }
   \biggr)
 \label{eq.lemma.zeng.linear.1}
\ee
where $r_1 r_2 = \zp \zv$ and $r_1 + r_2 = \zda + \zdd$
analogously to Lemma~\ref{lemma.zeng.cyclic}.
This function satisfies the differential equation
\be
   G'(t)  \;=\; \zp \,+\, (\zda+\zdd) G(t) \,+\, \zv G(t)^2
   \;.
 \label{eq.lemma.zeng.linear.2}
\ee
\end{lemma}

We can now put these ingredients together to determine
the exponential generating functions
for the matrix \reff{def.coeffmat.gen.second.gen}.
A Laguerre digraph $G \in \LD_{n,k}$
consists of a permutation (that is, a collection of disjoint cycles)
on some subset $S \subseteq [n]$
together with $k$ disjoint paths on $[n] \setminus S$.
Each of these paths can be considered as a permutation written in word form.
%% In the generating function \reff{eq.lemma.zeng.linear.1},
%% isolated vertices get a weight $\zp t$
%% and the paths of length $\ge 1$ get a total weight $G(t) - \zp t$;
%% but we want to multiply the former by $w_0$ and the latter by $w_{\ge 1}$.
%% The resulting generating function is thus
%% \be
%%    \widehat{G}(t;\zp,\zv,\zda,\zdd,w_0,w_{\ge 1})
%%    \;\eqdef\;
%%    w_{\ge 1} \, G(t;\zp,\zv,\zda,\zdd)
%%      \:+\: (w_0 - w_{\ge 1}) \zp t
%%    \;.
%% \ee
%% {\bf But now I worry that $\widehat{G}$ does not satisfy any obvious
%%    autonomous differential equation, so I won't know what the A-series is
%%    when $w_0 \neq w_{\ge 1}$.
%%    But in my OLD section (next) it looks like I \emph{did} get an A-series
%%    for $w_0 \neq w_{\ge 1}$, albeit for a kind of row-generating matrix
%%    that has an additional variable $x$.  WHAT IS GOING ON HERE???
%%    Maybe the point is that those previous formulae did not insist on
%%    what the meaning of $k$ was, and implicitly made $k$ be the number
%%    of paths \emph{of length $\ge 1$} (i.e.~threw the isolated vertices
%%    into $F$), while here we have declared that $k$ is the total number
%%    of paths.}
By the exponential formula, the exponential generating function
for the $k$th column of the matrix \reff{def.coeffmat.gen.second.gen} is then
\begin{eqnarray}
   & &
   \sum_{n=0}^\infty
      \sfLtilde^{(\alpha)}(\yp,\yv,\yda,\ydd,\yfp, \zp,\zv,\zda,\zdd)_{n,k}
      \: {t^n \over n!}
             \nonumber \\[3mm]
   & &
   \qquad =\;
   F(t;\yp,\yv,\yda,\ydd,\yfp,1)^{1+\alpha}
   \;
   {G(t;\zp,\zv,\zda,\zdd)^k \over k!}
   \;,
   \qquad\qquad
\end{eqnarray}
where the $1/k!$ comes because the paths are indistinguishable.
The bivariate exponential generating function is therefore
\begin{eqnarray}
   & &
   \sum_{n=0}^\infty \sum_{k=0}^n
      \sfLtilde^{(\alpha)}(\yp,\yv,\yda,\ydd,\yfp, \zp,\zv,\zda,\zdd)_{n,k}
      \: {t^n \over n!} \, u^k
             \nonumber \\[3mm]
   & &
   \qquad =\;
   F(t;\yp,\yv,\yda,\ydd,\yfp,1)^{1+\alpha}
   \;
   \exp\bigl[ u \, G(t;\zp,\zv,\zda,\zdd) \bigr]
   \;.
   \qquad\qquad
\end{eqnarray}

\subsection{Alternating Laguerre digraphs}   \label{subsec.altlag}

We define an \textbfit{alternating Laguerre digraph}
to be a Laguerre digraph in which there are no
double ascents, double descents or fixed points.
It follows that all the cycles are of even length
and consist of cycle peaks and cycle valleys in alternation.
Furthermore, since we are using 0--0 boundary conditions,
all the paths are of odd length and consist of peaks and valleys
in alternation, starting and ending with a peak:
$0 < x_1 > x_2 < x_3 > \ldots < x_{2m+1} > 0$.
That is, the paths are simply alternating (down-up) permutations
\cite{Stanley_10} of odd length written in word form:
$x_1 > x_2 < x_3 > \ldots < x_{2m+1}$.
We write $\LD^{\rm alt}_{n,k}$ for the set of
alternating Laguerre digraphs on the vertex set $[n]$ with $k$ paths.
Since each path has odd length and each cycle has even length,
it follows that $\LD^{\rm alt}_{n,k} \neq \emptyset$ only if $n+k$ is even.
We have the obvious bijection $\LD^{\rm alt}_{n,0} \simeq \Sym^{\rm ca}_n$.
In this way, alternating Laguerre digraphs generalize
cycle-alternating permutations, to which they reduce when $k=0$.

We define a generating polynomial for alternating Laguerre digraphs
by specializing the one for Laguerre digraphs
to forbid double ascents, double descents and fixed points:
\begin{subeqnarray}
   \sfLtilde^{(\alpha)\rm alt}(\yp,\yv, \zp,\zv)_{n,k}
   & \eqdef &
      \sfLtilde^{(\alpha)}(\yp,\yv,0,0,0, \zp,\zv,0,0)_{n,k}
     \\[2mm]
   & = &
   \sum_{G \in \LD^{\rm alt}_{n,k}}
       y_{\rm p}^{\pkcyc(G)} y_{\rm v}^{\valcyc(G)}
       \,
       z_{\rm p}^{\pkpa(G)} z_{\rm v}^{\valpa(G)}
         \:\times
       \, (1+\alpha)^{\cyc(G)}
   \:.
      \nonumber \\[-5mm]
 \label{def.coeffmat.gen.second.gen.alt}
\end{subeqnarray}
This polynomial is homogeneous of degree $n$ in $\yp,\yv,\zp,\zv$.

Specializing Lemma~\ref{lemma.zeng.cyclic} to $\yda = \ydd = \yfp = 0$,
we have $r_1 = +i \sqrt{\yp \yv}$ and $r_2 = -i \sqrt{\yp \yv}$,
hence
\be
   F(t;\yp,\yv,0,0,0,\lambda)
   \;=\;
   \bigl[ \sec\big(\sqrt{\yp \yv} \, t\big) \bigr]^{\!\lambda}
   \;.
\ee
Similarly, specializing Lemma~\ref{lemma.zeng.linear} to $\zda = \zdd = 0$,
we have
\be
   G(t;\zp,\zv,0,0)
   \;=\;
   \sqrt{\zp/\zv} \: \tan\big(\sqrt{\zp \zv} \, t\big)
   \;.
\ee
Therefore, the exponential generating function
for the $k$th column of the matrix \reff{def.coeffmat.gen.second.gen.alt} is
\be
   \sum_{n=0}^\infty
      \sfLtilde^{(\alpha)\rm alt}(\yp,\yv, \zp,\zv)_{n,k}
      \: {t^n \over n!}
             \nonumber \\[3mm]
   \;=\;
   \bigl[\sec\big(\sqrt{\yp \yv} \, t\big)\bigr]^{\! 1+\alpha}
   \;
   {\bigl[ \sqrt{\zp/\zv} \: \tan\big(\sqrt{\zp \zv} \, t\big)\bigr]^k \over k!}
   \;.
\ee

In order to make the connection with continued fractions,
we will want to give path vertices and cycle vertices
the same weights, i.e.~$\zp = \yp$ and $\zv = \yv$.
Once we have made this specialization,
there is no loss of generality in taking $\yp = \yv = 1$,
since those variables would simply make a trivial rescaling.
We therefore have
\be
   \sum_{n=0}^\infty
      \sfLtilde^{(\alpha)\rm alt}(1,1,1,1)_{n,k}
      \: {t^n \over n!}
   \;=\;
   (\sec t)^{1+\alpha}
   \;
   {(\tan t)^k \over k!}
   \;.
 \label{eq.columnegf.altlaguerre}
\ee
In particular, the zeroth column $\sfLtilde^{(\alpha)\rm alt}(1,1,1,1)_{n,0}$
has the \textbfit{secant power polynomials} at even $n$
and zero entries at odd $n$,
where the secant power polynomials are defined by
\be
   (\sec t )^\lambda
   \;=\;
   \sum_{n=0}^\infty E_{2n}(\lambda) \: {t^{2n} \over (2n)!}
   \;.
\ee

The first 8 rows and columns of the matrix
$\sfLtilde^{(\alpha)\rm alt}(1,1,1,1)$
(written for simplicity in terms of $\lambda = 1+\alpha$) are
\smallskip
\begin{eqnarray}
%% (J_{n,k})_{0 \le n,k \le 7} \;=\;
\Scale[0.68]{
\begin{bmatrix}
1 &   &   &   &   &   &   &   \\
0  & 1 &   &   &   &   &   &   \\
 \lambda & 0 & 1 &   &   &   &   &   \\
0  & 2 + 3\lambda & 0 & 1 &   &   &   &   \\
 2\lambda + 3\lambda^2 & 0 & 8 + 6\lambda & 0 & 1 &   &   &   \\
0  & 16 + 30\lambda + 15\lambda^2 & 0 & 20 + 10\lambda & 0 & 1 &   &   \\
 16\lambda + 30\lambda^2 + 15\lambda^3 & 0 & 136 + 150\lambda + 45 \lambda^2 & 0 & 40 + 15\lambda & 0 & 1 &   \\
0  & 272 + 588\lambda + 420\lambda^2 + 105 \lambda^3 & 0 & 616 + 490\lambda + 105 \lambda^2 & 0 & 70 + 21\lambda & 0 & 1 \\
\end{bmatrix}
}
 \hspace*{-1cm}
     \nonumber \\[1mm]
\end{eqnarray}

\subsection{Continued fraction for secant power polynomials}
   \label{subsec.contfrac.secpower}

Following Stieltjes \cite{Stieltjes_1889} and Rogers \cite{Rogers_07},
we can use Rogers' addition formula (Theorem~\ref{thm.rogers} above)
to obtain the continued fraction for the secant power polynomials.
%% \cite[Example~9.10]{Sokal_alg_contfrac}.
{}From the high-school angle-addition formula
\begin{subeqnarray}
   \cos(t+u)
   & = &
   (\cos t)(\cos u) \,-\, (\sin t)(\sin u)
       \\[1mm]
   & = &
   (\cos t)(\cos u) \, [ 1 \,-\, (\tan t)(\tan u) ]
\end{subeqnarray}
we obtain
\be
   [\sec(t+u)]^\lambda
   \;=\;
   (\sec t)^\lambda (\sec u)^\lambda \:
      \sum_{k=0}^\infty \binom{\lambda+k-1}{k} \, (\tan t)^k \, (\tan u)^k
   \;,
\ee
which is of the form \reff{eq.thm.rogers.1}/\reff{eq.thm.rogers.2} with
\be
   \beta_k \:=\; k(\lambda+k-1)
   \,,\qquad
   F_k(t) \:=\: {(\sec t)^\lambda \, (\tan t)^k \over k!}
          \:=\: {t^k \over k!} \,+\, O(t^{k+2})
   \;,
\ee
so that $\mu_k = 0$ and hence $\gamma_k = 0$.
Theorem~\ref{thm.rogers} then implies that the ordinary generating function
of the secant power polynomials is given by the J-fraction
\be
   \sum_{n=0}^\infty E_{2n}(\lambda) \, t^{2n}
   \;=\;
   \cfrac{1}{1 - \cfrac{1 \cdot \lambda t^2}{1 - \cfrac{2(\lambda+1) t^2}{1 - \cfrac{3(\lambda+2) t^2}{1 - \cdots}}}}
   \;.
 \label{eq.contfrac.secpower}
\ee
After renaming $t^2 \to t$, this is actually an S-fraction
with coefficients $\alpha_n = n(\lambda+n-1)$.
So far, this is all well known.
And comparing \reff{eq.contfrac.secpower}
with Theorem~\ref{thm.second.Sfrac},
we see that $E_{2n}(\lambda)$ enumerates
cycle-alternating permutations of $[2n]$
with a weight $\lambda$ for each cycle.

But Theorem~\ref{thm.rogers} gives more:
it tells us that $F_k(t) = (\sec t)^\lambda \, (\tan t)^k / k!$
is the exponential generating function
for the $k$th column of the generalized Jacobi--Rogers matrix
$\sfJ = \big( J_{n,k}(\bbeta,\bgamma) \big)_{\! n,k \ge 0}$
with $\beta_n = n(\lambda+n-1)$ and $\bgamma = \bzero$.
On the other hand, we know from \reff{eq.columnegf.altlaguerre}
that this $F_k(t)$ is precisely the exponential generating function
for alternating Laguerre digraphs with $k$ paths,
$\sfLtilde^{(\alpha)\rm alt}(1,1,1,1)_{n,k}$,
with $\lambda = 1+\alpha$:
\be
   \sfLtilde^{(\alpha)\rm alt}(1,1,1,1)_{n,k}
   \;=\;
   \Bigl[ {t^n \over n!} \Bigr] \:  {(\sec t)^{1+\alpha} \, (\tan t)^k \over k!}
   \;.
\ee
We have therefore shown:

\begin{proposition}
   \label{prop.JR.cycle}
The generalized Jacobi--Rogers polynomial $J_{n,k}(\bbeta,\bgamma)$
with $\beta_n = n(n+\alpha)$ and $\bgamma = \bzero$
enumerates alternating Laguerre digraphs on $[n]$ with $k$ paths
with a weight $1+\alpha$ for each cycle.
\end{proposition}

By extracting the submatrices $(2n,2k)$ and $(2n+1,2k+1)$
using \reff{eq.Snk.Jnk},
we can equivalently obtain the generalized Stieltjes--Rogers polynomials
$S_{n,k}(\balpha)$ and $S'_{n,k}(\balpha)$ with $\alpha_n = n(n+\alpha)$:
\begin{subeqnarray}
   S_{n,k}(\balpha)
   & = &
   \Bigl[ {t^{2n} \over (2n)!} \Bigr] \:
      {(\sec t)^{1+\alpha} \, (\tan t)^{2k} \over (2k)!}
           \\[3mm]
   S'_{n,k}(\balpha)
   & = &
   \Bigl[ {t^{2n+1} \over (2n+1)!} \Bigr] \:
      {(\sec t)^{1+\alpha} \, (\tan t)^{2k+1} \over (2k+1)!}
\end{subeqnarray}

\begin{corollary}
   \label{cor.SR.cycle}
The generalized Stieltjes--Rogers polynomial $S_{n,k}(\balpha)$
[resp.~$S'_{n,k}(\balpha)$] with $\alpha_n = n(n+\alpha)$
enumerates alternating Laguerre digraphs on $[2n]$ with $2k$ paths
[resp.~on $[2n+1]$ with $2k+1$ paths]
with a weight $1+\alpha$ for each cycle.
\end{corollary}

\subsection{Generalized Stieltjes--Rogers polynomials for the second master S-fraction}
   \label{subsec.contfrac.secondmaster}

We now propose to show how the generalized Stieltjes--Rogers polynomials
of the first and second kinds for the second master S-fraction
for cycle-alternating permutations,
\reff{eq.thm.permutations.Stype.final1.cycle-alternating.second}/%
\reff{eq.thm.permutations.Stype.final1.cycle-alternating.weights.second},
can be understood as generating polynomials
for alternating Laguerre digraphs.
This will lead to a far-reaching generalization of
Proposition~\ref{prop.JR.cycle} and Corollary~\ref{cor.SR.cycle},
and at the same time of \cite[Theorem~2.23]{Sokal-Zeng_masterpoly}.

To do this, it is convenient to use a different convention 
about the {\em boundary conditions}\/ at the two ends of a path
of a Laguerre digraph.
Namely, we will here use $\infty$--$\infty$ (rather than $0$--$0$)
boundary conditions:
that is, we extend the Laguerre digraph $G$ on the vertex set $[n]$
to a digraph $\widehat{G}$ on the vertex set $[n] \cup \{\infty\}$
by decreeing that any vertex $i \in [n]$
that has in-degree (resp.~out-degree) 0 in $G$
will receive an incoming (resp.~outgoing) edge
from (resp.~to) the vertex $\infty$.
However, our definitions of peak, valley, double ascent, double descent
and fixed point will be the same as before.
Note that, because of the $\infty$--$\infty$ boundary conditions,
an isolated vertex is always a valley,
the initial vertex of a path is always a valley or double descent,
and the final vertex of a path is always a valley or double ascent;
moreover, each path contains at least one valley.
We also recall that, for a given vertex $i\in [n]$, 
we use $p(i)$ to denote its predecessor in $\widehat{G}$,
and $s(i)$ to denote its successor.
Later we will explain why we have found it convenient to use
$\infty$--$\infty$ boundary conditions.

With these boundary conditions, an alternating Laguerre digraph
will still have cycles of even length with 
cycle peaks and cycle valleys in alternation,
and paths of odd length consisting of peaks and valleys in alternation.
However, the paths will start and end with a valley:
$\infty > x_1 < x_2 > x_3 < \ldots > x_{2m+1} < \infty$.
The paths are now alternating (up-down) permutations
of odd length written in word form:
$x_1 < x_2 > x_3 < \ldots > x_{2m+1}$.
We write $\LD^{\rm alt\infty}_{n,k}$ for the set of
alternating Laguerre digraphs on the vertex set $[n]$ with $k$ paths
using $\infty$--$\infty$ boundary conditions.
Of course, there is a natural bijection between
$\LD^{\rm alt}_{n,k}$ and $\LD^{\rm alt\infty}_{n,k}$,
by mapping vertices $i \mapsto n+1-i$.

%%We define an \textbfit{alternating Laguerre digraph}
%%to be a Laguerre digraph in which there are no
%%double ascents, double descents or fixed points.
%%It follows that all the cycles are of even length
%%and consist of cycle peaks and cycle valleys in alternation.
%%Furthermore, since we are using 0--0 boundary conditions,
%%all the paths are of odd length and consist of peaks and valleys
%%in alternation, starting and ending with a peak:
%%$0 < x_1 > x_2 < x_3 > \ldots < x_{2m+1} > 0$.
%%That is, the paths are simply alternating (down-up) permutations
%%\cite{Stanley_10} of odd length written in word form:
%%$x_1 > x_2 < x_3 > \ldots < x_{2m+1}$.
%%We write $\LD^{\rm alt}_{n,k}$ for the set of
%%alternating Laguerre digraphs on the vertex set $[n]$ with $k$ paths.
%%Since each path has odd length and each cycle has even length,
%%it follows that $\LD^{\rm alt}_{n,k} \neq \emptyset$ only if $n+k$ is even.
%%We have the obvious bijection $\LD^{\rm alt}_{n,0} \simeq \Sym^{\rm ca}_n$.
%%In this way, alternating Laguerre digraphs generalize
%%cycle-alternating permutations, to which they reduce when $k=0$.

Now we shall define some new statistics on Laguerre digraphs
which extend the crossing and nesting statistics on permutations
that were defined in Section~\ref{subsec.statistics.2}.
We say that a quadruplet $1\leq i < j < k < l\leq \infty$ forms a
\begin{itemize}
   \item {\em lower crossing}\/ (lcross) if $p(i) = k$ and $p(j) = l$;
   \item {\em lower nesting}\/  (lnest)  if $p(j)=k$ and $p(i) = l$.
\end{itemize}
Notice that in such a quadruplet,
the vertex $k$ must be either a peak or a double descent.
When $l < \infty$, these conditions could equivalently be written
in terms of the successor function:
a quadruplet is a lower crossing if $i = s(k)$ and $j = s(l)$,
and a lower nesting if $j = s(k)$ and $i = s(l)$.
But we have written them in terms of the predecessor function
in order to handle smoothly the case $l = \infty$,
i.e.~the case when $j$ (resp.~$i$) has in-degree 0 in $G$.
This definition now allows us to define the
index-refined {\em lower}\/ crossing and nesting statistics
for Laguerre digraphs,
as a generalization of (\ref{def.ucrossnestjk}c,d):
\begin{subeqnarray}
   \lcross(k,G)
   & = &
	\#\{ 1\leq i<j<k<l\leq \infty \colon\: p(i) = k \hbox{ and } p(j) = l \}
         \\[2mm]
   \lnest(k,G)
   & = &
   \#\{ 1\leq i<j<k<l\leq \infty \colon\: p(j)=k \hbox{ and } p(i) = l \}
 \label{def.lcrossnestAltLaguerre}
\end{subeqnarray}
Note that because these statistics are defined in terms
of the index in {\em third}\/ position (namely, $k$),
we certainly have $k < \infty$ (since $k \in [n]$).
When the Laguerre digraph $G$ has no paths
and hence arises as the digraph of a permutation $\sigma$,
this definition manifestly coincides with the definitions
(\ref{def.ucrossnestjk}c,d).

However, things are going to be more subtle for the
{\em upper}\/ crossings and nestings,
because these statistics will be defined
[cf.~(\ref{def.ucrossnestjk}a,b)]
in terms of the index in {\em second}\/ position (namely, $j$),
and there is then no guarantee that $k < \infty$.
Rather, we need to reckon with the possibility
that $i < j < k = l = \infty$,
and this means that (at least in this case)
we will be unable to distinguish between crossings and nestings.
We therefore define a new statistic $\ulev$
that plays the role of $\ucross + \unest$.
That is, for $j\in [n]$ we define
\be
    \ulev(j,G)
    \;=\;
    \#\{ i<j \colon\: j < s(i) \le \infty \hbox{ and }
                                    j < s(j) \le \infty \}
   \;.
 \label{def.ulev}
\ee
We shall say that $\ulev(j,G)$ is the {\em upper level}\/ of the vertex $j$.
Note that $\ulev(j,G)$ can be nonzero only when
$j$ is either a valley or a double ascent.
When the Laguerre digraph $G$ has no paths
and hence arises as the digraph of a permutation $\sigma$,
we have
\be
  \ulev(j,G) \;=\; \ucross(j,\sigma) \,+\, \unest(j,\sigma)
  \;.
\label{eq.upper.lev.crossnest}
\ee

We can also define analogously a {\em lower level}\/ for the vertex $k$:
\be
    \llev(k,G)
    \;=\;
    \#\{ m<k \colon\: k < p(m) \le \infty \hbox{ and }
                                    s(k) < k \}
   \;.
\label{def.llev}
\ee
But this definition is superfluous, because it is concocted to satisfy
\be
\llev(k,G) \; = \; \lcross(k,G) \, + \, \lnest(k,G)
   \;;
\label{eq.lower.lev.crossnest}
\ee
to see this it suffices to look at the two cases
$m=i$ (nesting) and $m=j$ (crossing).

We now restrict attention to {\em alternating}\/ Laguerre digraphs,
and define a generating polynomial in two infinite families of indeterminates
$\bsfa = (\sfa_{\ell})_{\ell\geq 0}$
and $\bsfb = (\sfb_{\ell,\ell'})_{\ell,\ell'\geq 0}$
that uses the crossing, nesting and level statistics:
\be
   \widetilde{Q}_{n,k}(\bsfa,\bsfb,\lambda)
   \;=\;
   \sum_{G \in \LD^{\rm alt\infty}_{n,k}}
   \;\:\lambda^{\cyc(G)}
   \prod\limits_{i \in {\rm Val}(G)}  \! \sfa_{\ulev(i,G)}
   \prod\limits_{i \in {\rm Peak}(G)} \!\! \sfb_{\lcross(i,G),\,\lnest(i,G)}
 \label{def.Qnhat.altLaguerre}
\ee
where ${\rm Val}(G)$ denotes the set of valleys of $G$
and ${\rm Peak}(G)$ denotes the set of peaks of $G$.
For $k=0$, we have $\LD^{\rm alt}_{2n,0} \simeq \Sym_{2n}^{\rm ca}$,
i.e.~an alternating Laguerre digraph $G$ with no paths
is the digraph of a cycle-alternating permutation $\sigma$,
and we have
\be
  \widetilde{Q}_{2n,0} \;=\; \widehat{Q}_n
\ee
where $\widehat{Q}_n$ is the second master polynomial 
for cycle-alternating permutations 
defined in \reff{def.Qnhat.secondmaster.cycle-alternating}
\cite[(2.137)]{Sokal-Zeng_masterpoly}.
We therefore call $\widetilde{Q}_{n,k}(\bsfa,\bsfb,\lambda)$
the \textbfit{second master polynomial for alternating Laguerre digraphs}.

Recall that there is a natural bijection between
$\LD^{\rm alt\infty}_{n,k}$ and $\LD^{\rm alt}_{n,k}$
by mapping vertices $i \mapsto n+1-i$.
This bijection sends cycle peaks to cycle valleys, path peaks to path valleys,
and vice-versa.
Thus, under the specializations
$\sfa_{\ell} = y_v$, $\sfb_{\ell,\ell'} = y_p$ and $\lambda = 1+\alpha$,
the polynomial $\widetilde{Q}_{n,k}(\bsfa,\bsfb,\lambda)$
is equal to the polynomial 
$\sfLtilde^{(\alpha)\rm alt}(\yv,\yp, \yv, \yp)_{n,k}$
obtained by substituting 
$\yp = \zp$
to $\yv$
and $\yv = \zv$
to $\yp$
in \reff{def.coeffmat.gen.second.gen.alt}.

The main technical result of this subsection
is that the polynomials $\widetilde{Q}_{n,k}(\bsfa,\bsfb,\lambda)$ 
satisfy the following recurrence:

\begin{proposition}[Recurrence for second master polynomial for alternating Laguerre digraphs]
   \label{prop.Qnkhat.recurrence}
The polynomials $\widetilde{Q}_{n,k}(\bsfa,\bsfb,\lambda)$ defined in
\reff{def.Qnhat.altLaguerre} satisfy the recurrence
\be
\widetilde{Q}_{n+1,k} \; = \; 
	\sfa_{k-1}\,\widetilde{Q}_{n,k-1} \:+\:
	(\lambda+k)\left( \sum_{i=0}^{k} \sfb_{i,k-i} \right)\,
           \widetilde{Q}_{n,k+1}\;.
\label{eq.Qnhat.recurrence}
\ee
\end{proposition}

\proof
In order to prove the recurrence \reff{eq.Qnhat.recurrence},
we first start with an alternating Laguerre digraph
$G\in \LD^{\rm alt\infty}_{n+1,k}$.
We then look at the status of the vertex $n+1$,
and for each possibility 
we consider the induced subgraph on the vertex set $[n]$, call it $G'$.
(Clearly $G'$ is a Laguerre digraph;
 we will also see that it is alternating,
 but that takes a bit more work to prove.)
For each situation, we then count the number of the ways 
in which $n+1$ could be inserted into $G'$,
and evaluate what this does to the weights.

Let us begin by observing that, for every $i \in [n]$,
\begin{subeqnarray}
   \lcross(i,G)  & = & \lcross(i,G') \\[2mm]
   \lnest(i,G)   & = & \lnest(i,G') \\[2mm]
   \ulev(i,G)    & = & \ulev(i,G')
\end{subeqnarray}
This is because, if the predecessor (resp.~successor) of $i$ in $G$ is $n+1$,
then its predecessor (resp.~successor) in $G'$ is $\infty$;
and in the inequalities defining $\lcross,\lnest,\ulev$,
the values $n+1$ and $\infty$ play identical roles.
(That is why we have chosen here to use
 $\infty$--$\infty$ boundary conditions.)
It follows that all vertices $i \in [n]$
get the same weight $\sfa_\ell$ or $\sfb_{\ell,\ell'}$
in $G$ as they do in $G'$.
Therefore, we only need to worry about
the weight associated to the vertex $n+1$,
as well as the possible change in the number of cycles.

We begin by considering the case when $n+1$ is a valley.
For this to happen, we must have $p(n+1) = s(n+1) = \infty$;
in other words, $n+1$ is an isolated vertex in $G$.
Removing $n+1$, we see that
the resulting digraph $G'$ has $k-1$ paths and is alternating.
On the other hand if we begin with an alternating Laguerre digraph 
$G'\in \LD^{\rm alt\infty}_{n,k-1}$,
there is exactly one way in which the vertex $n+1$ can be inserted into $G'$
as an isolated vertex.
Note that, in $G$, the vertex $n+1$
has upper level $\ulev(n+1,G) = k-1$:
this is because, in \reff{def.ulev} with $j=n+1$,
we have $s(j) = \infty$ (because $j=n+1$ is a valley),
and the contributing vertices $i$ are those that have $s(i) = \infty$;
this happens when (and only when) $i$ ($\neq j$) is the end vertex of a path,
and there are $k-1$ such paths.
Therefore, the vertex $n+1$ gets a weight $\sfa_{k-1}$
in its contribution to $\widetilde{Q}_{n+1,k}$.
Therefore, the cases in which $n+1$ is a valley
contribute the summand $\sfa_{k-1}\,\widetilde{Q}_{n,k-1}$.

When $n+1$ is a peak, its predecessor $p(n+1)$ and successor $s(n+1)$
lie in $[n]$;
therefore, it cannot be the initial or final vertex of a path in $G$.
If $n+1$ is a peak belonging to a cycle of $G$,
then removing it turns its component into a path 
and increases the number of paths by $1$.
If $n+1$ is a peak belonging to a path of $G$
(and hence an internal vertex of that path),
removing it causes that path to be divided into two nonempty paths,
and hence again increases the number of paths by $1$.
In both cases the resulting Laguerre digraph $G'$
will be alternating,
because the vertices that had $n+1$ as predecessor or successor in $G$
will now have $\infty$ as their predecessor or successor in $G'$,
maintaining the required inequalities.
Hence $G'\in \LD^{\rm alt\infty}_{n,k+1}$.

We now begin with an alternating Laguerre digraph
$G'\in \LD^{\rm alt\infty}_{n,k+1}$,
and count the number of ways in which $n+1$ can be inserted as peak.
To do this, we need to count the various possibilities
in which the predecessor and successor of $n+1$ can be chosen.
The successor of $n+1$ can be chosen from any of the 
$k+1$ initial vertices of the paths.
Let $1\leq f_0<f_1<\ldots< f_k \leq n$ be the initial vertices of the paths
of $G'$, ordered in increasing order. 
If $s(n+1) = f_i$, any resulting alternating Laguerre digraph $G$ will have
$p(f_j) = \infty$ for all $j \neq i$.
In particular, for all $j < i$,
the quadruplets $f_j < f_i < n+1 < \infty$
will be lower nestings;
and these will be the only lower nestings with third index $n+1$.
Similarly, for all $j > i$,
the quadruplets $f_i < f_j < n+1 < \infty$
will be lower crossings;
and these will be the only lower crossings with third index $n+1$.
Therefore, we will have $\lnest(n+1, G) = i$ and $\lcross(n+1, G) = k-i$.
Thus, whenever $s(n+1) = f_i$ ($0 \leq i\leq k$),
the vertex $n+1$ gets weight $\sfb_{i,k-i}$ in $G$.

The predecessor of $n+1$ can be chosen from any of the
$k+1$ final vertices of the paths.
However, for each choice of successor $f_i$,
there is exactly one choice of final vertex that forms a cycle
(namely, when the successor and predecessor vertices of $n+1$
 are the initial and final vertices of the same path);
in this situation we get an additional weight $\lambda$.
In the other $k$ cases we get a weight $1$.
Thus, the case when $n+1$ is a peak contributes the summand
$(\lambda+k)\left( \sum_{i=0}^{k} \sfb_{i,k-i} \right)\, \widetilde{Q}_{n,k+1}$.
\qed

The recurrence \reff{eq.Qnhat.recurrence}
combined with Proposition~\ref{prop.JR.bothkinds} yields:

\begin{theorem}[Generalized Jacobi-Rogers polynomials]
   \label{thm.JRmat}
The generalized Jacobi-Rogers polynomial $J_{n,k}(\balpha,\bzero)$ with weights
\be
   \alpha_n \;=\;
   (\lambda + n-1) \sfa_{n-1}
   \biggl( \sum_{\ell=0}^{n-1} \sfb_{\ell,n-1-\ell} \biggr) 
  \label{eq.thm.JRpoly}
\ee
[cf.~\reff{eq.thm.permutations.Stype.final1.cycle-alternating.weights.second}]
satisfies
\be
     \widetilde{Q}_{n,k}(\bsfa, \bsfb, \lambda) 
     \;=\;
     \left(\prod\limits_{i=0}^{k-1}\sfa_i \right) J_{n,k}(\balpha, \bzero) 
  \label{eq.thm.JRmat}
\ee
where $\widetilde{Q}_{n,k}(\bsfa, \bsfb, \lambda)$ 
is defined in \reff{def.Qnhat.altLaguerre}.
\end{theorem}

Since the J-fraction in Theorem~\ref{thm.JRmat} has $\bgamma = \bzero$,
the even and odd submatrices of $\sfJ(\balpha, \bzero)$
can be interpreted by \reff{eq.Snk.Jnk}
as being the generalized Stieltjes--Rogers polynomials
of the first and second kinds for the second master S-fraction
for cycle-alternating permutations,
\reff{eq.thm.permutations.Stype.final1.cycle-alternating.second}/%
\reff{eq.thm.permutations.Stype.final1.cycle-alternating.weights.second}:

\begin{corollary}[Generalized Stieltjes-Rogers polynomials]
   \label{cor.JRmat}
With $\balpha$ given by \reff{eq.thm.JRpoly}, we have
\begin{subeqnarray}
   \widetilde{Q}_{2n,2k}(\bsfa, \bsfb, \lambda)
   & = &
   \left(\prod\limits_{i=0}^{2k-1}\sfa_i \right) S_{n,k}(\balpha)
      \\[2mm]
   \widetilde{Q}_{2n+1,2k+1}(\bsfa, \bsfb, \lambda)
   & = &
   \left(\prod\limits_{i=0}^{2k}\sfa_i \right) S'_{n,k}(\balpha)
 \label{eq.cor.JRmat}
\end{subeqnarray}
\end{corollary}

\subsection{Specialization to the case of
   Theorems~\ref{thm.second.Sfrac} and \ref{thm.second.Sfrac.pqgen}}
  \label{subsec.contfrac.specialization}

By specializing Corollary~\ref{cor.JRmat},
we can obtain the generalized Stieltjes-Rogers polynomials
corresponding to the S-fractions of
Theorems~\ref{thm.second.Sfrac} and \ref{thm.second.Sfrac.pqgen}.
To do this, we will first need to extend the concepts of records and antirecords
from permutations to Laguerre digraphs.
This extension will be motivated by \cite[Lemma~2.10]{Sokal-Zeng_masterpoly},
which we have recalled in Lemma~\ref{lemma.SZ.nesting}.

%\begin{lemma} {\hspace*{-2mm}\bf \cite[Lemma~2.10]{Sokal-Zeng_masterpoly}\ }
%   \label{lemma.SZ.nesting}
%Let $\sigma$ be a permutation.
%\vspace*{-2mm}
%\begin{itemize}
%   \item[(a)]  If $i$ is a cycle valley or cycle double rise,
%       then $i$ is a record if and only if $\unest(i,\sigma) = 0$;
%       and in this case it is an exclusive record.
%   \item[(b)]  If $i$ is a cycle peak or cycle double fall,
%       then $i$ is an antirecord if and only if $\lnest(i,\sigma) = 0$;
%       and in this case it is an exclusive antirecord.
%\end{itemize}
%\end{lemma}

We now extend the concepts of record and antirecord from permutations
to arbitrary words in a totally ordered alphabet:

\begin{definition}
Let $w_1 \ldots w_n$ be a word in a totally ordered alphabet.
Then:
\begin{itemize}
   \item[(a)] $i$ is a {\em record position}\/
      (and $w_i$ is a {\em record value}\/)
      in case, for all $j < i$, we have $w_j < w_i$.
   \item[(b)] $i$ is an {\em antirecord position}\/
      (and $w_i$ is an {\em antirecord value}\/)
      in case, for all $j > i$, we have $w_j > w_i$.
\end{itemize}
\end{definition}

When the word $w_1 \ldots w_n$ is a permutation of $[n]$,
these concepts are connected as follows:

\begin{lemma}
   \label{lemma.records}
Let $\sigma$ be a permutation of $[n]$,
and let $\sigma^{-1}$ be the inverse permutation.  Then:
\begin{itemize}
   \item[(a)] $i$ is a record position in the word
      $\sigma(1) \,\ldots\, \sigma(n)$
      if and only if $i$ is an antirecord value in the word
      $\sinv(1) \,\ldots\, \sinv(n)$.
   \item[(b)] $i$ is an antirecord position in the word
      $\sigma(1) \,\ldots\, \sigma(n)$
      if and only if $i$ is a record value in the word
      $\sinv(1) \,\ldots\, \sinv(n)$.
\end{itemize}
\end{lemma}

\proof
(a) By definition, $i$ is a record position in the word
$\sigma(1) \,\ldots\, \sigma(n)$ in case the following holds true:
whenever $j < i$, we have $\sigma(j) < \sigma(i)$.
Similarly, $i$ is an antirecord value in the word
$\sinv(1) \,\ldots\, \sinv(n)$ in case the following holds true:
whenever $k > \sigma(i)$, we have $\sinv(k) > i$.
Defining $j = \sinv(k)$, this latter statement can be rewritten as:
whenever $\sigma(j) > \sigma(i)$, we have $j > i$.
For each pair $i \neq j$,
the third statement is the contrapositive of the first.

(b) The proof is similar.
\qed

Using Lemma~\ref{lemma.records}, we can rewrite Lemma~\ref{lemma.SZ.nesting}
as follows:

\begin{lemma}
   \label{lemma.SZ.nesting.bis}
Let $\sigma$ be a permutation.
\vspace*{-2mm}
\begin{itemize}
   \item[(a)]  If $i$ is a cycle valley or cycle double rise,
       then $i$ is an antirecord value in the word
       $\sinv(1) \,\ldots\, \sinv(n)$
       if and only if $\unest(i,\sigma) = 0$.
   \item[(b)]  If $i$ is a cycle peak or cycle double fall,
       then $i$ is a record value in the word
       $\sinv(1) \,\ldots\, \sinv(n)$
       if and only if $\lnest(i,\sigma) = 0$.
\end{itemize}
\end{lemma}

We now generalize Lemma~\ref{lemma.SZ.nesting.bis}(b)
from permutations to Laguerre digraphs.
To do this, we first observe that to any Laguerre digraph $G$
on the vertex set $[n]$,
we can associate a \textbfit{predecessor word} $p(1) \,\ldots\, p(n)$
and a \textbfit{successor word} $s(1) \,\ldots\, s(n)$,
each taking values in the totally ordered alphabet $[n] \cup \{\infty\}$.
Of course, if $G$ is the digraph associated to a permutation $\sigma$,
then the predecessor word is $\sinv(1) \,\ldots\, \sinv(n)$
and the successor word is $\sigma(1) \,\ldots\, \sigma(n)$,
which of course take values in $[n]$.
We then have the following generalization of
Lemma~\ref{lemma.SZ.nesting.bis}(b):

\begin{lemma}
   \label{lemma.laguerre.nesting}
Let $G$ be a Laguerre digraph on the vertex set $[n]$.
If $k$ is a peak or double descent,
then $k$ is a record value in the predecessor word $p(1) \,\ldots\, p(n)$
if and only if $\lnest(k,G) = 0$.
\end{lemma}

\proof
% If the letter $k$ occurs in the predecessor word $p(1) \,\ldots\, p(n)$,
% we must of course have $k = p(j)$ for some $j \in [n]$;
% and when $k < \infty$, there can be only one such $j$
% --- namely, $j = s(k)$ ---
% since $k$ has out-degree at most~1 in $G$.
% 
% Since $k$ is here a peak or double descent, we have $s(k) < k < \infty$.
%
Since $k$ is a peak or double descent, we have $s(k) < k < \infty$.
Therefore, $k$ has out-degree~1 in $G$.
It follows that the letter $k$ occurs exactly once
in the predecessor word $p(1) \,\ldots\, p(n)$:
namely, $k = p(j)$ if and only if $j = s(k)$.

By definition, $k$ is a record value
in the predecessor word $p(1) \,\ldots\, p(n)$
in case, for all $i < s(k)$, we have $p(i) < k$.
Equivalently, $k$ is a record value 
in the predecessor word $p(1) \,\ldots\, p(n)$
in case there does not exist $i < s(k)$ such that $p(i) > k$
(equality $p(i) = k$ is excluded because it would imply $i = s(k)$).
But, by the definition (\ref{def.lcrossnestAltLaguerre}b) of $\lnest$,
this is precisely the statement that $\lnest(k,G) = 0$.
\qed

With these preliminaries, we are now able to state
the generalized Stieltjes--Rogers polynomials
corresponding to the S-fractions of
Theorems~\ref{thm.second.Sfrac} and \ref{thm.second.Sfrac.pqgen}.
Let us start with Theorem~\ref{thm.second.Sfrac},
which concerns the polynomials \reff{def.Qnhat.cycle-alternating.evenodd}
specialized to $\ye=\ve$ and $\yo=\vo$:
this means that for cycle peaks (but not for cycle valleys)
we are distinguishing between those that are antirecords
(i.e.~antirecord {\em positions}\/)
and those that are not.
By Lemma~\ref{lemma.records}(b),
$i$ is an antirecord position in $\sigma$
if and only if $i$ is a record value in $\sinv$.
We then generalize the latter condition
to Laguerre digraphs $G$ on the vertex set $[n]$ by defining
\begin{subeqnarray}
   & &
   \Arecpeak^\star(G)
   \;\eqdef\;
   \{ i \in [n] \colon\: i \in \Peak(G) \hbox{ and }
          \nonumber \\
   & & \qquad
      \hbox{$i$ is a record value in the predecessor word
            $p(1) \,\ldots\, p(n)$}  \},
         \\[2mm]
   & &
   \Nrpeak^\star(G)
   \;\eqdef\;
   \{ i \in [n] \colon\: i \in \Peak(G) \hbox{ and }
          \nonumber \\
   & & \qquad
      \hbox{$i$ is not a record value in the predecessor word
            $p(1) \,\ldots\, p(n)$}  \},
      \qquad
\end{subeqnarray}
where we recall that we have defined
\begin{subeqnarray}
   \Peak(G)
   & \eqdef &
   \{ i \in [n] \colon\: \hbox{$i$ is a peak in $G$ using
      $\infty$--$\infty$ boundary conditions} \},
    \qquad\qquad
        \\[2mm]
   \Val(G)
   & \eqdef &
   \{ i \in [n] \colon\: \hbox{$i$ is a valley in $G$ using
      $\infty$--$\infty$ boundary conditions} \}.
    \qquad\qquad
\end{subeqnarray}
We then define
\begin{subeqnarray}
   \eareccpeakeven^\star(G) & \eqdef & |\Arecpeak^\star(G) \cap \Even|
       \\[2mm]
   \nrcpeakeven^\star(G) & \eqdef & |\Nrpeak^\star(G) \cap \Even|
       \\[2mm]
   \cvaleven^\star(G)  & \eqdef & |\Val(G) \cap \Even|
\end{subeqnarray}
and likewise for the odd ones.
(The names for these quantities are chosen to coincide with
the previous ones when $G$ is the digraph arising from a permutation,
even though the names may seem a bit bizarre when compared with
the actual definitions;
we have used a star to stress this change of emphasis.)
Having done this, we then extend the definition
\reff{def.Qnhat.cycle-alternating.evenodd},
specialized to $\ye=\ve$ and $\yo=\vo$,
from permutations to Laguerre digraphs:
\begin{eqnarray}
   & & \hspace*{-6mm}
   \widehat{Q}_{n,k}(\xe,\ye,\ue,\xo,\yo,\uo,\lambda)
   \;=\;
        \nonumber \\[3mm]
   & & \qquad
   \sum_{G \in \LD^{\rm alt\infty}_{n,k}}
   \xe^{\eareccpeakeven^\star(G)}
   \ye^{\cvaleven^\star(G)}
   \ue^{\nrcpeakeven^\star(G)}
          \:\times
       \qquad\qquad
       \nonumber \\[0mm]
   & & \qquad\qquad\qquad\!
   \xo^{\eareccpeakodd^\star(G)}
   \yo^{\cvalodd^\star(G)}
   \uo^{\nrcpeakodd^\star(G)}
   \lambda^{\cyc(G)}
   \;.
 \label{def.Qnhat.cycle-alternating.evenodd.G}
\end{eqnarray}
We can then state:

\begin{theorem}[Generalized Stieltjes--Rogers polynomials for
   Theorem~\ref{thm.second.Sfrac}]
 \label{thm.second.Sfrac.genSR}
\hfill\break
With $\balpha$ given by
\reff{def.weights.perm.Stype.cycle-alternating.evenodd.second},
we have
\begin{subeqnarray}
   \widehat{Q}_{2n,2k}(\xe,\ye,\ue,\xo,\yo,\uo,\lambda)
   & = &
   \ye^k \yo^k \:
   S_{n,k}(\balpha)
       \\[2mm]
   \widehat{Q}_{2n+1,2k+1}(\xe,\ye,\ue,\xo,\yo,\uo,\lambda)
   & = &
   \ye^k \yo^{k+1} \:
   S'_{n,k}(\balpha)
\end{subeqnarray}
\end{theorem}

As preparation for this proof, we need a generalization
of Lemma~\ref{lemma.cycle-alt} from cycle-alternating permutations
to alternating Laguerre digraphs:

\begin{lemma}
   \label{lemma.cycle-alt.G}
If $G$ is an alternating Laguerre digraph on the vertex set $[N]$, then
\begin{subeqnarray}
   \hbox{\rm valleys:}
   & \!\!  &
   \ulev(i,G)
   \;\equiv\;
   i-1 \: \pmod{2}
         \\[1mm]
   \hbox{\rm peaks:}
   & \!\!  &
   \llev(i,G)
   \;\equiv\;
   i \: \pmod{2}
 \label{eq.lemma.cycle-alt.G}
\end{subeqnarray}
for all $i \in [N]$.
Here $\ulev$ and $\llev$ were defined in
\reff{def.ulev} and \reff{def.llev}, respectively.
\end{lemma}

We introduce some new notation before providing the proof.
Given a Laguerre digraph $G$ on the vertex set $[n]$
and a set $S\subseteq [n]$,
let $\left. G\right|_S$ denote the induced subgraph
on the vertex set $S$ (which is, of course, also a Laguerre digraph).
Also, let $h^{S}$ denote the number of paths in the graph
$\left. G\right|_S$.

\proofof{Lemma~\ref{lemma.cycle-alt.G}}
To begin with, let $G$ be any Laguerre digraph, not necessarily alternating.
Let $j$ be a valley in $G$.
Then, in the graph $\left. G\right|_{[j]}$,
$j$ is an isolated vertex,
i.e., $s(j) = p(j) = \infty$.
From the definition \reff{def.ulev} we have
\begin{eqnarray}
    \ulev(j, \left. G\right|_{[j]})
    &=&
    \#\{ i<j \colon\: s(i) = \infty \}\nonumber\\
    &=& h^{[j]} -1
   \;.\label{eq.ulev.Gi}
\end{eqnarray}
where the first equality holds because $j$ is the largest-numbered vertex
in $\left. G\right|_{[j]}$, so that $s(i) > j$ implies $s(i) = \infty$;
and the second equality follows by noticing
that the vertices $i$ with $s(i) = \infty$ correspond
to the final vertices of the paths in $\left. G\right|_{[j]}$
different from $j$.

Similarly, when $k$ is a peak in $G$,
it follows from \reff{def.llev} that
\begin{eqnarray}
    \llev(k, \left. G\right|_{[k]})
    &=&
    \#\{m<k \colon\: p(m) = \infty \}\nonumber\\
	&=& h^{[k]} \;,
	\label{eq.llev.Gi}
\end{eqnarray}
where the vertices $m$ with $p(m) = \infty$ correspond
to the initial vertices of the paths in $\left. G\right|_{[k]}$.

Next, observe that
\begin{subeqnarray}
	h^{[i]} &=& h^{[i-1]} + 1 \quad \hbox{when  $i$ is a valley} \\
	h^{[i]} &=& h^{[i-1]} - 1 \quad \hbox{when  $i$ is a peak}
\label{eq.hi.path}
\end{subeqnarray}
The first formula is true because a valley $i$
is an isolated vertex in $\left. G\right|_{[i]}$,
and removing it decreases the number of paths by $1$.
On the other hand, if $i$ is a peak,
it has both a predecessor and a successor in $\left. G\right|_{[i]}$;
so if $i$ is a cycle peak, 
removing it changes a cycle to a path,
whereas if $i$ is a path peak,
removing it splits the path into two paths;
either way the number of paths is increased by 1.

Now let $G$ be an {\em alternating}\/ Laguerre digraph
on the vertex set $[N]$, with $k$ paths.
Consider the lattice path $\omega$ starting at $(0,0)$ and ending at $(N,k)$
with steps $s_i = (1,1)$ when $i$ is a valley
and $s_i = (1,-1)$ when $i$ is a peak.
By definition, the heights $h_i$ in the path $\omega$
satisfy the same recurrence \reff{eq.hi.path}
satisfied by the $h^{[i]}$;
and of course $h_0 = h^{[0]} = 0$ (since $[0] = \emptyset$).
It follows that $h_i = h^{[i]} \ge 0$,
and hence that $\omega$ is a partial Dyck path.
This also implies that
\be
h^{[i]} \;\equiv\; i \pmod{2}.
\label{eq.path.heights}
\ee
The proof of \reff{eq.lemma.cycle-alt.G} is completed by combining 
\reff{eq.ulev.Gi}, \reff{eq.llev.Gi} and \reff{eq.path.heights}
and noticing that
$\ulev(i,\left. G\right|_{[i]} ) = \ulev(i,G) $
and
$\llev(i,\left. G\right|_{[i]} ) = \llev(i,G) $.
\qed

\proofof{Theorem~\ref{thm.second.Sfrac.genSR}}
We obtain the polynomials
$\widehat{Q}_{n,k}(\xe,\ye,\ue,\xo,\yo,\uo,\lambda)$
defined in
\reff{def.Qnhat.cycle-alternating.evenodd.G}
by using Lemmas~\ref{lemma.laguerre.nesting} and \ref{lemma.cycle-alt.G}
and eq.~\reff{eq.lower.lev.crossnest},
and specializing the polynomials $\widetilde{Q}_{n,k}(\bsfa,\bsfb,\lambda)$
defined in \reff{def.Qnhat.altLaguerre} to
\begin{subeqnarray}
   \sfa_{\ell}
   & = &
   \begin{cases}
       \yo    &  \textrm{if $\ell$ is even}  \\
       \ye    &  \textrm{if $\ell$ is odd}  \\
   \end{cases}
       \\[2mm]
   \sfb_{\ell,\ell'}
   & = &
   \begin{cases}
       \xe    &  \textrm{if $\ell' = 0$ and $\ell+\ell'$ is even}  \\
       \ue    &  \textrm{if $\ell' \ge 1$ and $\ell+\ell'$ is even} \\
       \xo    &  \textrm{if $\ell' = 0$ and $\ell+\ell'$ is odd}  \\
       \uo    &  \textrm{if $\ell' \ge 1$ and $\ell+\ell'$ is odd}
   \end{cases}
\end{subeqnarray}
(This is just \reff{eq.sfasfb.pq} with all the $p,q$ variables set to 1.)
Inserting these into Corollary~\ref{cor.JRmat} 
completes the proof of Theorem~\ref{thm.second.Sfrac.genSR}.
\qed

Next we show the $p,q$-generalization of all this.
We extend the definition \reff{def.Qnhat.cycle-alternating.evenodd.pqgen},
specialized to $\ye=\ve$, $\yo=\vo$, $p_{+1} = q_{+1}$ and $p_{+2}=q_{+2}$,
from permutations to Laguerre digraphs:
\begin{eqnarray}
& & \hspace*{-6mm}
	\widehat{Q}_{n,k}(\xe,\ye,\ue,\xo,\yo,\uo,p_{-1},p_{-2},p_{+1},p_{+2},q_{-1},q_{-2},\lambda)
   \;=\;
        \nonumber \\[4mm]
   & & \qquad
   \sum_{G \in \LD^{\rm alt\infty}_{n,k}}
   \xe^{\eareccpeakeven^\star(G)}
   \ye^{\cvaleven^\star(G)}
   \ue^{\nrcpeakeven^\star(G)}
          \:\times
       \qquad\qquad
       \nonumber \\[0mm]
   & & \qquad\qquad\qquad\!
   \xo^{\eareccpeakodd^\star(G)}
   \yo^{\cvalodd^\star(G)}
   \uo^{\nrcpeakodd^\star(G)}
        \:\times
       \qquad\qquad
       \nonumber \\[2mm]
   & & \;\quad\qquad\qquad\;\:
   p_{-1}^{{\rm lcrosscpeakeven}^\star(G)}
   q_{-1}^{{\rm lnestcpeakeven}^\star(G)}
   p_{+2}^{{\rm ulevcvaleven}^\star(G)}
   \:\times
       \qquad\qquad
       \nonumber \\[2mm]
   & & \;\quad\qquad\qquad\;\:
   p_{-2}^{{\rm lcrosscpeakodd}^\star(G)}
   q_{-2}^{{\rm lnestcpeakodd}^\star(G)}
   p_{+1}^{{\rm ulevcvalodd}^\star(G)}
   \;,
\end{eqnarray}
where 
\be
{\rm lcrosscpeakeven}^\star(G) \;=\; \sum_{k\in \Peak(G) \,\cap \,{\rm Even}} \lcross(k,G)
\ee
and likewise for the others.
We can then state:

\begin{theorem}[Generalized Stieltjes--Rogers polynomials for
   Theorem~\ref{thm.second.Sfrac.pqgen}]
 \label{thm.second.Sfrac.genSR.pqgen}
\hfill\break
With $\balpha$ given by
\reff{def.weights.perm.Stype.cycle-alternating.evenodd.pqgen.second},
we have
\begin{subeqnarray}
	\widehat{Q}_{2n,2k}(\xe,\ye,\ue,\xo,\yo,\uo, p_{-1},q_{-1},p_{+2},p_{-2},q_{-2},p_{+1},\lambda)\nonumber\\
  \; = \;
    p_{+1}^{k(k-1)}  p_{+2}^{k^2} \: \ye^k \yo^k \:
   S_{n,k}(\balpha)\qquad\qquad\quad
       \\[2mm]
	\widehat{Q}_{2n+1,2k+1}(\xe,\ye,\ue,\xo,\yo,\uo, p_{-1},q_{-1},p_{+2},p_{-2},q_{-2},p_{+1},\lambda)\nonumber\\
  \; = \;
  p_{+1}^{k(k+1)}  p_{+2}^{k^2} \: \ye^k \yo^{k+1} \:
   S'_{n,k}(\balpha)\qquad\qquad\quad
\end{subeqnarray}
\end{theorem}

\proof
The proof of Theorem~\ref{thm.second.Sfrac.genSR.pqgen}
is almost identical to that of Theorem~\ref{thm.second.Sfrac.genSR},
but uses the full specialization \reff{eq.sfasfb.pq}.
\qed

\section*{Acknowledgments}
\addcontentsline{toc}{section}{Acknowledgments}

We wish to thank Einar Steingr\'{\i}msson and Jiang Zeng
for helpful discussions.
We also thank the organizers of the
Heilbronn Workshop on Positivity Problems Associated to Permutation
Patterns (Lancaster, UK, June~2022)
for making these discussions possible.

The first author is currently supported by the DIMERS project ANR-18-CE40-0033
funded by Agence Nationale de la Recherche (ANR, France).

\appendix
\section{Nonexistence of J-fractions with polynomial coefficients}

In this appendix we show that it is not possible to obtain nice J-fractions
for polynomials that take account of the record status of both
cycle peaks and cycle valleys and also count the number of cycles.
Nor is it possible to obtain nice J-fractions for analogous polynomials
restricted to alternating cycles.
We show these two facts in the two following subsections.

\subsection{J-fraction for the polynomials \reff{eq.Pn.xylam}}
    \label{app.Pnxylam}

% From Bishal's file ugly_Jfrac.tex

Consider the polynomials
\be
   P_n(x,y,\lambda)
   \;=\;
   \sum_{\sigma \in \Sym^{\rm ca}_{2n}}
       x^{\eareccpeak(\sigma)}
       y^{\ereccval(\sigma)}
        \lambda^{\cyc(\sigma)}
 \label{eq.Pn.xylam.bis}
\ee
that were introduced in \reff{eq.Pn.xylam}.
By the reversal map $i \mapsto 2n+1-i$,
%% the weight $x^{\eareccpeak(\sigma)} y^{\ereccval(\sigma)}$
%% is equivalent to $x^{\ereccval(\sigma)} y^{\eareccpeak(\sigma)}$, so
we have the symmetry $x \leftrightarrow y$:
$P_n(x,y,\lambda) = P_n(y,x,\lambda)$.
The first few polynomials $P_n(x,y,\lambda)$ are
\begin{subeqnarray}
P_0(x,y,\lambda) &=& 1 \\[1mm]
P_1 ( x,y, \lambda ) & = & \lambda x y\\[1mm]
P_2(x,y,\lambda) & = & \lambda x y \, \left[\lambda\,(1+2xy)+ x+y \right] \\[1mm]
P_3(x,y,\lambda) &=&  \lambda x y \, \left[\lambda^2 \, (3 + 7 x y + 5 x^2 y^2)\right. \nonumber\\
&& + \; \lambda \, (2 + 5 x + x^2 + 5 y + 4 x y + 6 x^2 y + y^2 + 6 x y^2)\nonumber \\
&& + \; \left.(3 x + 2 x^2 + 3 y + 4 x y + x^2 y + 2 y^2 + x y^2)\right].
\end{subeqnarray}
We now give the first few J-fraction coefficients for this polynomial sequence: they are
\be
\gamma_0 \;=\; \lambda xy, \quad \beta_1 \; =\; \lambda xy \, \left[\lambda\,(1+xy)+ x+y \right]
\ee
followed by
\be
\gamma_1 \;=\; 
{\scalebox{0.75}{$\dfrac{\lambda^2  (1 + x y) (3 + 2 x y)
\,+ \, \lambda  (2 + 5 x + x^2 + 5 y + 4 x y + 4 x^2 y + y^2 + 4 x y^2) \, +\, (x + y) (3 + 2 x + 2 y + x y)}{\lambda\,(1+xy)+ x+y}$}}
 \;.
\ee
It can then be shown that
\begin{itemize}
\item[(a)] $\gamma_1$ is not a polynomial in $x$ (when $y$ and $\lambda$ are given fixed real values) unless  $\lambda \in \{-1,0,+1\}$ or $y \in  \{-1,+1\}$.

\item[(b)] $\gamma_1$ is not a polynomial in $y$ unless $\lambda \in \{-1,0,+1\}$ or $x \in  \{-1,+1\}$.

\item[(c)] $\gamma_1$ is not a polynomial in $\lambda$ unless $x \in \{-1,+1\}$ or $y \in  \{-1,+1\}$ or $x=-y$.
\end{itemize}

We state the polynomials $\gamma_1$ obtained from the foregoing
specializations:

\begin{center}
\begin{tabular}{c|c}
Specialization & $\gamma_1$\\
\hline
$\lambda=-1$ & $-1 + x + y - 2 x y$ \\
$\lambda=0$ & $3 + 2 x + 2 y + x y$ \\
$\lambda=+1$ & $5 + 3 x + 3 y + 2 x y$ \\
$x=-1$ & $\lambda\, (3-2y) + 1+y$ \\
$x=+1$ & $\lambda\, (3 + 2 y) + 5 + 3 y$\\
$y=-1$ & $\lambda\, (3-2x) + 1+x$ \\
$y=+1$ & $\lambda\, (3 + 2 x) + 5 + 3 x$\\
$x=-y$ & $ \lambda \,(3 - 2 y^2)+2$ \\
\end{tabular}
\end{center}

These give rise to continued fractions as follows:
\medskip

$\bm{\lambda=+1}$.
By \cite[Theorem~2.18]{Sokal-Zeng_masterpoly}
specialized to $\lambda = u_1 = v_1 = 1$, $x_1 = x$, $y_1=y$,
we obtain an S-fraction with
$\alpha_{n} = (x + n-1)(y + n-1)$, 
and hence by contraction
%% \cite[p.~21]{Wall_48} \cite[p.~V-31]{Viennot_83}
\reff{eq.contraction_even.coeffs}
a J-fraction with $\gamma_0 = xy$,
$\gamma_n = (x+2n-1)(y+2n-1)+(x+2n)(y+2n)$ for $n \ge 1$,
$\beta_n = (x+2n-2)(x+2n-1)(y+2n-2)(y+2n-1)$.

\medskip

$\bm{y=+1}$.
By \cite[Theorem~2.21]{Sokal-Zeng_masterpoly}
specialized to $y_1 = u_1 = v_1 =1$, $x_1 = x$,
we obtain an S-fraction with
$\alpha_{n} =  (\lambda + n - 1) (x + n - 1)$.
By contraction this yields a J-fraction with
$\gamma_0 = \lambda x$,
$\gamma_n =  (\lambda + 2n - 1) (x + 2n - 1) + (\lambda + 2n) (x + 2n)$ for $n \ge 1$,
$\beta_n = (\lambda+2n-2)(\lambda+2n-1)(x+2n-2)(x+2n-1)$.

\medskip

$\bm{x=+1}$.
This follows from the case $y=+1$ by using the symmetry
$x \leftrightarrow y$.
% By the reversal map $i \mapsto 2n+1-i$,
% the weight $x^{\eareccpeak(\sigma)}
%        y^{\ereccval(\sigma)}$
% is equivalent to $x^{\ereccval(\sigma)}
%        y^{\eareccpeak(\sigma)}$.
% Now apply \cite[Theorem~2.21]{Sokal-Zeng_masterpoly}
% specialized to $y_1 = u_1 = v_1 =1$,
% $x_1 = y$:
% we obtain an S-fraction with
% $\alpha_{n} =  (\lambda + n - 1) (y + n - 1)$.
% By contraction this yields a J-fraction with
% $\gamma_0 = \lambda y$,
% $\gamma_n =  (\lambda + 2n - 1) (y + 2n - 1) + (\lambda + 2n) (y + 2n)$ for $n \ge 1$,
% $\beta_n = (\lambda+2n-2)(\lambda+2n-1)(y+2n-2)(y+2n-1)$.

\medskip

$\bm{\lambda=0}$.
For all $n \ge 1$ and all $\sigma \in \Sym_n$,
each permutation contains at least one cycle.
Therefore, setting $\lambda=0$ suppresses all permutations for $n \ge 1$,
and we have $P_n(x,y,0) = \delta_{n0}$ (Kronecker delta).
The value of $\gamma_1$ stated in the table above is completely irrelevant,
because $\beta_1 = 0$.

\medskip

$\bm{\lambda=-1}$.
We have found empirically an S-fraction with
$\alpha_{2k-1} = -xy$, $\alpha_{2k}  = -(1-x)(1-y)$ for $k \ge 1$,
and have verified it through $\alpha_6$.
More generally, we have found empirically an S-fraction
involving the family of polynomials
\reff{def.Qnhat.cycle-alternating.evenodd},
which we recall are
\begin{eqnarray}
   & & \hspace*{-6mm}
	\widehat{Q}_n(\xe,\ye,\ue,\ve,\xo,\yo,\uo,\vo,,\lambda)
	\;=\;
        \nonumber \\[4mm]
   & & \qquad
   \sum_{\sigma \in \Sym^{\rm ca}_{2n}}
   \xe^{\eareccpeakeven(\sigma)}
   \ye^{\ereccvaleven(\sigma)}
   \ue^{\nrcpeakeven(\sigma)}
   \ve^{\nrcvaleven(\sigma)}
          \:\times
       \qquad\qquad
       \nonumber \\[-2mm]
   & & \qquad\qquad\;\,
   \xo^{\eareccpeakodd(\sigma)}
   \yo^{\ereccvalodd(\sigma)}
   \uo^{\nrcpeakodd(\sigma)}
   \vo^{\nrcvalodd(\sigma)}
   \lambda^{\cyc(\sigma)}
   \;.
\end{eqnarray}
Namely, after specializing $\lambda = -1$
we get an S-fraction with coefficients:
\be
   \alpha_{2k-1}=-\xe \yo \,,\qquad
   \alpha_{2k}  = -(\xo - \uo) (\ye-\ve) \qquad \text{for $k \ge 1$}.
\ee
(Note that $\ue$ and $\vo$ do not appear in this formula.)
We have also verified this through $\alpha_6$.
{\bf Note added:}
We have now proven these S-fractions;
they will be published elsewhere \cite{Deb-Sokal_lambda=-1}.

\medskip

$\bm{x=-1}$, $\bm{y=-1}$.
These are polynomials through $\gamma_1$ and $\beta_2$.
However $\gamma_2$ fails to be a polynomial.

\medskip

$\bm{x=-y}$.
These are polynomials through $\beta_1$ and $\gamma_1$. 
However, $\beta_2$ fails to be a polynomial.

\subsection{J-fraction for the polynomials \reff{eq.Pn.xylam}
   restricted to alternating cycles}

We now consider the polynomials
\be
   P_n^\Altcyc(x,y)
   \;=\;
   \sum_{\sigma \in \Altcyc_{2n}}
       x^{\eareccpeak(\sigma)}
       y^{\ereccval(\sigma)}
   \;,
 \label{eq.Pn.altcyc.xy}
\ee
which are simply the coefficient of $\lambda^1$
in the just-studied polynomials \reff{eq.Pn.xylam.bis}.
We again have the symmetry $x \leftrightarrow y$:
$P^\Altcyc_n(x,y) = P^\Altcyc_n(y,x)$.

We now give the first few J-fraction coefficients
for the polynomial sequence
$\left(P_{n+1}^\Altcyc(x,y)/(xy)\right)_{n\geq 0}$:
they are
\be
\gamma_0 \;=\; x+y, \quad \beta_1 \; =\; (x + y) (3 + x + y + x y) 
\ee
followed by
\be
\gamma_1 \;=\;
{\scalebox{0.75}{$\dfrac{
	30 (x + y) + 18 (x^2 + y^2) + 4 (x^3 + y^3) +
 x y (34 + 29 (x + y) + 5 (x^2 + y^2)) +  x^2 y^2 (8 + x + y)
}{(x + y) (3 + x + y + x y)}$}}
 \;.
\ee
It can then be shown that
\begin{itemize}
\item[(a)] $\gamma_1$ is not a polynomial in $x$ (when $y$ is a given fixed complex value) unless  $y \in  \{-1,+1\}$.

\item[(b)] $\gamma_1$ is not a polynomial in $y$ unless $x \in  \{-1,+1\}$.

\end{itemize}

The specializations $x,y\in \{-1,+1\}$ give rise to polynomials as follows:
\medskip

$\bm{x=+1 \text{ or } y=+1}$.
For the case $y=+1$, 
we specialize  Corollary~\ref{cor.second.Sfrac.altcyc}
to $\ye=\yo=\ue=\uo=\ve=\vo=1$, $\xe = \xo = x$,
to obtain an S-fraction with
$\alpha_{n} =   n (x + n)$.
By contraction this yields a J-fraction with
$\gamma_0 = x+1$,
$\gamma_n =   2n (x + 2n) + (2n+1) (x + 2n+1)$ for $n \ge 1$,
$\beta_n = (2n-1)(2n)(x+2n-1)(x+2n)$.

The case $x=+1$ follows by using the symmetry $x \leftrightarrow y$.

\medskip

$\bm{x=-1 \text{ or } y=-1}$.
These are polynomials through $\beta_2$ and $\gamma_1$.
However, $\gamma_2$ fails to be a polynomial.

\end{document}